\documentclass[12pt]{article}
\usepackage{verbatim,amsmath,amsthm,amssymb,graphics,graphicx}
\author{J\"org Brendle and Saka\'e Fuchino}
\title{Coloring ordinals 
	by reals}
\date{March 24, 2007}
\newif\iftesting
\newif\ifextended
\let\Label\label%
\iftesting
\def\label#1{\marginpar{{\renewcommand{\baselinestretch}{0.6}\tiny 
			#1}}\Label{#1}\ignorespaces}%
\fi

\newcounter{frml}[section]
\def\thefrml{{\arabic{section}.\arabic{frml}}}
\def\frmlabel#1{\refstepcounter{frml}{\def\baka{#1}\ifx\baka\empty\else\Label{#1}\fi}%
{\rm({\thefrml})\hfill\hfill\hfill}}
\def\xitem[#1]{\item[\frmlabel{#1}]\mbox{}%
	\iftesting\marginpar{{\renewcommand{%
				\baselinestretch}{0.6}\tiny#1}}\fi\ignorespaces}
\def\xxitem[#1][#2]{\item[(\ref{#1}{\makebox[1.4ex][c]{#2}})]\mbox{}%
	\iftesting\marginpar{{\renewcommand{%
				\baselinestretch}{0.6}\tiny\{#1\}\{#2\}}}\fi\ignorespaces}
\def\xitemof#1{{\rm({\ref{#1}})}}
\def\xitemdof#1{$\mbox{\rm(\ref{#1})}'$}
\def\xxitemof#1#2{(\ref{#1}#2)}
\def\xitemabove{{\rm(\thefrml)}}
\newenvironment{xitemize}{\begin{list}{}{\parsep=0.5\smallskipamount%
			\itemindent=-0.0ex%
			\itemsep=0.5\smallskipamount\leftmargin=8ex\labelwidth=6ex\labelsep=1.4ex}}%
							 {\end{list}}
\newenvironment{xxitemize}{\begin{list}{}{\parsep=0.5\smallskipamount%
			\itemindent=-0.4ex%
			\itemsep=0.5\smallskipamount\itemindent=0pt\leftmargin=2ex%
			\labelwidth=6.4ex\labelsep=1.4ex}}%
							 {\end{list}}
\def\assert#1{\noindent\ \ \makebox[5ex][l]{\rm (#1)}\ignorespaces}
\def\assertof#1{{\rm (#1)}}
\newenvironment{assertion}[1]{\begin{trivlist}
\newbox\assertbox
\dimen255=\textwidth
\setbox\assertbox=\hbox{\hspace*{\parindent}{#1}\hspace{\labelsep}}
\advance\dimen255 by -1\wd\assertbox
\advance\dimen255 by -1ex
\item[]\unhbox\assertbox\hfill
\begin{minipage}[t]{\dimen255}}%
{\end{minipage}\end{trivlist}}
\newcommand{\bbd}[1]{{\mathbb{#1}}}
\renewcommand{\baselinestretch}{1.2}

\newcommand{\ctentenc}{,\mbox{\hspace{0.04ex}{.}{.}{.}\hspace{0.1ex},\,}}

\newcommand{\Ctenten}{{\cdot}{\cdot}{\cdot}}
\newcommand{\xmbox}[1]{ $\relax{\rm #1}\relax$ }
\newcommand{\cardof}[1]{\mathopen{|\,}#1\mathclose{\,|}}
\newcommand{\ellof}[1]{\cardof{#1}}

\newcommand{\setof}[2]{\{#1\,:\,#2\}}
\newcommand{\ssetof}[1]{\{#1\}}

\newcommand{\seqof}[2]{\langle#1\,:\,#2\rangle}
\newcommand{\pairof}[1]{\langle#1\rangle}
\newcommand{\psetof}[1]{{\mathcal P}\/(#1)}
\newcommand{\psof}[1]{{\mathcal P}\/(#1)}

\newcommand{\mapping}[3]{#1:#2\rightarrow #3}
\newcommand{\isom}[3]{#1:#2\stackrel{\cong}{\rightarrow} #3}
\newcommand{\imageof}{{}^{\,{\prime}{\prime}}}
\newcommand{\fnsp}[2]{\mbox{}^{#1^{\mbox{}\!}}#2}
\newcommand{\circleq}{\mathrel{{\leq}%
		\hspace{-0.86ex}{\lower-0.5ex\hbox{$\scriptscriptstyle\circ$}}}}
\newcommand{\restr}{\restriction}
\newcommand{\dotcup}{\mathrel{\dot{\cup}}}

\newcommand{\forces}[2]{\,\|\hspace{-.35ex}\mbox{\sf--}_{\,#1\,}%
\mbox{\rm``}\,#2\,\mbox{\rm''}}
\newcommand{\notforces}[2]{\rlap{\ 
/}\|\hspace{-.35ex}\mbox{\sf--}_{\,#1\,}%
\mbox{\rm``}\,#2\,\mbox{\rm''}}
\newcommand{\utildeT}[1]{%
	\hspace{0.25ex}\hbox to 0pt%
				 {$\mathop{#1}\limits_{\raise0.24ex\hbox{$\scriptstyle\sim$}}$\hss}%
		\relax\phantom{\underline{#1\hspace{0.2ex}}}}
\newcommand{\utildeS}[1]{%
	\hspace{0.25ex}\hbox to 0pt%
				 {$\mathop{\scriptstyle #1}\limits_{%
						 \raise0.44ex\hbox{$\scriptscriptstyle\sim$}}$\hss}%
		\hspace{0.25ex}\relax\phantom{\underline{#1\hspace{0.25ex}}}}
\newcommand{\utildeSS}[1]{%
	\hbox to 0pt{$\mathop{\scriptscriptstyle #1}%
		\limits_{\raise0.44ex\hbox{$\scriptscriptstyle\sim$}}$\hss}%
		\relax\phantom{\underline{#1\hspace{0.25ex}}}}

\newcommand{\reals}{\bbd{R}}
\newcommand{\continuum}{{\mathfrak c}}

\newcommand{\non}{{\sf non}}
\newcommand{\cov}{{\sf cov}}

\newcommand{\meas}{\mbox{\small\sf null}}
\newcommand{\cat}{\mbox{\small\sf meager}}

\newcommand{\injprod}[1]{\mathopen{{(}\!{(}}#1\mathclose{{)}\!{)}}}
\newcommand{\On}{{\rm On}}
\newcommand{\Reg}{{\rm Reg}}
\newcommand{\Card}{{\rm Card}}
\newcommand{\dom}{\mathop{\rm dom}}
\newcommand{\codom}{\mathop{\rm codom}}
\newcommand{\range}{\mathop{\rm rng}}
\newcommand{\supp}{\mathop{\rm supp}}
\newcommand{\cf}{\mathop{\rm cf}}
\newcommand{\Fn}{\mathop{\rm Fn}}

\newcommand{\pfeil}{\mathrel{{\rightarrow\,}\llap{$\rightarrow$}}}

\newcommand{\GCH}{{\rm GCH}}
\newcommand{\CH}{{\rm CH}}
\newcommand{\CPA}{{\rm CPA}}
\newcommand{\ZF}{{\rm ZF}}
\newcommand{\ZFC}{{\rm ZFC}}
\newcommand{\HP}{{\rm HP}}
\newcommand{\IP}{{\rm IP}}
\newcommand{\SEP}{{\rm SEP}}
\newcommand{\WFN}{{\rm WFN}}
\newcommand{\Princ}{\mbox{\sc Princ}}
\newcommand{\Cs}{\mbox{\rm C}^{\rm s}}
\newcommand{\hCs}{{\vphantom{C}\smash{\hat{\mbox{\rm C}}}}^{\rm s}}
\newcommand{\Csast}{{\mbox{${}^*$\rm C}}^{\rm s}}
\newcommand{\Fs}{\mbox{\rm F}^{\rm s}}

\newcommand{\DOset}{{\mathfrak{D}\mathfrak{O}}}
\newcommand{\donum}{{\mathfrak{d}\mathfrak{o}}}
\newcommand{\Es}{\mathfrak{S}}
\newcommand{\boundingno}{{\mathfrak b}}
\newcommand{\dominatingno}{{\mathfrak d}}

\newtheorem{Thm}{{\bf Theorem}}[section]

\newtheorem{Cor}[Thm]{{\bf Corollary}}
\newtheorem{Prop}[Thm]{{\bf Proposition}}
\newtheorem{Lemma}[Thm]{{\bf Lemma}}

\newtheorem{Claim}{{\bf\bf Claim}}[Thm]
\newtheorem{Subclaim}{{\bf\bf Subclaim}}[Claim]
\newtheorem{Problem}{{\bf Problem}}

\newcommand{\Thmof}[1]{{\rm Theorem \ref{#1}}}
\newcommand{\bfThmof}[1]{{\bf Theorem \ref{#1}}}
\newcommand{\Corof}[1]{{\rm Corollary \ref{#1}}}
\newcommand{\Propof}[1]{{\rm Proposition \ref{#1}}}
\newcommand{\Lemmaof}[1]{{\rm Lemma \ref{#1}}}

\newcommand{\Claimof}[1]{{\rm Claim \rm\ref{#1}}}

\newcommand{\sectionof}[1]{{\rm Section \ref{#1}}}

\newcommand{\Thmabove}{{\rm Theorem \number\theThm}}
\newcommand{\Corabove}{{\rm Corollary \number\theThm}}

\newcommand{\Lemmaabove}{{\rm Lemma \number\theThm}}

\newcommand{\Claimabove}{{\rm Claim \rm\number\theClaim}}

\newcommand{\prf}{\noindent{\bf Proof.\ }\ignorespaces}
\newcommand{\prfof}[1]{\noindent{\bf Proof of #1:\,\ }\ignorespaces}
\newcommand{\prfofClaim}{\noindent\raisebox{-.4ex}{\Large $\vdash$\ \ }}
\newsavebox{\qedbox}\sbox{\qedbox}{% QED symbol (by U.Fuchs modified by S.F.)
{\unitlength=0.065mm \begin{picture}(40,60)
\put(0,0){\framebox(30,44)[lb]{}}
\put(30,-6){\rule{6\unitlength}{44\unitlength}}
\put(10,-6.2){\rule{26\unitlength}{6\unitlength}}
\end{picture}}}
\renewcommand{\qed}{\mbox{}\hfill\usebox{\qedbox}}
\newcommand{\smallqed}%
{\mbox{}\smallskip\hfill\raisebox{-.4ex}{\Large $\dashv$}\\}
\newcommand{\qedof}[1]%
{\mbox{} \hspace*{\fill}{\usebox{\qedbox}{\rm~(#1)}}%
\mbox{}}
\newcommand{\Qedof}[1]%
{\mbox{} \hspace*{\fill}{\usebox{\qedbox}%
{\rm~(#1~\number\theThm)}}}
\newcommand{\qedofThm}{\Qedof{Theorem}}
\newcommand{\qedofCor}{\Qedof{Corollary}}
\newcommand{\qedofProp}{\Qedof{Proposition}}
\newcommand{\qedofLemma}{\Qedof{Lemma}}

\newcommand{\qedskip}{\medskip\smallskip}
\newcommand{\qedofClaim}%
{\mbox{}\hfill\raisebox{-.4ex}{\Large $\dashv$ }\nolinebreak%
\mbox{\rm~({\rm Claim}~\rm\number\theClaim)}}
\newcommand{\qedofSubclaim}%
{\mbox{}\hfill\raisebox{-.4ex}{\Large $\dashv$ }\nolinebreak%
\mbox{\rm~({\rm Subclaim}~\rm\number\theSubclaim)}}
\newcommand{\gma}{{\mathfrak a}}

\newcommand{\calA}{{\mathcal A}}

\newcommand{\calD}{{\mathcal D}}
\newcommand{\calF}{{\mathcal F}}

\newcommand{\calH}{{\mathcal H}}

\newcommand{\calL}{{\mathcal L}}

\newcommand{\calX}{{\mathcal X}}
\newcommand{\poB}{\bbd{B}}
\newcommand{\poC}{\bbd{C}}

\newcommand{\poD}{\bbd{D}}
\newcommand{\poP}{\bbd{P}}
\newcommand{\poS}{\bbd{S}}
\newcommand{\bbbone}{{\mathchoice {\rm 1\mskip-4mu l} {\rm 1\mskip-4mu l}
{\rm 1\mskip-4.5mu l} {\rm 1\mskip-5mu l}}}
\newcommand{\poQ}{\bbd{Q}}

\newcommand{\st}{such that}
\newcommand{\wolog}{without loss of generality}
\newcommand{\Wolog}{Without loss of generality}
\newcommand{\wrt}{with respect to}
\newcommand{\po}{poset}
\newcommand{\pos}{posets}
\newcount\minute	% Current minute within the hour
\newcount\hour		% Current hour (24-hour type)
\newcount\hourMins  % Temporary for taking hour modulo 60
\def\now%
{% Displays today's time as ``hh:mm''
  \minute=\time    % Number of minutes since midnight
  \hour=\time \divide \hour by 60 % Get hours
  \hourMins=\hour \multiply\hourMins by 60
  \advance\minute by -\hourMins % Hours modulo 60
  \zeroPadTwo{\the\hour}:\zeroPadTwo{\the\minute}%
}% --- \now ---
\def\zeroPadTwo#1%
{% Left zero pad of the argument to 2 digits.  The argument
  \ifnum #1<10 0\fi    % Conditionally output a zero
  #1%    Then output the argument
}% --- \zeroPadTwo ---

\begin{document}
\maketitle
\footnotetext{The first author is partially supported by Grant-in-Aid for 
	Scientific Research (C) 17540116, Japan Society for the Promotion of 
	Science. The second author is partially supported by Chubu 
	University grant 16IS55A.}
\begin{abstract}
	We study combinatorial principles we call Homogeneity Principle 
	$\HP(\kappa)$ and Injectivity Principle $\IP(\kappa,\lambda)$ for regular
	$\kappa>\aleph_1$ and $\lambda\leq\kappa$ which are 
	formulated in terms of coloring the ordinals $<\kappa$ by reals. 

	These principles are strengthenings of $\Cs(\kappa)$ and $\Fs(\kappa)$ of 
	I.\ Juh\'asz, L.\ Soukup and Z.\ Szentmikl\'ossy \cite{ju-so-sz}. 
	Generalizing, their results, we show e.g.\ that 
	$\IP(\aleph_2,\aleph_1)$ (hence also $\IP(\aleph_2,\aleph_2)$ as well as   
	$\HP(\aleph_2)$) holds in a generic extension of a model of \CH\ by Cohen 
	forcing  
	and $\IP(\aleph_2,\aleph_2)$ (hence also $\HP(\aleph_2)$) holds in 
	a generic extension by 
	countable support side-by-side product of Sacks or Prikry-Silver forcing 
	(\Corof{fs-prod-IP}). We also show that the latter result is optimal 
	(\Thmof{prikry-silver}).  

	Relations between these principles and their 
	influence on the values of the variations
	$\boundingno^\uparrow$, $\boundingno^h$, $\boundingno^*$, $\donum$
	of the bounding 
	number $\boundingno$ are studied. 

	One of the consequences of $\HP(\kappa)$ besides $\Cs(\kappa)$ is that 
	there is no projective well-ordering of length $\kappa$ on any subset of 
	$\fnsp{\omega}{\omega}$. 
	We construct a model in which there is no projective well-ordering of 
	length $\omega_2$ on any subset of $\fnsp{\omega}{\omega}$
	($\donum=\aleph_1$ in our terminology) while $\boundingno^*=\aleph_2$ 
	(\Thmof{new-model}).  
\end{abstract}
\begin{center}
	bounding number, 
	projective well-ordering,\\ Cohen forcing, Brendle-LaBerge forcing, 
	Prikry-Silver forcing
\end{center}

\ifextended
\noindent\mbox{}\hfill
\framebox{\mbox{\,\parbox{0.7\textwidth}{This version of the paper
		contains some details which will be omitted in the version to be 
		published. 

		To obtain the version for publication, comment out the line 
		containing 
		{\tt \textbackslash extendedtrue} in the \LaTeX\ source file and 
		recompile it.}\,}}\hfill\mbox{}\\
\fi
\section{Introduction}
The Cohen model which is obtained by adding at least $\aleph_2$ Cohen
reals over a model of $\GCH$ was the first and simplest model for the
negation of $\CH$, and it is still one of the most important. 
A plethora of statements have been shown to be consistent with
$\ZFC$ by adjoining Cohen reals, and it is therefore natural to look
for axioms which hold in the Cohen model and from which many 
such statements can be decided, that is, axioms which capture as much
as possible of the combinatorial structure of the Cohen extension.
Something similar has been done for the iterated Sacks model
by Ciesielski and Pawlikowski who devised the Covering Property Axiom
$\CPA$ \cite{CP}.

For Cohen models, several such axioms have been proposed 
in the past. Some of them are {\em homogeneity type} statements,
that is, they assert that given at least $\omega_2$ many reals,
many of them ``look similar". Examples are the 
combinatorial principles $\Cs(\kappa)$, $\hCs(\kappa)$, and
$\Fs(\kappa)$ introduced by I. Juh\'asz, L. Soukup and
Z. Szentmikl\'ossy \cite{ju-so-sz} who showed they hold in Cohen models
(see Section 2 below for definitions).

On the other hand, rather different-looking statements have
been also investigated in connection with Cohen models, for example, the 
axiom $\WFN$ asserting that 
$\pairof{\psof{\omega}, {\subseteq}}$ has the weak Freese-Nation Property
(see \cite{fu-ko-sh}, \cite{fu-so} and \cite{fu-ge-so}). Here a partial 
ordering $\pairof{P, \leq}$ has the 
weak Freese-Nation Property if there is a mapping
$\mapping{f}{P}{[P]^{\aleph_0}}$ such that for all $p,q \in P$,
$p \leq_P q$ holds 
if and only if there is an $r \in f(p) \cap f(q)$ such that $p \leq_P r \leq_P q$.

In \cite{fu-ko-sh}, it is shown that $\WFN$ holds in a Cohen model for 
adding $\aleph_n$ Cohen reals for any $n<\omega$. If we start e.g.\ from 
$V=L$ then $\WFN$ holds even after adding any number of Cohen reals 
(\cite{fu-so}). 
In \cite{fu-ge-so}, it was shown that $\WFN$ implies many of the known
combinatorial properties of Cohen models and so it may be seen
as an axiomatization of the combinatorial structure of the
Cohen extension.
Since $\WFN$ can be reformulated in terms
of elementary submodels, $\WFN$ as well as some closely related
statements have come to be known as {\em elementary submodel
type} axioms (see \cite{juhasz-kunen} for this).

At first glance it seemed that there would be no connection between these 
two types of axioms except that they both hold in a Cohen model. 
Surprisingly enough though 
S.\ Shelah \cite{shelah} showed that $\Cs(\aleph_2)$ follows from the 
combinatorial principle he called $\Princ$,  which is a consequence of 
\WFN. The proof can be easily recast 
to show that $\WFN$ implies $\Cs(\kappa)$ for all regular $\kappa>\aleph_1$ 
(see \cite{fuchino-x} for more details). 

In this paper, we introduce some new principles of the homogeneity type,
namely, the Homogeneity Principle $\HP(\kappa)$ and the Injectivity
Principle $\IP (\kappa, \lambda)$ which are formulated in terms of
homogeneity of colorings of ordinals below the cardinal $\kappa$ by reals.
We establish that these axioms hold in Cohen models and address
the question in which other models these axioms hold as well.
It turns out that, in fact, these principles seem to capture a good deal of the
combinatorial features of models of set theory obtained by forcing by the
side-by-side (finite or countable support) product of copies of a fixed
relatively small partial ordering (see \Thmof{forcing-IP*} and \Corof{fs-prod-IP}).

Though the relation of these principles to \WFN\ is not yet completely clear, 
our principles imply the principles of I.\ Juh\'asz, L.\ Soukup and Z.\ 
Szentmikl\'ossy (\Thmof{HP-implies-Cs})
and thus can be seen as natural strengthenings of these 
principles. 

Our paper is organized as follows.

In \sectionof{the-principles}, we review the principles $\Cs(\kappa)$, 
$\hCs(\kappa)$ and $\Fs(\kappa)$ of I.\ Juh\'asz, L.\ Soukup and Z.\ 
Szentmikl\'ossy, and introduce our principles $\HP(\kappa)$ 
and $\IP(\kappa,\lambda)$. Some basic facts in \ZFC\ concerning these 
principles are also proved. In particular, we show 
that $\Cs(\kappa)$ and $\hCs(\kappa)$ follow from $\HP(\kappa)$ 
(\Thmof{HP-implies-Cs}), $\Fs(\kappa)$ follows from $\IP(\kappa,\aleph_1)$ 
(\Thmof{IP-implies-Fs}) and 
$\HP(\kappa)$ follows from $\IP(\kappa,\kappa)$ (\Thmof{IP-implies-HP}). 

After reviewing some cardinal invariants introduced in \cite{fuchino-x} 
which are variants of the bounding number $\boundingno$ and the shrinking 
number $\boundingno^*$ in \cite{eda-kada-yuasa}, we study in 
\sectionof{bounding-no} 
the effect of the combinatorial principles $\Cs(\kappa)$, $\hCs(\kappa)$ 
and $\HP(\kappa)$ on the values of these cardinal invariants. 

In \sectionof{forcing-construction} we give a 
forcing construction of models of $\IP(\kappa,\lambda)$ 
(\Thmof{forcing-IP*}) and its applications (\Corof{fs-prod-IP}). 
The results in this section improve consistency results in \cite{ju-so-sz}.

As a further application of \Thmof{forcing-IP*}
we show in \sectionof{IP-aleph2-aleph1} the consistency of 
$\neg\IP(\aleph_2,\aleph_1)$ and $\IP(\aleph_2,\aleph_2)$. 

One of the consequences of $\HP(\aleph_2)$ discussed in 
\sectionof{bounding-no} is that there is no definable well-ordering of 
length $\omega_2$ on any subset of $\fnsp{\omega}{\omega}$ (or
$\donum=\aleph_1$ in our notation). 
Refining a forcing extension of Brendle and LaBerge 
\cite{brendle-laberge}, we prove in \sectionof{Omega} the consistency of
$\donum=\aleph_1$ with $\boundingno^*=\aleph_2$ (\Thmof{new-model}). 
We also show that the model of $\donum=\aleph_1$ 
and $\boundingno^*=\aleph_2$ we construct in this section satisfies a 
strong negation of $\Cs(\aleph_2)$.  

\sectionof{forcing(A)} is devoted to 
the consistency proof of the combinatorial principle used in the 
proof of \Thmof{new-model}. 

In \sectionof{zusammenfassung-und-offene-frage}, we summarize the 
consistency results obtained in this paper together with other consistency 
results established by 
some previously known constructions. We discuss also some open problems at 
the end of the section.  
%% The forcing constructions also show that the implications of the principles 
%% on the values of the variants of bounding number established in 
%% \sectionof{bounding-no} are optimal. 

\section{Combinatorial principles formulated in terms of coloring of ordinals by reals}
\label{the-principles}
For any set $X$, let
\begin{xitemize}
	\xitem[] $\injprod{X}^n=\setof{\vec{x}\in X^n}{\vec{x}\xmbox{ is injective}}$ 
	and 
\end{xitemize}
\begin{xitemize}
	\xitem[] $\injprod{X}^{<\omega}=\bigcup_{n<\omega}\injprod{X}^{n}$. 
\end{xitemize}
Likewise, for any sets $X_0$\ctentenc  $X_{n-1}$, let 
\begin{xitemize}
	\xitem[] $\injprod{X_0\ctentenc X_{n-1}}=
	\setof{\vec{x}\in X_0\times\cdots\times X_{n-1}}{\vec{x}\xmbox{ is injective}}$. 
\end{xitemize}

For a cardinal $\kappa$, 
the following principle $\Cs(\kappa)$ was introduced by I.\ Juh\'asz, L.\ 
Soukup and  
Z.\ Szentmikl\'ossy in \cite{ju-so-sz}.

\begin{assertion}{$\Cs(\kappa)$: }\it 
For any matrix $\seqof{a_{\alpha,n}}{\alpha\in\kappa,\,n\in\omega}$ of subsets 
of $\omega$ and 
$T\subseteq\fnsp{\omega>}{\omega}$, at least one of the following holds:
\begin{xitemize}
\item[\assertof{c0}] there is a stationary $S\subseteq\kappa$ \st\ 
	$\bigcap_{n<|t|}a_{\alpha_n,t(n)}\not=\emptyset$ for 
	all $t\in T$ and
	$\pairof{\alpha_0\ctentenc \alpha_{|t|-1}}\in\injprod{S}^{<\omega}$;
\item[\assertof{c1}] there exist $t\in T$ and stationary $S_0$\ctentenc 
	$S_{|t|-1}\subseteq \kappa$ \st\
	$\bigcap_{n<|t|}a_{\alpha_n,t(n)}=\emptyset$ for all
	$\pairof{\alpha_0\ctentenc \alpha_{|t|-1}}
	\in\injprod{S_0\ctentenc S_{|t|-1}}$. 
\end{xitemize}
\end{assertion}
For any cardinal $\kappa$ it is easy to see that $\Cs(\kappa)$ holds if and 
only if $\Cs(\cf\kappa)$ holds. Thus it is enough to consider $\Cs(\kappa)$ 
for regular $\kappa$. 
The corresponding assertion is also true for other 
combinatorial principles we are going to introduce in this section. 
Hence, in the rest of this section, we 
shall assume that $\kappa$ is a regular cardinal unless mentioned otherwise. 

The combinatorial principle $\hCs(\kappa)$, a sort of dual 
of the principle $\Cs(\kappa)$, is 
also considered in \cite{ju-so-sz}: 
\begin{assertion}{$\hCs(\kappa)$:}\it
For any $T\subseteq\omega^{<\omega}$ and any matrix 
$\seqof{a_{\alpha,n}}{\alpha<\kappa,\,n\in\omega}$ of subsets of 
$\omega$, at least one of the following holds:
\begin{xitemize}
\item[\assertof{\^c0}] 
	there is a stationary $S\subseteq\kappa$ \st\ 
	$\cardof{\bigcap_{n<|t|}a_{\alpha_n,t(n)}}<\aleph_0$ for every $t\in T$ and
	$\pairof{\alpha_0\ctentenc \alpha_{\cardof{t}-1}}\in\injprod{S}^{|t|}$;
\item[\assertof{\^c1}] there exist $t\in T$ and stationary $S_0$\ctentenc 
	$S_{|t|-1}\subseteq\kappa$ \st\   
	$\cardof{\bigcap_{n<|t|}a_{\alpha_n,t(n)}}=\aleph_0$ for every 
	$\pairof{\alpha_0\ctentenc \alpha_{\cardof{t}-1}}
	\in \injprod{S_0\ctentenc S_{|t|-1}}$.
\end{xitemize}
\end{assertion}

The following is easily seen: 
\begin{Lemma}{\rm (I.\ Juh\'asz, L.\ Soukup and 
		Z.\ Szentmikl\'ossy \cite{ju-so-sz})}\label{Cs-aleph-1}

	\assert{a} Neither of $\Cs(\aleph_1)$ and $\hCs(\aleph_1)$ holds.

	\assert{b} 
	$\Cs(\kappa)$ and $\hCs(\kappa)$ 
	hold for any regular $\kappa>2^{\aleph_0}$.\qed
\end{Lemma}

Let us call 
a subset $A$ of $\calH(\aleph_1)$ definable if there are a formula 
$\varphi$ and $a\in\calH(\aleph_1)$ \st\ 
$A=\setof{x\in\calH(\aleph_1)}{\pairof{\calH(\aleph_1),{\in}}\models\varphi(x,a)}$. 
Note that for any $n\in\omega$, $A\subseteq\reals^n$ is projective 
if and only if it is definable in our sense. Note also since all elements of 
$\calH(\aleph_1)$ can be coded by elements of $\fnsp{\omega}{\omega}$ we 
may assume that $a$ as above is an element of $\fnsp{\omega}{\omega}$. 

In \Thmof{HP-implies-Cs}, we show that the following {\em Homogeneity 
	Principle $\HP(\kappa)$} implies both of $\Cs(\kappa)$ and 
$\hCs(\kappa)$. 
\begin{assertion}{$\HP(\kappa)$: }\it{}For any
	$\mapping{f}{\kappa}{\psof{\omega}}$ and any definable
	$A\subseteq\injprod{\psof{\omega}}^{<\omega}$, at least one of the following holds:
	\begin{xitemize}
		\item[\assertof{h0}] there is a stationary $S\subseteq\kappa$ \st\
			$\injprod{f\imageof S}^{<\omega}\setminus\ssetof{\emptyset}\subseteq A$;
		\item[\assertof{h1}] there are $k\in\omega\setminus 1$ and 
			stationary $S_0$\ctentenc 	$S_{k-1}\subseteq \kappa$ \st\ 
			$\injprod{f\imageof S_0\ctentenc f\imageof S_{k-1}}\cap A=\emptyset$. 
	\end{xitemize}
\end{assertion}

\noindent
Note that $\psof{\omega}$ in the definition of $\HP(\kappa)$ above can be 
replaced by $\reals$, $\fnsp{\omega}{\omega}$,  $(\psof{\omega})^n$ or  
$(\fnsp{\omega}{\omega})^n$ etc.\ since these spaces can be coded 
as definable subsets of $\psof{\omega}$ and vice versa. 

As for $\Cs(\kappa)$ (and $\hCs(\kappa)$), it is enough to consider
$\HP(\kappa)$ for regular $\kappa$. 
\Lemmaof{Cs-aleph-1} is also true for $\HP(\kappa)$:
\begin{Lemma}\label{HP-aleph-1}
	\assert{a} $\HP(\aleph_1)$ does not hold.

\assert{b} $\HP(\kappa)$ holds for any regular $\kappa>2^{\aleph_0}$. 
\end{Lemma}
\prf \assertof{a}: This follows from \Lemmaof{Cs-aleph-1} and 
\Thmof{HP-implies-Cs}. 

\assertof{b}: Let $\kappa>2^{\aleph_0}$ be a regular cardinal. Suppose that 
$\mapping{f}{\kappa}{\psetof{\omega}}$ and $A$ are as in the definition of 
$\HP(\kappa)$. 
Then there is a stationary 
$S\subseteq\kappa$ \st\ $f\restr S$ is constant. If \assertof{h0} in the 
definition of $\HP(\kappa)$ does not hold then we must have 
$\injprod{f\imageof S}^{1}\cap A=\emptyset$ since
$\injprod{f\imageof S}^{n}=\emptyset$ for all $n>1$.  
Hence \assertof{h1} holds with $n=1$ and 
$S_0=S$. \qedofLemma\qedskip

The following combinatorial principle $\Fs(\kappa)$ is also introduced in 
\cite{ju-so-sz}: 

\begin{assertion}{$\Fs (\kappa)$: }\it{}For any $T\subseteq\omega^{<\omega}$ 
	and any matrix  
	$\seqof{a_{\alpha,n}}{\alpha<\kappa,n\in\omega}$ of subsets of 
	$\omega$, at least one of the following holds: 
	\begin{xitemize}
	\item[\assertof{f\,0}] there is a stationary $S\subseteq\kappa$ \st\\
		$\cardof{\setof{{\bigcap}_{n<|t|}a_{\alpha_n,t(n)}}%
		{t\in T\mbox{ and }\pairof{\alpha_0\ctentenc \alpha_{\cardof{t}-1}}
			\in \injprod{S}^{|t|}}}\leq\aleph_0$\,;
	\item[\assertof{f\,1}] there are $t\in T$ and stationary $S_0$\ctentenc  
		$S_{|t|-1}\subseteq\kappa$ \st\ for every 
		$\pairof{\alpha_0\ctentenc \alpha_{\cardof{t}-1}}, 
		\pairof{\beta_0\ctentenc\beta_{\cardof{t}-1}}
		\in \injprod{S_0\ctentenc S_{|t|-1}}$, if $\alpha_n\not=\beta_n$ for all 
		$n<|t|$, then 
		$\bigcap_{n<|t|}a_{\alpha_n,t(n)}\not=\bigcap_{n<|t|}a_{\beta_n,t(n)}$. 
	\end{xitemize}
\end{assertion}

\begin{Lemma}{\rm (I.\ Juh\'asz, L.\ Soukup and 
		Z.\ Szentmikl\'ossy \cite{ju-so-sz})}\label{Fs-aleph-1}
	
	\assert{a} $\Fs(\aleph_1)$ does not hold.

	\assert{b} $\Fs(\kappa)$ holds for every regular $\kappa>2^{\aleph_0}$.

	\assert{c} $\Fs(\kappa)$ implies $\hCs(\kappa)$. \qed
\end{Lemma}

A combinatorial principle in terms of coloring of ordinals by reals 
corresponding naturally to $\Fs(\kappa)$ might be the following {\em 
	Injectivity Principle}  
$\IP(\kappa,\lambda)$ for cardinals $\kappa$ and $\lambda$ with
$\lambda\leq\kappa$: 
\begin{assertion}{$\IP(\kappa,\lambda)$: }
\em For any $\mapping{f}{\kappa}{\psetof{\omega}}$ and definable 
$\mapping{g}{\injprod{\psetof{\omega}}^{<\omega}}{\psetof{\omega}}$, 
at least one of the following holds: 
\begin{xitemize}
\item[\assert{i0}] there is a stationary $S\subseteq\kappa$ \st\ 
	$\cardof{g\imageof\injprod{f\imageof S}^n}<\lambda$ for every $n\in\omega$; 
\item[\assert{i1}] there are $k\in\omega\setminus 1$ and stationary 
	$S_0$\ctentenc $S_{k-1}\subseteq\kappa$ \st\ for any 
	$\pairof{x_0\ctentenc x_{k-1}}$, 
	$\pairof{y_0\ctentenc y_{k-1}}\in
	\injprod{f\imageof S_0\ctentenc f\imageof S_{k-1}}$, if $x_n\not= y_n$ for 
	all $n<k$, then we have $g(x_0\ctentenc x_{k-1})\not=g(y_0\ctentenc y_{k-1})$. 
\end{xitemize}
\end{assertion}
Again here, we may replace $\psof{\omega}$ in the definition 
of $\IP(\kappa,\kappa)$ above by $\reals$, $\fnsp{\omega}{\omega}$,
$(\psof{\omega})^n$ or $(\fnsp{\omega}{\omega})^n$ etc.

\begin{Lemma}\label{IP-aleph-1}
	\assert{a} For $\lambda\leq\lambda'\leq\kappa$, $\IP(\kappa,\lambda)$ 
	implies $\IP(\kappa,\lambda')$.  

	\assert{b} $\IP(\aleph_1,\aleph_1)$ does not hold.
\end{Lemma}
\prf \assertof{a}: Immediate from the definition. 

\assertof{b}: 
By \Lemmaof{HP-aleph-1},\,\assertof{a} and \Thmof{IP-implies-HP}.
\qedofLemma\qedskip

$\IP(\kappa,\aleph_0)$ for a regular cardinal $\kappa$ is equivalent to the 
cardinal inequality $2^{\aleph_0}<\kappa$. 
\begin{Prop}\label{continuum-less-than-kappa}
	For a regular cardinal $\kappa$ the following are equivalent:
\\
	\assertof{a} $\IP(\kappa,\aleph_0)$ holds;\quad
	\assertof{b} $2^{\aleph_0}<\kappa$;\quad
	\assertof{c} $\IP(\kappa,2)$ holds.
\end{Prop}
\prf \assertof{a} $\Rightarrow$ \assertof{b}: Suppose that
$2^{\aleph_0}\geq\kappa$. We show that 
$\IP(\kappa,\aleph_0)$ does not hold. Let 
$\mapping{f}{\kappa}{\psetof{\omega}}$ be any injective mapping and 
$\mapping{g}{\injprod{\psetof{\omega}}^{<\omega}}{\psetof{\omega}}$ be defined by 
$g(\emptyset)=\emptyset$, $g(\pairof{x})=\emptyset$ for all $x\in\psof{\omega}$ 
and 
\begin{xitemize}
\item[] $g(\pairof{x_0\ctentenc x_{n-1}})=\min\setof{m\in\omega}%
	{m\in x_0\not\leftrightarrow m\in x_1}$ 
\end{xitemize}
for $\pairof{x_0\ctentenc x_{n-1}}\in\injprod{\psof{\omega}}^n$ with $n\geq2$.
Let $S$ be any stationary subset of $\kappa$. Then 
$\cardof{g\imageof\injprod{f\imageof S}^2}\geq\aleph_0$: Suppose not and let
$k\in\omega$ be \st\ 
$g\imageof\injprod{f\imageof S}^2\subseteq k$. Since $\psetof{k+1}$ is 
finite, there are  
$\alpha$, $\beta\in S$, $\alpha\not=\beta$ \st\ 
$f(\alpha)\cap (k+1)=f(\beta)\cap (k+1)$. But then, by definition of 
$g$, it follows that $g(\pairof{f(\alpha), f(\beta)})>k$. This is a contradiction. 

Thus \assertof{i0} does not hold for these $f$ and $g$. 
On the other hand, 
for arbitrary stationary subsets $S_0$\ctentenc $S_{n-1}$ of $\kappa$, as 
there are only countably many values of $g$, we can find 
$\pairof{x_0\ctentenc x_{n-1}}$, 
$\pairof{y_0\ctentenc y_{n-1}}
\in\injprod{f\imageof S_0\ctentenc f\imageof S_{n-1}}$  \st\ $x_i\not= y_i$ for 
all $i<n$ and
$g(\pairof{x_0\ctentenc x_{n-1}})=g(\pairof{y_0\ctentenc y_{n-1}})$. Thus 
\assertof{i1} neither holds. \smallskip

\assertof{b} $\Rightarrow$ \assertof{c}: Suppose $2^{\aleph_0}<\kappa$. 
For $\mapping{f}{\kappa}{\psof{\omega}}$ and 
$\mapping{g}{\injprod{\psof{\omega}}^{<\omega}}{\psof{\omega}}$ as in  
the definition of $\IP(\kappa,2)$, there is a stationary 
$S\subseteq\kappa$ \st\ $f$ is constant on $S$. This $S$ witnesses that 
\assertof{i0} holds. \smallskip

\assertof{c} $\Rightarrow$ \assertof{a}: This follows from 
\Lemmaof{IP-aleph-1},\,\assertof{a}. 
\qedofProp
\begin{Cor}
	$\IP(\aleph_2,\aleph_0)$ is equivalent to \CH.\qed
\end{Cor}
$\IP(\aleph_2,\aleph_1)$ and $\IP(\aleph_2,\aleph_2)$ are thus the first 
two non-trivial instances of $\IP(\kappa,\lambda)$. 

For $\kappa\geq\aleph_2$, the principles introduced in this section and some 
other principles  
discussed in \cite{fuchino-x} can be put together in the following 
diagram:
\bigskip\\
\mbox{}\hspace{-1.6cm}
%%\input{coloring-x-fig1}
%WinTpicVersion3.08
\unitlength 0.1in
\begin{picture}( 46.3500, 26.0000)( -3.9000,-28.0000)
% VECTOR 0 0 3 0
% 4 1390 650 900 1700 1640 820 1640 820
% 
\special{pn 20}%
\special{pa 1390 650}%
\special{pa 900 1700}%
\special{fp}%
\special{sh 1}%
\special{pa 900 1700}%
\special{pa 946 1648}%
\special{pa 924 1652}%
\special{pa 910 1632}%
\special{pa 900 1700}%
\special{fp}%
\special{pa 1640 820}%
\special{pa 1640 820}%
\special{fp}%
% VECTOR 0 0 3 0
% 2 1530 650 1780 1040
% 
\special{pn 20}%
\special{pa 1530 650}%
\special{pa 1780 1040}%
\special{fp}%
\special{sh 1}%
\special{pa 1780 1040}%
\special{pa 1762 974}%
\special{pa 1752 996}%
\special{pa 1728 996}%
\special{pa 1780 1040}%
\special{fp}%
% STR 2 0 3 0
% 3 1490 370 1490 470 5 0
% $\IP(\kappa,\aleph_1)$
\put(14.9000,-4.7000){\makebox(0,0){$\IP(\kappa,\aleph_1)$}}%
% STR 2 0 3 0
% 3 880 1790 880 1890 5 0
% $\Fs(\kappa)$
\put(8.8000,-18.9000){\makebox(0,0){$\Fs(\kappa)$}}%
% STR 2 0 3 0
% 3 2285 1585 2285 1685 5 0
% $\HP(\kappa)$
\put(22.8500,-16.8500){\makebox(0,0){$\HP(\kappa)$}}%
% VECTOR 0 0 3 0
% 2 2185 1885 1585 2685
% 
\special{pn 20}%
\special{pa 2186 1886}%
\special{pa 1586 2686}%
\special{fp}%
\special{sh 1}%
\special{pa 1586 2686}%
\special{pa 1642 2644}%
\special{pa 1618 2642}%
\special{pa 1610 2620}%
\special{pa 1586 2686}%
\special{fp}%
% VECTOR 2 0 3 0
% 2 990 2060 1390 2700
% 
\special{pn 8}%
\special{pa 990 2060}%
\special{pa 1390 2700}%
\special{fp}%
\special{sh 1}%
\special{pa 1390 2700}%
\special{pa 1372 2634}%
\special{pa 1362 2656}%
\special{pa 1338 2654}%
\special{pa 1390 2700}%
\special{fp}%
% STR 2 0 3 0
% 3 1485 2785 1485 2885 5 0
% $\hCs(\kappa)$
\put(14.8500,-28.8500){\makebox(0,0){$\hCs(\kappa)$}}%
% VECTOR 0 0 3 0
% 2 2310 1880 2890 2700
% 
\special{pn 20}%
\special{pa 2310 1880}%
\special{pa 2890 2700}%
\special{fp}%
\special{sh 1}%
\special{pa 2890 2700}%
\special{pa 2868 2634}%
\special{pa 2860 2656}%
\special{pa 2836 2658}%
\special{pa 2890 2700}%
\special{fp}%
% STR 2 0 3 0
% 3 2985 2785 2985 2885 5 0
% $\Cs(\kappa)$
\put(29.8500,-28.8500){\makebox(0,0){$\Cs(\kappa)$}}%
% VECTOR 2 0 3 0
% 2 4085 2085 3085 2685
% 
\special{pn 8}%
\special{pa 4086 2086}%
\special{pa 3086 2686}%
\special{fp}%
\special{sh 1}%
\special{pa 3086 2686}%
\special{pa 3152 2668}%
\special{pa 3132 2658}%
\special{pa 3132 2634}%
\special{pa 3086 2686}%
\special{fp}%
% STR 2 0 3 0
% 3 4285 1865 4285 1965 5 0
% $\Princ(\kappa,\kappa)$
\put(42.8500,-19.6500){\makebox(0,0){$\Princ(\kappa,\kappa)$}}%
% VECTOR 2 0 3 0
% 2 4235 1385 4235 1785
% 
\special{pn 8}%
\special{pa 4236 1386}%
\special{pa 4236 1786}%
\special{fp}%
\special{sh 1}%
\special{pa 4236 1786}%
\special{pa 4256 1718}%
\special{pa 4236 1732}%
\special{pa 4216 1718}%
\special{pa 4236 1786}%
\special{fp}%
% STR 2 0 3 0
% 3 4235 1085 4235 1185 5 0
% $\SEP(\kappa,\kappa)$
\put(42.3500,-11.8500){\makebox(0,0){$\SEP(\kappa,\kappa)$}}%
% STR 2 0 3 0
% 3 100 100 100 200 0 0
% 
\put(1.0000,-2.0000){\makebox(0,0)[lb]{}}%
% STR 2 0 3 0
% 3 780 1140 780 1240 5 0
% \tiny\Thmof{IP-implies-Fs}
\put(7.8000,-12.4000){\makebox(0,0){\tiny\Thmof{IP-implies-Fs}}}%
% STR 2 0 3 0
% 3 2490 1340 2490 1440 5 0
% \tiny\Thmof{IP-implies-HP}
\put(24.9000,-14.4000){\makebox(0,0){\tiny\Thmof{IP-implies-HP}}}%
% STR 2 0 3 0
% 3 2160 2300 2160 2400 5 0
% \tiny\Thmof{HP-implies-Cs}
\put(21.6000,-24.0000){\makebox(0,0){\tiny\Thmof{HP-implies-Cs}}}%
% STR 2 0 3 0
% 3 1020 2360 1020 2460 5 0
% \tiny\cite{ju-so-sz}
\put(10.2000,-24.6000){\makebox(0,0){\tiny\cite{ju-so-sz}}}%
% STR 2 0 3 0
% 3 3690 2330 3690 2430 5 0
% \tiny\cite{shelah} (see also \cite{fuchino-x})
\put(36.9000,-24.3000){\makebox(0,0){\tiny\cite{shelah} (see also \cite{fuchino-x})}}%
% STR 2 0 3 0
% 3 4430 1470 4430 1570 5 0
% \tiny\cite{fuchino-x}
\put(44.3000,-15.7000){\makebox(0,0){\tiny\cite{fuchino-x}}}%
% STR 2 0 3 0
% 3 1900 1110 1900 1210 5 0
% $\IP(\kappa,\kappa)$
\put(19.0000,-12.1000){\makebox(0,0){$\IP(\kappa,\kappa)$}}%
% VECTOR 0 0 3 0
% 2 1970 1330 2130 1570
% 
\special{pn 20}%
\special{pa 1970 1330}%
\special{pa 2130 1570}%
\special{fp}%
\special{sh 1}%
\special{pa 2130 1570}%
\special{pa 2110 1504}%
\special{pa 2100 1526}%
\special{pa 2076 1526}%
\special{pa 2130 1570}%
\special{fp}%
% VECTOR 2 0 3 0
% 2 3900 2100 1830 2700
% 
\special{pn 8}%
\special{pa 3900 2100}%
\special{pa 1830 2700}%
\special{fp}%
\special{sh 1}%
\special{pa 1830 2700}%
\special{pa 1900 2702}%
\special{pa 1882 2686}%
\special{pa 1888 2662}%
\special{pa 1830 2700}%
\special{fp}%
% STR 2 0 3 0
% 3 2800 2000 2800 2100 5 0
% \tiny\Thmof{HP-implies-Cs}
\put(28.0000,-21.0000){\makebox(0,0){\tiny\Thmof{HP-implies-Cs}}}%
% VECTOR 2 0 3 0
% 2 4240 620 4240 1020
% 
\special{pn 8}%
\special{pa 4240 620}%
\special{pa 4240 1020}%
\special{fp}%
\special{sh 1}%
\special{pa 4240 1020}%
\special{pa 4260 954}%
\special{pa 4240 968}%
\special{pa 4220 954}%
\special{pa 4240 1020}%
\special{fp}%
% STR 2 0 3 0
% 3 4240 370 4240 470 5 0
% $\WFN$
\put(42.4000,-4.7000){\makebox(0,0){$\WFN$}}%
% STR 2 0 3 0
% 3 4440 690 4440 790 5 0
% \tiny\cite{fuchino-geschke}
\put(44.4000,-7.9000){\makebox(0,0){\tiny\cite{fuchino-geschke}}}%
\end{picture}%
\nopagebreak\bigskip\\
\mbox{}\hfill\hfill {\bf fig.\ 1}\hfill\mbox{}
\bigskip

%% By \Lemmaof{Cs-aleph-1}\,\assertof{a}, the diagram in fig.1 makes sense 
%% only for regular $\kappa\geq\aleph_2$. Also, by 
%% \Propof{continuum-less-than-kappa} and by the definition of
%% $\SEP(\kappa,\kappa)$ all argumented principles in the diagram hold 
%% for regular $\kappa>2^{\aleph_0}$. 
In the rest of the section,  we shall prove the implications indicated by 
the thick arrows in fig.1. 
\begin{Thm}\label{HP-implies-Cs}
	For a regular cardinal $\kappa$, $\HP(\kappa)$ implies both $\Cs(\kappa)$ 
	and $\hCs(\kappa)$.
\end{Thm}
\prf \ifextended
Let us prove first that $\HP(\kappa)$ implies $\Cs(\kappa)$.
\else
We prove that $\HP(\kappa)$ implies $\Cs(\kappa)$. The other implication can 
be proved similarly. 
\fi

By \Lemmaof{Cs-aleph-1},\,\assertof{b}, we may assume that
$\kappa\leq 2^{\aleph_0}$.  
Let $\seqof{t_i}{i\in\omega}$ be an enumeration of $\fnsp{\omega>}{\omega}$ 
\st\ 
\begin{xitemize}
	\xitem[enum-T] $\cardof{t_i}\leq i$ for all $i<\omega$
\end{xitemize}
and let 
$\mapping{\iota}{\psof{\omega}}{\psof{\omega}^\omega}$ be a definable 
bijection. For each $x\in\psof{\omega}$ and $i<\omega$, let $(x)_i$ denote 
the i'th component of $\iota(x)$.  

Suppose that $T\subseteq\fnsp{\omega>}{\omega}$ and 
$\calA=\seqof{a_{\alpha,n}}{\alpha<\kappa,n\in\omega}$ is a matrix of subsets of
$\omega$. 
We show that either \assertof{c0} or \assertof{c1} holds for these $\calA$ 
and $T$. 

Let $\mapping{g}{\kappa}{\psof{\omega}}$ be a fixed 
injective mapping which exists by $\kappa\leq 2^{\aleph_0}$. 
Let $\mapping{f}{\kappa}{\psof{\omega}}$ be defined by
\begin{xitemize}
	\xitem[]
	$f(\alpha)=\iota^{-1}(\seqof{a'_{\alpha,n}}{n\in\omega})$ 
\end{xitemize}
where 
\begin{xitemize}
	\xitem[A-dashed] 
	 $a'_{\alpha,n}={}\left\{\,%
	\begin{array}[c]{@{}ll}
		g(\alpha), &\mbox{if }n=0,\\
		a_{\alpha,n-1}, &\mbox{otherwise}.
	\end{array}\right.
	$
\end{xitemize}
Note that $f$ is injective by ``if $n=0$'' clause of \xitemof{A-dashed}. 
For $i<\omega$, let  
\begin{xitemize}
	\xitem[A-star-i]
	$A^*_i={}\left\{\,%
	\begin{array}[c]{@{}ll}
		\setof{\pairof{x_0\ctentenc x_{i-1}}\in\injprod{\psof{\omega}}^i}{
		\bigcap_{n<\cardof{t_i}}(x_n)_{t_i(n)+1}\not=\emptyset}, &\mbox{if }
		t_i\in T,\\
		\injprod{\psof{\omega}}^i, &\mbox{otherwise} 
	\end{array}\right.
	$
\end{xitemize}
and 
\begin{xitemize}
	\xitem[] $A=\bigcup_{i<\omega}A^*_i$. 
\end{xitemize}
It is easy to see that 
$A$ is definable noting that $T\in\calH(\aleph_1)$ and hence $T$ 
can be used as a parameter in the definition of $A$.  By $\HP(\kappa)$, 
we have either \assertof{h0} or \assertof{h1} for these $A$ and $f$. 

Assume first that 
\assertof{h0} holds. Then there is a stationary 
$S\subseteq\kappa$ \st\
$\injprod{f\imageof S}^{<\omega}\setminus\ssetof{\emptyset}\subseteq A$.  
We show that this $S$ witnesses \assertof{c0} for $T$ and $\calA$: for
$t\in T$, let $i\in\omega$ be \st\ $t=t_i$.  
By \xitemof{enum-T}, we have $\cardof{t}\leq i$. For 
$s\in\injprod{S}^{\cardof{t}}$, let $s'\in\injprod{S}^i$ be 
an end-extension 
of $s$. Then 
$\pairof{f(s'(0))\ctentenc f(s'(i-1))}\in\injprod{f\imageof S}^i$ 
since $f$ is injective. Hence
$\pairof{f(s'(0))\ctentenc f(s'(i-1))}\in A^*_i$ by the assumption on $S$. By 
\xitemof{A-star-i}, we have
\begin{xitemize}
\item[] $\emptyset\not=\bigcap_{n<\cardof{t_i}}(f(s'(n)))_{t_i(n)+1}
	=\bigcap_{n<\cardof{t_i}}a'_{s'(n),t_i(n)+1}=\bigcap_{n<\cardof{t}}a_{s(n),t(n)}$.
\end{xitemize}
Thus $T$ and $\calA$ satisfy \assertof{c0}.

Assume now that \assertof{h1} holds. In this case, there are $i\in\omega$ 
and stationary $S_0$\ctentenc  $S_{i-1}\subseteq\kappa$ \st\ 
\begin{xitemize}
\xitem[A-0] $\injprod{f\imageof S_0\ctentenc f\imageof S_{i-1}}
\subseteq\injprod{\psof{\omega}}^i\setminus A^*_i$. 
\end{xitemize}
Let $t=t_i$. Then $t\in T$ by \xitemof{A-0} and ``otherwise'' clause of 
\xitemof{A-star-i}. For 
$s\in\injprod{S_0\ctentenc S_{|t|-1}}$, let
$s'\in\injprod{S_0\ctentenc S_{i-1}}$ be an end extension of $s$. Then we 
have 
$\pairof{f(s'(0))\ctentenc f(s'(i-1))}\in\injprod{f\imageof S_0\ctentenc  
f\imageof S_{i-1}}$. It follows that 
\begin{xitemize}
\item[] 
	$\pairof{f(s'(0))\ctentenc f(s'(i-1))}\in\injprod{\psof{\omega}}^i\setminus A^*_i$ 
\end{xitemize}
by \xitemof{A-0}. 
Hence, by \xitemof{A-star-i}, we have 
\begin{xitemize}
\item[] $\textstyle\emptyset=\bigcap_{n<\cardof{t}}(f(s'(n)))_{t(n)+1}
	=\bigcap_{n<\cardof{t_i}}a'_{s'(n),t(n)+1}
	=\bigcap_{n<\cardof{t}}a_{s(n),t(n)}$.
\end{xitemize}
Thus, $T$ and $\calA$ satisfy \assertof{c1} in this case.

The proof of $\hCs(\kappa)$ from $\HP(\kappa)$ is exactly like the proof 
above with \xitemof{A-star-i} replaced by 
\begin{xitemize}
	\item[\xitemdof{A-star-i}] 
	$A^*_i={}\left\{\,%
	\begin{array}[c]{@{}ll}
		\setof{\pairof{x_0\ctentenc x_{i-1}}\in\injprod{\psof{\omega}}^i}{
		\cardof{\bigcap_{n<\cardof{t_i}}(x_n)_{t_i(n)+1}}<\aleph_0}, &\mbox{if }
		t_i\in T,\\
		\injprod{\psof{\omega}}^i, &\mbox{otherwise.} 
	\end{array}\right.
	$
\end{xitemize}
\qedofThm
\qedskip
\\
$\HP(\kappa)$ also imply other variants of $\Cs(\kappa)$. For example, let
\begin{assertion}{$\Csast(\kappa)$: }\it 
For any matrix $\seqof{a_{\alpha,n}}{\alpha\in\kappa,\,n\in\omega}$ of subsets 
of $\omega$ and 
$T\subseteq\fnsp{\omega>}{\omega}$, at least one of the following holds:
\begin{xitemize}
\item[\assertof{${}^*$c0}] there is a stationary $S\subseteq\kappa$ \st\ 
	$\bigcap_{n<|t|}a_{\alpha_n,t(n)}$ is infinite for 
	all $t\in T$ and
	$\pairof{\alpha_0\ctentenc \alpha_{|t|-1}}\in\injprod{S}^{<\omega}$;
\item[\assertof{${}^*$c1}] there exist $t\in T$ and stationary $S_0$\ctentenc 
	$S_{|t|-1}\subseteq \kappa$ \st\
	$\bigcap_{n<|t|}a_{\alpha_n,t(n)}$ is finite for all
	$\pairof{\alpha_0\ctentenc \alpha_{|t|-1}}
	\in\injprod{S_0\ctentenc S_{|t|-1}}$. 
\end{xitemize}
\end{assertion}
It is easy to see by a proof similar to that of \Thmabove\ 
that $\HP(\kappa)$ implies $\Csast(\kappa)$ as well. 

The following can also be proved 
similarly to \Thmabove.
\begin{Thm}\label{IP-implies-Fs}
	$\IP(\kappa,\aleph_1)$ implies $\Fs(\kappa)$. \ifextended\else\qed\fi
\end{Thm}
\ifextended
\noindent
{\bf [\,The following proof will be omitted in the version for 
		publication.]} 
\smallskip\\
{\footnotesize
\prf 
Since $\Fs(\kappa)$ holds for all $\kappa>2^{\aleph_0}$ (see \cite{ju-so-sz}), 
we may assume that 
$\kappa\leq 2^{\aleph_0}$.  
Let $\seqof{t_i}{i<\omega}$,  
$\mapping{\iota}{\psof{\omega}}{\psof{\omega}^\omega}$ and $(x)_i$ for 
$x\in\psof{\omega}$ and $i<\omega$ be as in the proof of \Thmof{HP-implies-Cs}. 

Suppose that $T\subseteq\fnsp{\omega>}{\omega}$ and 
$\calA=\seqof{a_{\alpha,n}}{\alpha<\kappa,n\in\omega}$ is a matrix of subsets of
$\omega$. 
We show that either \assertof{f\,0} or \assertof{f\,1} holds for these $\calA$ 
and $T$. 

Let $\mapping{h}{\kappa}{\psof{\omega}}$ be a fixed 
injective mapping which exists by $\kappa\leq 2^{\aleph_0}$. 
Let $\mapping{f}{\kappa}{\psof{\omega}}$ be defined by
\begin{xitemize}
	\xitem[]
	$f(\alpha)=\iota^{-1}(\seqof{a'_{\alpha,n}}{n\in\omega})$ 
\end{xitemize}
where 
\begin{xitemize}
	\xitem[A-dashed++] 
	 $a'_{\alpha,n}={}\left\{\,%
	\begin{array}[c]{@{}ll}
		h(\alpha), &\mbox{if }n=0,\\
		a_{\alpha,n-1}, &\mbox{otherwise}.
	\end{array}\right.
	$
\end{xitemize}
$f$ is injective by ``if $n=0$'' clause of \xitemof{A-dashed++}. 
For $i<\omega$, let
$\mapping{g^*_i}{\injprod{\psof{\omega}}^i}{\psof{\omega}}$ be defined by
\begin{xitemize}
	\xitem[A-star-i++]
	$g^*_i(\pairof{x_0\ctentenc x_{i-1}})={}\left\{\,%
	\begin{array}[c]{@{}ll}
		\bigcap_{n<\cardof{t_i}}(x_n)_{t_i(n)+1}, &\mbox{if }
		t_i\in T,\\
		\emptyset, &\mbox{otherwise} 
	\end{array}\right.
	$
\end{xitemize}
for all $\pairof{x_0\ctentenc x_{i-1}}\in\injprod{\psof{\omega}}^i$ 
and 
\begin{xitemize}
	\xitem[] $g=\bigcup_{i<\omega}g^*_i$. 
\end{xitemize}
Then $g$ is definable.  By $\IP(\kappa,\aleph_1)$, 
we have either \assertof{i0} or \assertof{i1} for 
$\lambda=\aleph_1$, for these $f$ and $g$.  

Assume first that 
\assertof{i0} holds. Then there is a stationary 
$S\subseteq\kappa$ \st\ 
$\cardof{X}\leq\aleph_0$ for  
\begin{xitemize}
\item[] $X=g\imageof \injprod{f\imageof S}^{<\omega}$.
\end{xitemize}
We show that $S$ witnesses \assertof{f\,0} for $T$ and $\calA$: for $t\in T$ and
$s\in\injprod{S}^{\cardof{t}}$, let $i\in\omega$ be \st\ $t=t_i$. 
By \xitemof{enum-T}, we have $\cardof{t}\leq i$. Let $s'\in\injprod{S}^i$ be 
an end-extension 
of $s$. $\pairof{f(s'(0))\ctentenc f(s'(i-1))}\in\injprod{f\imageof S}^i$ since $f$ is 
injective. Hence $g(\pairof{f(s'(0))\ctentenc f(s'(i-1))})\in X$. By 
\xitemof{A-star-i++}, we have
\begin{xitemize}
	\item[] $X\ni g(\pairof{f(s'(0))\ctentenc f(s'(i-1))})
		=g^*_i(\pairof{f(s'(0))\ctentenc f(s'(i-1))})$\\[\jot]
		$\phantom{X\ni}=\bigcap_{n<\cardof{t_i}}(f(s'(n)))_{t_i(n)+1}
		=\bigcap_{n<\cardof{t_i}}a'_{s'(n),t_i(n)+1}
		=\bigcap_{n<\cardof{t}}a_{s(n),t(n)}$.
\end{xitemize}
Thus
$\setof{\bigcap_{n<\cardof{t}}a_{s(n),t(n)}}{s\in\injprod{f\imageof S}^{|t|}}
\subseteq X$. 
Since $T$ is countable it follows that 
$\cardof{\setof{{\bigcap}_{n<|t|}a_{s(n),t(n)}}%
			{t\in T\mbox{ and }s\in \injprod{S}^{|t|}}}\leq\aleph_0$. 

Assume now that \assertof{i1} holds. In this case, there are $i\in\omega$ 
and stationary $S_0$\ctentenc  $S_{i-1}\subseteq\kappa$ \st\ 
\begin{xitemize}
	\xitem[A-0++] $g^*_i(x_0\ctentenc x_{i-1})=g(x_0\ctentenc x_{i-1})
	\not=g(y_0\ctentenc y_{i-1})=g^*_i(y_0\ctentenc y_{i-1})$
\end{xitemize}
for all 	$\pairof{x_0\ctentenc x_{i-1}}$, 
	$\pairof{y_0\ctentenc y_{i-1}}\in
	\injprod{f\imageof S_0\ctentenc f\imageof S_{i-1}}$ with $x_n\not= y_n$ for 
	all $n<i$. 
Let $t=t_i$. Then $|t|\leq i$ by \xitemof{enum-T}. 
Also we have $t\in T$ by \xitemof{A-0++} and ``otherwise'' clause of 
\xitemof{A-star-i++}. For 
$s_0$, $s_1\in\injprod{S_0\ctentenc S_{|t|-1}}$, with 
$s_0(n)\not= s_1(n)$ for 
all $n<|t|$, let $s'_0$, $s'_1\in \injprod{S_0\ctentenc S_{i-1}}$ be end 
extensions of $s_0$ and $s_1$ respectively \st\ $s'_0(n)\not= s'_1(n)$ for 
all $n<i$. Then we have 
\begin{xitemize}
\item[] $\bigcap_{n<|t|}a_{s_0(n),t(n)}=g^*_i(s'_0)
	\not=g^*_i(s'_1)=\bigcap_{n<|t|}a_{s_1(n),t(n)}$
\end{xitemize}
by \xitemof{A-0++}. Thus \assertof{f\,1} holds. 
\qedofThm
}\fi

\begin{Thm}\label{IP-implies-HP}
	$\IP(\kappa,\kappa)$ implies $\HP(\kappa)$. 
\end{Thm}
\prf Suppose that $A\subseteq\injprod{\psetof{\omega}}^{<\omega}$ is 
definable and $\mapping{f}{\kappa}{\psetof{\omega}}$. If 
$f^{-1}[\ssetof{x}]$ is stationary for some $x\in\psetof{\omega}$, then,   
either \assertof{h0} holds for $S=f^{-1}[\ssetof{x}]$ or \assertof{h1} holds 
for $n=1$ and $S_0=f^{-1}[\ssetof{x}]$ depending 
on whether $x\in A$ or not. 
Otherwise let 
$\mapping{g}{\injprod{\psetof{\omega}}^{<\omega}}{\psetof{\omega}}$ be defined 
by 
\begin{xitemize}
	\xitem[] $g(\emptyset)=\emptyset$;
	\xitem[A-1] $g(\pairof{x})=\emptyset$ for all $x\in\psetof{\omega}$ and
	\xitem[A-2] $g(\pairof{x_0\ctentenc x_{n-1},x})=\,\left\{
	\begin{array}[c]{@{}ll}
		\emptyset,\qquad &\mbox{if }\pairof{x_0\ctentenc x_{n-1}}\in A,\\[\jot]
		x &\mbox{otherwise}
	\end{array}
	\right.$
\end{xitemize}
for all $\pairof{x_0\ctentenc x_{n-1},x}\in\injprod{\psetof{\omega}}^{n+1}$. 
If \assertof{i0} holds for this $g$ with $S$ as in \assertof{i0}, then, 
by \xitemof{A-2}, we must have 
$g\imageof\injprod{f\imageof S}^{<\omega}=\ssetof{\emptyset}$. Hence 
$\injprod{f\imageof S}^{<\omega}\setminus\ssetof{\emptyset}\subseteq A$. On 
the other hand,  
if \assertof{i1} holds for some $n<\omega$ and $S_0$\,\ctentenc $S_{n-1}$, 
then we should have $n\geq 2$ by \xitemof{A-1} 
and 
$g(\pairof{x_0\,\ctentenc x_{n-2},x_{n-1}})=x_{n-1}$ for all
$x_i\in f\imageof S_i$, $i<n$ by \xitemof{A-2}. It follows that 
$\injprod{f\imageof S_0\ctentenc f\imageof S_{n-2}}
\subseteq \injprod{\psetof{\omega}}^{n-1}\setminus A$ 
by \xitemof{A-2}. 
\qedofThm

\section{The bounding number and its variations}
\label{bounding-no}
In this section, we show that the combinatorial principles introduced in the 
last section make some of the cardinal invariants from \cite{fuchino-x} 
small. %% This will be used in the next section to separate the principles in 
%% some generic extensions. Let us begin with reviewing the definition and the 
%% main facts about the cardinal invariants.  

Adopting the notation of \cite{fuchino-x}, we consider the following 
spectra of cardinal numbers in connection with a partial ordering
$\pairof{P,{\leq}}$; unbounded spectrum, 
hereditary unbounded spectrum and the spectrum of length of $P$\,:
\[
\begin{array}{r@{{}={}}l}
	\Es(P) &\setof{\cardof{X}}{X\subseteq P,X\mbox{ is unbounded in }P,
	\forall B\in[X]^{{<}\cardof{X}}(B\mbox{ is bounded in }P)}\,,
	\\[\jot]
	\Es^h(P) &\setof{\cardof{X}}{X\subseteq P,\,
	\forall B\subseteq X\ (B\mbox{ is bounded in }P\leftrightarrow
	\cardof{B}<\cardof{X})}\,,\\[\jot]
	\Es^\uparrow(P) &\setof{\cf(C)}{C\subseteq P,\,C\mbox{ is an unbounded chain}}\,.
\end{array}
\]\noindent
Clearly, we have 
\begin{xitemize}
	\xitem[Es-0] $\Es^\uparrow(P)\subseteq\Es^h(P)\subseteq\Es(P). $
\end{xitemize}
For $P=\pairof{\fnsp{\omega}{\omega},{\leq^*}}$, we shall simply 
write $\Es^\uparrow$, $\Es^h$ and $\Es$ in place of
$\Es^\uparrow(\pairof{\fnsp{\omega}{\omega},{\leq^*}})$,
$\Es^h(\pairof{\fnsp{\omega}{\omega},{\leq^*}})$ and 
$\Es(\pairof{\fnsp{\omega}{\omega},{\leq^*}})$, respectively. 

Recall that the bounding number $\boundingno$ is defined by 
\begin{xitemize}
	\item[] $\boundingno =\min\setof{\cardof{X}}{X\subseteq\fnsp{\omega}{\omega}
		\mbox{ is unbounded with respect to }{\leq^*}}$. 
\end{xitemize}
The variant $\boundingno^*$ of $\boundingno$ was introduced and studied in 
\cite{eda-kada-yuasa} and \cite{kada-yuasa} where    
\begin{xitemize}
	\item[] $\boundingno^* =\min\setof{\kappa}{\forall X\subseteq\fnsp{\omega}{\omega}\,
	\left(X\mbox{ is unbounded }\rightarrow 
		\exists X'\in[X]^{\leq\kappa}(X'\mbox{ is unbounded})\right)}$. 
\end{xitemize}
$\boundingno$ and $\boundingno^*$ can be characterized in terms of $\Es^\uparrow$, 
$\Es^h$ and $\Es$ as follows:
\begin{Lemma}\label{bounding-0}
	\assert{a} $\boundingno=\min\Es^\uparrow=\min\Es^h=\min\Es$.

\assert{b} $\boundingno^*=\sup\Es$.\ifextended\else\qed\fi
\end{Lemma}
\ifextended
\noindent
{\bf [\,The following proof will be omitted in the version for 
		publication.]} 
\smallskip\\
{\footnotesize
\prf 
\assertof{a}: $\boundingno\leq\min\Es\leq\min\Es^h\leq\min\Es^\uparrow$ is 
clear by definition. It is also  
easily seen that there is an increasing sequence in 
$\pairof{\fnsp{\omega}{\omega},{\leq^*}}$ of length $\boundingno$. Hence
$\min\Es^\uparrow\leq\boundingno$.

\assertof{b}: For any cardinal $\kappa$, we have
\[
\begin{array}{r@{\ {}\Leftrightarrow{}\ }l}
	\kappa<\boundingno^* 
	&\exists X\subseteq\fnsp{\omega}{\omega}\ (X\mbox{ is unbounded }\land\ 
	\forall X'\in[X]^{\leq\kappa}\ (X'\mbox{ is bounded}))\\[\jot]
	&\exists \lambda>\kappa\ (\lambda\in\Es)\\[\jot]
	&\kappa<\sup\Es.
\end{array}
\]\noindent
\mbox{}\vspace{-2\baselineskip}\\
\qedofLemma
\qedskip
}\fi

In analogy to \Lemmaabove,\,\assertof{b}, let 
\begin{xitemize}
	\xitem[] 
	$\boundingno^\uparrow =\sup\Es^\uparrow,\qquad \boundingno^h =\sup\Es^h$.
\end{xitemize}
Recall also that the dominating number $\dominatingno$ is defined as:
\begin{xitemize}
\item[] 
	$\dominatingno=\min\setof{\cardof{X}}{X\subseteq\fnsp{\omega}{\omega},\,
	X\mbox{ dominates }\fnsp{\omega}{\omega}}$. 
\end{xitemize}

By \xitemof{Es-0} and \Lemmaof{bounding-0}, we have 
\begin{Lemma}\label{b-b*}
	 $\boundingno\leq\boundingno^\uparrow\leq\boundingno^h\leq\boundingno^*
	\leq\dominatingno$.\qed
\bigskip
\end{Lemma}

\noindent
\mbox{}\hfill\hspace{-1.6cm}
%%\input{coloring-x-fig2}
%WinTpicVersion3.08
\unitlength 1cm
\begin{picture}(  8.9225,  3.7775)( -2.1675, -3.6650)
% LINE 2 2 3 0
% 2 270 1466 2702 1466
% 
\special{pn 8}%
\special{pa 266 1443}%
\special{pa 2660 1443}%
\special{dt 0.045}%
% ELLIPSE 2 0 3 0
% 4 1416 966 1646 1459 1646 1459 1646 1459
% 
\special{pn 8}%
\special{ar 1394 951 227 486  0.0000000 6.2831853}%
% ELLIPSE 2 0 3 0
% 4 738 1075 930 1459 930 1459 930 1459
% 
\special{pn 8}%
\special{ar 727 1059 189 378  0.0000000 6.2831853}%
% ELLIPSE 2 0 3 0
% 4 2190 858 2459 1466 2446 1466 2446 1466
% 
\special{pn 8}%
\special{ar 2156 845 265 599  0.0000000 6.2831853}%
% STR 2 0 3 0
% 3 738 1018 738 1082 5 0
% $\Es^\uparrow$
\put(1.8450,-2.7050){\makebox(0,0){$\Es^\uparrow$}}%
% STR 2 0 3 0
% 3 1410 909 1410 973 5 0
% $\Es^h$
\put(3.5250,-2.4325){\makebox(0,0){$\Es^h$}}%
% STR 2 0 3 0
% 3 2190 806 2190 870 5 0
% $\Es$
\put(5.4750,-2.1750){\makebox(0,0){$\Es$}}%
% STR 2 0 3 0
% 3 104 1395 104 1459 5 0
% $\boundingno$
\put(0.2600,-3.6475){\makebox(0,0){$\boundingno$}}%
% LINE 2 2 3 0
% 2 245 691 1051 691
% 
\special{pn 8}%
\special{pa 242 681}%
\special{pa 1035 681}%
\special{dt 0.045}%
% LINE 2 2 3 0
% 2 245 467 1781 467
% 
\special{pn 8}%
\special{pa 242 460}%
\special{pa 1753 460}%
\special{dt 0.045}%
% LINE 2 2 3 0
% 2 251 250 2696 250
% 
\special{pn 8}%
\special{pa 248 247}%
\special{pa 2654 247}%
\special{dt 0.045}%
% STR 2 0 3 0
% 3 123 634 123 698 5 0
% $\boundingno^\uparrow$
\put(0.3075,-1.7450){\makebox(0,0){$\boundingno^\uparrow$}}%
% STR 2 0 3 0
% 3 123 397 123 461 5 0
% $\boundingno^h$
\put(0.3075,-1.1525){\makebox(0,0){$\boundingno^h$}}%
% STR 2 0 3 0
% 3 123 186 123 250 5 0
% $\boundingno^*$
\put(0.3075,-0.6250){\makebox(0,0){$\boundingno^*$}}%
% LINE 2 2 3 0
% 2 250 40 2695 40
% 
\special{pn 8}%
\special{pa 247 40}%
\special{pa 2653 40}%
\special{dt 0.045}%
% STR 2 0 3 0
% 3 120 -24 120 40 5 0
% $\dominatingno$
\put(0.3000,-0.1000){\makebox(0,0){$\dominatingno$}}%
\end{picture}%
\hfill\nopagebreak\bigskip\\
\mbox{}\hfill {\bf fig.\ 2}\hfill\mbox{}

Let 
\begin{xitemize}
\item[] $\DOset=\setof{\cf(otp(\pairof{X,R\restr X}))}{
	\begin{array}[t]{@{}l}
		X\subseteq\fnsp{\omega}{\omega},\,R\mbox{ is a definable binary relation }\\[\jot]
		\mbox{and }R\cap X^2\mbox{ well orders }X}
	\end{array}$
\end{xitemize}
and 
\begin{xitemize}
\item[] $\donum=\sup\DOset$.
\end{xitemize}
By definition, $\Es^\uparrow\subseteq\DOset$. Hence 
\begin{Lemma}
	$\boundingno^\uparrow\leq\donum$. 
\end{Lemma}
\noindent
\Lemmaof{b-b*} and \Lemmaabove\ may be put together into the following diagram:
\newpage
\[
\begin{array}{@{}r@{}l}
	&\donum\\
	&\rotatebox{90}{$\leq$}\\
	\boundingno\leq{}&\boundingno^\uparrow\leq\boundingno^h\leq\boundingno^*\leq\dominatingno
\end{array}\nopagebreak
\]\nopagebreak\mbox{}\hfill {\bf fig.\ 3}\hfill\mbox{}\bigskip

If $\Es^\uparrow$ has a maximal element then we have
$\boundingno^\uparrow=\max\Es^\uparrow$. In such case we shall say that 
$\boundingno^\uparrow$ is {\em attained}. Also we shall say that $\boundingno^*$,
$\boundingno^h$ or $\donum$ is attained if the corresponding set has a
maximal element.  

In the following, $\Reg$ denotes the class of regular cardinals. 
The following lemma can be proved similarly to 
\Lemmaof{bounding-5},\,\assertof{c}. 
\begin{Lemma}
{\rm (\cite{fuchino-x})}
\label{bounding-1}
$\Es^h\cap \Reg\subseteq\DOset$. 
%% In particular, if 
%% $2^{\aleph_0}<\aleph_\omega$ then $\Es^h\subseteq\DOset$.
\ifextended\else\qed\fi
\end{Lemma}
\ifextended
\noindent
{\bf [\,The following proof will be omitted in the version for 
		publication.]} 
\smallskip\\
{\footnotesize
\prf
For $\kappa\in\Es^h\cap \Reg$, let 
$X\in[\fnsp{\omega}{\omega}]^\kappa$ be as in the definition of $\Es^h$. 
Then we can construct $a_\alpha$, $b_\alpha$, $\alpha<\kappa$ inductively 
\st\ 
\begin{xitemize}
	\xitem[bounding-1.5] $a_\alpha\in X$ and $b_\alpha\in\fnsp{\omega}{\omega}$ for
	$\alpha<\kappa$. 
	\xitem[bounding-2] $a_\alpha\leq^* b_\beta$ $\Leftrightarrow$ $\alpha<\beta$.
\end{xitemize}
\xitemof{bounding-2} is possible since, at the $\beta$'th step in the 
inductive construction, we have 
$\cardof{\setof{a\in X}{a\leq^* b_\alpha\,\mbox{ for some }\alpha<\beta}}<\kappa$.  
Note that we need here that $\kappa$ is regular. 

Let $Y=\setof{\pairof{a_\alpha,b_\alpha}}{\alpha<\kappa}$ and let $R$ be the 
binary relation defined by 
\[ \pairof{a,b}\mathrel{R}\pairof{c,d}\ \Leftrightarrow\ a\leq^* d.
\]\noindent
for $\pairof{a,b}$, $\pairof{c,d}\in(\fnsp{\omega}{\omega})^2$
Clearly $R$ is projective and orders $Y$ in order type $\kappa$.\qedofLemma
\qedskip
\\
}\fi
\begin{Cor}
If $\Es^h\cap \Reg$ is cofinal in $\Es^h$ then
$\boundingno^h\leq\donum$. 	\qed
\end{Cor}
Note that the condition ``$\Es^h\cap \Reg$ is cofinal in $\Es^h$'' holds 
if $2^{\aleph_0}<\aleph_\omega$ or if $\boundingno^h$ is regular and 
attained. Under this condition, we can thus improve the diagram in fig.3 to the 
following:
\[
\begin{array}{@{}r@{}l}
	&\donum\\
	&\rotatebox{90}{$\leq$}\\
	\boundingno\leq\boundingno^\uparrow\leq{}&\boundingno^h\leq\boundingno^*\leq\dominatingno
\end{array}
\]\nopagebreak
\mbox{}\hfill{\bf fig.\ 4}\hfill\mbox{}\bigskip

For an ideal $I$ over a set $X$, $\non(I)$ and $\cov(I)$ denote, as usual, 
the uniformity 
and the covering number of $I$, respectively. More exactly 
\begin{xitemize}
	\item[] $\non(I)=\min\setof{\cardof{A}}{A\in\psof{X}\setminus I}$ and 
	\item[] $\cov(I)=\min\setof{\cardof{\calA}}{\calA\subseteq I,\,\bigcup\calA=X}$.
\end{xitemize}
$\cat$ and $\meas$ denote the ideal of meager sets and the ideal of null 
sets (over $\reals$) respectively. 

\begin{Lemma}\label{min-non-cov-leq-do} 
	Suppose that $I$ is an ideal over $\reals$ with Borel basis. Then we have
	$\min\ssetof{\non(I),\cov(I)}\leq\donum$. In particular, we have 
	$\min\ssetof{\non(\cat), \cov(\cat)}\leq\donum$ and 
	$\min\ssetof{\non(\meas), \cov(\meas)}\leq\donum$.
\end{Lemma}
\prf
Suppose that $I\subseteq\psof{\reals}$ is an ideal with a Borel basis and 
$\kappa=\min\ssetof{\non(I), \cov(I)}$. 
We can construct inductively a sequence
$\seqof{\pairof{f_\alpha, g_\alpha}}{\alpha<\kappa}$ \st\ 
\begin{xitemize}
	\xitem[] $f_\alpha$, $g_\alpha\in\fnsp{\omega}{\omega}$ for 
	all $\alpha<\kappa$; 
	\xitem[bounding-3-0] $g_\alpha$ codes a Borel set
	$X_\alpha\subseteq\fnsp{\omega}{\omega}$ \st\ $X_\alpha\in I$ and 
	$\setof{f_\beta}{\beta<\alpha}\subseteq X_\alpha$ ;
	\xitem[bounding-3-1] $f_\alpha\not\in \bigcup_{\beta<\alpha}X_\beta$ for 
	all $\alpha<\kappa$.  
\end{xitemize}
Note that \xitemof{bounding-3-0} is possible by $\kappa\leq\non(I)$ and 
\xitemof{bounding-3-1} by $\kappa\leq\cov(I)$. 

The sequence $\seqof{\pairof{f_\alpha, g_\alpha}}{\alpha<\kappa}$ is well 
ordered in order type $\kappa$ by the definable ordering: 
\[\pairof{f',g'}\leq\pairof{f,g}\ 
 \ \Leftrightarrow\ \ f'\mbox{ is an element of the Borel set coded by }g.
\]\noindent
It follows that $\kappa\leq\donum$. 
\qedofLemma
\qedskip

The following lemma shows the relations of cardinal numbers 
$\boundingno$, $\boundingno^\uparrow$, $\boundingno^h$, 
$\donum$ to the combinatorial principles introduced in 
\sectionof{the-principles}. 

\begin{Lemma}\label{bounding-5}

	\assert{a} {\rm (I.\ Juh\'asz, L.\ Soukup and 
		Z.\ Szentmikl\'ossy \cite{ju-so-sz})} If there is a $\leq^*$-chain of 
	length $\kappa$ then 
	$\neg\Cs(\kappa)$ and $\neg\hCs(\kappa)$. In particular,
	$\kappa\in\Es^\uparrow$ implies 
	$\neg \Cs(\kappa)$ and $\neg\hCs(\kappa)$. 

	\assert{b} $\Cs(\kappa)$ (or $\hCs(\kappa)$) implies
	$\boundingno^\uparrow\leq\kappa$. 
	If $\boundingno^\uparrow$ is attained then $\Cs(\kappa)$ 
	(or $\hCs(\kappa)$) implies
	$\boundingno^\uparrow<\kappa$. 

	\assert{c} If $\kappa\leq\lambda$ for some $\lambda\in\Es^h$ with
	$\cf\lambda\geq\kappa$ then $\neg\Cs(\kappa)$ and $\neg\hCs(\kappa)$. 

	\assert{d} If $\Es^h\cap \Reg$ is cofinal in $\Es^h$ then $\Cs(\kappa)$ 
	(or $\hCs(\kappa)$) implies $\boundingno^h\leq\kappa$. 
	If $\boundingno^h$ is regular and attained then $\Cs(\kappa)$ 
	(or $\hCs(\kappa)$) implies $\boundingno^h<\kappa$. 

	\assert{e} $\kappa\in\DOset$ implies $\neg\HP(\kappa)$.

	\assert{f} $\HP(\kappa)$ implies $\donum\leq\kappa$.
	If $\donum$ is attained then $\HP(\kappa)$ implies $\donum<\kappa$.
\end{Lemma}
\prf 
\ifextended
{{\bf [\,The following proof of \assertof{a} will be omitted in the version for 
		publication.]} 
\smallskip

{\footnotesize
\assertof{a}: 
Assume that there is a $\leq^*$-chain  of length $\kappa$ 
in $\fnsp{\omega}{\omega}$.  Then there is a sequence 
$\seqof{b_\alpha}{\alpha<\kappa}$ of subsets of $\omega$ \st\ 
$b_\alpha\subseteq^*b_\beta$ and $b_\beta\not\subseteq^*b_\alpha$ for all
$\alpha<\beta<\kappa$. For  
$\alpha<\kappa$ and $n\in\omega$, let 
\begin{xitemize}
	\xitem[x-2] $a_{\alpha,n}=\left\{\,
	\begin{array}[c]{@{}ll}
		\omega\setminus b_\alpha, &\mbox{if }n=0\\
		b_\alpha\setminus n, &\mbox{otherwise}
	\end{array}\right.
	$
\end{xitemize}
and $\calA=\seqof{a_{\alpha,n}}{\alpha<\kappa,n\in\omega}$. Let
$T=\setof{\pairof{0,n}}{n\in\omega}$. Then it is easy to see that
$\pairof{\calA, T}\not\models$~\assertof{c0} and
$\pairof{\calA, T}\not\models$~\assertof{c1}.
\smallskip

}\else\assertof{a}: See \cite{ju-so-sz}.\smallskip

\fi 
\assertof{b}: This follows from \assertof{a}.
\smallskip

\assertof{c}: 
Suppose that $\kappa\leq\lambda\in\Es^h$ and $\kappa\leq\cf\lambda$. We 
show $\neg\Cs(\kappa)$. $\neg\hCs(\kappa)$ can be proved similarly from 
these assumptions. 

Let $X\subseteq\fnsp{\omega}{\omega}$ with $\cardof{X}=\lambda$ be as in 
the definition of $\Es^h$. Then we can find $f_\alpha\in X$ and 
$g_\alpha\in\fnsp{\omega}{\omega}$ for $\alpha<\kappa$ \st\ 
\begin{xitemize}
	\xitem[x-3] $f_\alpha\leq^*g_\beta$ for all $\alpha<\beta<\kappa$;
	\xitem[x-4] $f_\beta\not\leq^* g^+_\alpha$ for 
	all $\alpha\leq\beta<\kappa$ where $g^+_\alpha$ is defined by 
	$g^+_\alpha(k)=g_\alpha(k)+1$ for all $k\in\omega$. 
\end{xitemize}
Note that \xitemof{x-4} is possible since $\cf(\cardof{X})\geq\kappa$. 

For $\alpha<\kappa$, let $g_{\alpha,n}\in\fnsp{\omega}{\omega}$, 
$n\in\omega$ be \st\ 
\begin{xitemize}
	\xitem[x-5] 
	$\setof{g_{\alpha,n}}{n\in\omega}
	=\setof{g\in\fnsp{\omega}{\omega}}{g=^*g_\alpha}$.  
\end{xitemize}
Let 
\begin{xitemize}
	\xitem[x-6]
	$a_{\alpha,0}=\setof{\pairof{k,\ell}\in\omega^2}{\ell\leq f_\alpha(k)}$ and
	\xitem[x-7] 
	$a_{\alpha,n+1}=\setof{\pairof{k,\ell}\in\omega^2}{\ell> g_{\alpha,n}(k)}$
	for all $n\in\omega$. 
\end{xitemize}
We show that $\calA=\seqof{a_{\alpha,n}}{\alpha\in\kappa,\,n\in\omega}$ 
with $T=\setof{\pairof{0,n}}{n\in\omega\setminus 1}$ is a counter-example 
to $\Cs(\kappa)$. 

Suppose first that $S\subseteq\kappa$ is stationary. For any $\alpha\in S$, 
let $\beta\in S$ be \st\ $\alpha<\beta$. Then we have
$f_\alpha\leq^* g_\beta$ by \xitemof{x-3}. Hence there is $n\in\omega$ \st\
$f_\alpha\leq g_{\beta,n}$. By \xitemof{x-6} and \xitemof{x-7}, it follows 
that $a_{\alpha,0}\cap a_{\beta,n+1}=\emptyset$. This shows that
$\pairof{\calA,T}\not\models\assertof{c0}$. 

Suppose now that $S_0$, $S_1\subseteq\kappa$ are stationary and
$\pairof{0,n}\in T$. By the definition of $T$, it follows that 
$n\in\omega\setminus 1$. Let $\alpha\in S_0$ and $\beta\in S_1$ be \st\
$\beta<\alpha$. Then, by \xitemof{x-4}, we have
$f_\alpha\not\leq^* g^+_\beta$. Thus, by \xitemof{x-6} and \xitemof{x-7}, 
it follows that $a_{\alpha,0}\cap a_{\beta, n}\not=\emptyset$. 
 This shows that
$\pairof{\calA,T}\not\models\assertof{c1}$. \smallskip

\assertof{d}: This follows easily from \assertof{c}. 
\smallskip

\assertof{e}:
Suppose that $\kappa\in\DOset$ and let $\pairof{X,R}$ be \st\
$X\subseteq\psof{\omega}$,  
$R$ is a projective  binary relation and $otp(\pairof{X,R\cap X^2})=\kappa$. 
Let $\mapping{f}{\kappa}{\psof{\omega}}$ be the mapping sending 
$\alpha<\kappa$ to the $\alpha$'th element of $X$ \wrt\ $R$. Let 
\[\textstyle 
A= R
\cup\bigcup_{k\in\omega\setminus \ssetof{2}}\injprod{\psof{\omega}}^k\,.
\]\noindent
Then it is easily seen that $\pairof{f,A}\not\models\assertof{h0}$ and
$\pairof{f,A}\not\models\assertof{h1}$. 
\smallskip

\assertof{f}: This follows from \assertof{d} since $\DOset$ is downward closed. 
\qedofLemma

\begin{Cor}\label{HP-non-cat}
	\assert{a} $\HP(\kappa)$ implies
	$\min\ssetof{\non(I), \cov(I)}\leq\kappa$ for any ideal $I$ over $\reals$ 
	with Borel basis. In particular, it implies 
\begin{xitemize}
	\item[] 
	$\min\ssetof{\non(\cat), \cov(\cat)}\leq\kappa$ and 
	$\min\ssetof{\non(\meas), \cov(\meas)}\leq \kappa$. 
\end{xitemize}

	\assert{b} If $\donum$ is attained 
	then $\HP(\kappa)$ implies 	$\min\ssetof{\non(I), \cov(I)}<\kappa$ for 
	all any $I$ over $\reals$ with Borel basis. In particular, it implies 
\begin{xitemize}
\item[] 
	$\min\ssetof{\non(\cat), \cov(\cat)}<\kappa$ and
	$\min\ssetof{\non(\meas), \cov(\meas)}<\kappa$.  
\end{xitemize}
\end{Cor}
\prf By \Lemmaof{min-non-cov-leq-do} and 
\Lemmaof{bounding-5},\,\assertof{f}. 
\qedofCor

\begin{Cor}\label{c=aleph-2+...}
	\assert{a} $\Cs(\aleph_2)$ (or $\hCs(\aleph_2)$)
	implies $\boundingno^h=\aleph_1$.

	\assert{b} $\HP(\aleph_2)$ implies
\begin{xitemize}
\item[] 
	$\donum=\min\ssetof{\non(\cat),\cov(\cat)}
	=\min\ssetof{\non(\meas), \cov(\meas)}=\aleph_1$.  
\end{xitemize}
\end{Cor}
\prf 
\assertof{a}: By \Lemmaof{bounding-5},\,\assertof{d}. 

\assertof{b}: By 
\Lemmaof{bounding-5},\,\assertof{f} and \Corof{HP-non-cat}. \qedofCor
\bigskip\bigskip

\noindent
\mbox{}\hfill\hspace{-1.6cm}
%%\input{coloring-x-fig3}
%WinTpicVersion3.08
\unitlength 1cm
\begin{picture}( 11.5625,  6.5175)( -0.6750, -7.5925)
% VECTOR 1 0 3 0
% 2 1556 3032 1028 3032
% 
\special{pn 13}%
\special{pa 1532 2985}%
\special{pa 1012 2985}%
\special{fp}%
\special{sh 1}%
\special{pa 1012 2985}%
\special{pa 1078 3004}%
\special{pa 1064 2985}%
\special{pa 1078 2965}%
\special{pa 1012 2985}%
\special{fp}%
% VECTOR 1 0 3 0
% 2 2788 3032 2260 3032
% 
\special{pn 13}%
\special{pa 2745 2985}%
\special{pa 2225 2985}%
\special{fp}%
\special{sh 1}%
\special{pa 2225 2985}%
\special{pa 2291 3004}%
\special{pa 2277 2985}%
\special{pa 2291 2965}%
\special{pa 2225 2985}%
\special{fp}%
% VECTOR 1 0 3 0
% 2 4020 3032 3492 3032
% 
\special{pn 13}%
\special{pa 3957 2985}%
\special{pa 3438 2985}%
\special{fp}%
\special{sh 1}%
\special{pa 3438 2985}%
\special{pa 3503 3004}%
\special{pa 3490 2985}%
\special{pa 3503 2965}%
\special{pa 3438 2985}%
\special{fp}%
% VECTOR 1 0 3 0
% 2 3136 2306 3136 2834
% 
\special{pn 13}%
\special{pa 3087 2270}%
\special{pa 3087 2790}%
\special{fp}%
\special{sh 1}%
\special{pa 3087 2790}%
\special{pa 3107 2724}%
\special{pa 3087 2738}%
\special{pa 3067 2724}%
\special{pa 3087 2790}%
\special{fp}%
% VECTOR 1 0 3 0
% 2 3136 1514 3136 2042
% 
\special{pn 13}%
\special{pa 3087 1491}%
\special{pa 3087 2010}%
\special{fp}%
\special{sh 1}%
\special{pa 3087 2010}%
\special{pa 3107 1944}%
\special{pa 3087 1958}%
\special{pa 3067 1944}%
\special{pa 3087 2010}%
\special{fp}%
% VECTOR 1 0 3 0
% 2 3136 678 3136 1206
% 
\special{pn 13}%
\special{pa 3087 668}%
\special{pa 3087 1188}%
\special{fp}%
\special{sh 1}%
\special{pa 3087 1188}%
\special{pa 3107 1122}%
\special{pa 3087 1135}%
\special{pa 3067 1122}%
\special{pa 3087 1188}%
\special{fp}%
% STR 2 0 3 0
% 3 3136 2086 3136 2174 5 0
% $\donum=\aleph_1$
\put(7.8400,-5.4350){\makebox(0,0){$\donum=\aleph_1$}}%
% STR 2 0 3 0
% 3 3136 1250 3136 1338 5 0
% $\HP(\aleph_2)$
\put(7.8400,-3.3450){\makebox(0,0){$\HP(\aleph_2)$}}%
% STR 2 0 3 0
% 3 3140 427 3140 515 5 0
% $\IP(\aleph_2,\aleph_1)$
\put(7.8500,-1.2875){\makebox(0,0){$\IP(\aleph_2,\aleph_1)$}}%
% STR 2 0 3 0
% 3 3140 2944 3140 3032 5 0
% $\boundingno^h=\aleph_1$
\put(7.8500,-7.5800){\makebox(0,0){$\boundingno^h=\aleph_1$}}%
% STR 2 0 3 0
% 3 1908 2944 1908 3032 5 0
% $\boundingno^\uparrow=\aleph_1$
\put(4.7700,-7.5800){\makebox(0,0){$\boundingno^\uparrow=\aleph_1$}}%
% STR 2 0 3 0
% 3 720 2944 720 3032 5 0
% $\boundingno=\aleph_1$
\put(1.8000,-7.5800){\makebox(0,0){$\boundingno=\aleph_1$}}%
% VECTOR 1 0 3 0
% 2 4350 1160 4350 2840
% 
\special{pn 13}%
\special{pa 4282 1142}%
\special{pa 4282 2796}%
\special{fp}%
\special{sh 1}%
\special{pa 4282 2796}%
\special{pa 4302 2730}%
\special{pa 4282 2744}%
\special{pa 4262 2730}%
\special{pa 4282 2796}%
\special{fp}%
% STR 2 0 3 0
% 3 4354 2953 4354 3041 5 0
% $\boundingno^*=\aleph_1$
\put(10.8850,-7.6025){\makebox(0,0){$\boundingno^*=\aleph_1$}}%
% STR 2 0 3 0
% 3 4350 882 4350 970 5 0
% $\WFN$
\put(10.8750,-2.4250){\makebox(0,0){$\WFN$}}%
% VECTOR 1 0 3 0
% 2 2946 1501 2515 2038
% 
\special{pn 13}%
\special{pa 2900 1478}%
\special{pa 2476 2006}%
\special{fp}%
\special{sh 1}%
\special{pa 2476 2006}%
\special{pa 2532 1968}%
\special{pa 2508 1965}%
\special{pa 2501 1942}%
\special{pa 2476 2006}%
\special{fp}%
% VECTOR 1 0 3 0
% 2 2600 2280 2980 2852
% 
\special{pn 13}%
\special{pa 2560 2245}%
\special{pa 2934 2808}%
\special{fp}%
\special{sh 1}%
\special{pa 2934 2808}%
\special{pa 2914 2742}%
\special{pa 2904 2764}%
\special{pa 2880 2764}%
\special{pa 2934 2808}%
\special{fp}%
% STR 2 0 3 0
% 3 2524 2064 2524 2152 5 0
% $\Cs(\aleph_2)$
\put(6.3100,-5.3800){\makebox(0,0){$\Cs(\aleph_2)$}}%
% VECTOR 1 0 3 0
% 2 2850 1501 2181 1906
% 
\special{pn 13}%
\special{pa 2806 1478}%
\special{pa 2147 1876}%
\special{fp}%
\special{sh 1}%
\special{pa 2147 1876}%
\special{pa 2213 1860}%
\special{pa 2192 1849}%
\special{pa 2193 1825}%
\special{pa 2147 1876}%
\special{fp}%
% VECTOR 1 0 3 0
% 2 2181 2117 2850 2830
% 
\special{pn 13}%
\special{pa 2147 2084}%
\special{pa 2806 2786}%
\special{fp}%
\special{sh 1}%
\special{pa 2806 2786}%
\special{pa 2775 2725}%
\special{pa 2770 2748}%
\special{pa 2747 2751}%
\special{pa 2806 2786}%
\special{fp}%
% STR 2 0 3 0
% 3 2084 1906 2084 1994 5 0
% $\hCs(\aleph_2)$
\put(5.2100,-4.9850){\makebox(0,0){$\hCs(\aleph_2)$}}%
% VECTOR 1 0 3 0
% 2 3017 674 2594 1210
% 
\special{pn 13}%
\special{pa 2970 664}%
\special{pa 2554 1191}%
\special{fp}%
\special{sh 1}%
\special{pa 2554 1191}%
\special{pa 2610 1152}%
\special{pa 2586 1150}%
\special{pa 2579 1127}%
\special{pa 2554 1191}%
\special{fp}%
% VECTOR 1 0 3 0
% 2 2445 1413 2110 1835
% 
\special{pn 13}%
\special{pa 2407 1391}%
\special{pa 2077 1807}%
\special{fp}%
\special{sh 1}%
\special{pa 2077 1807}%
\special{pa 2133 1767}%
\special{pa 2110 1765}%
\special{pa 2103 1743}%
\special{pa 2077 1807}%
\special{fp}%
% STR 2 0 3 0
% 3 2533 1228 2533 1316 5 0
% $\Fs(\aleph_2)$
\put(6.3325,-3.2900){\makebox(0,0){$\Fs(\aleph_2)$}}%
% STR 2 0 3 0
% 3 4495 2319 4495 2407 5 0
% \tiny\cite{fu-ge-so}
\put(11.2375,-6.0175){\makebox(0,0){\tiny\cite{fu-ge-so}}}%
% STR 2 0 3 0
% 3 3180 1822 3180 1910 2 0
% \tiny\Corof{c=aleph-2+...},\assertof{b}
\put(7.9500,-4.7750){\makebox(0,0)[lb]{\tiny\Corof{c=aleph-2+...},\assertof{b}}}%
% STR 2 0 3 0
% 3 2220 2512 2220 2600 5 0
% \tiny\Corof{c=aleph-2+...},\assertof{a}
\put(5.5500,-6.5000){\makebox(0,0){\tiny\Corof{c=aleph-2+...},\assertof{a}}}%
% VECTOR 3 0 3 0
% 2 4200 1130 2350 1890
% 
\special{pn 4}%
\special{pa 4134 1113}%
\special{pa 2313 1861}%
\special{fp}%
\special{sh 1}%
\special{pa 2313 1861}%
\special{pa 2381 1854}%
\special{pa 2362 1841}%
\special{pa 2367 1817}%
\special{pa 2313 1861}%
\special{fp}%
% VECTOR 3 0 3 0
% 2 4240 1170 2670 2040
% 
\special{pn 4}%
\special{pa 4174 1152}%
\special{pa 2628 2008}%
\special{fp}%
\special{sh 1}%
\special{pa 2628 2008}%
\special{pa 2695 1994}%
\special{pa 2675 1983}%
\special{pa 2677 1959}%
\special{pa 2628 2008}%
\special{fp}%
\end{picture}%
\hfill\nopagebreak\bigskip\\
\mbox{}\hfill {\bf fig.\ 5}\hfill\mbox{}

\ifextended
\fi
\section{A forcing construction of models of $\IP(\kappa,\lambda)$}
\label{forcing-construction}
In this section, we shall prove that $\IP(\kappa,\lambda)$ holds in a 
generic extension by a homogeneous product of copies of a relatively small partial 
ordering (\Thmof{forcing-IP*}). 

Let us begin with definition of some notions needed for precise formulation of 
the theorem. 

For cardinals $\kappa$ and $\mu$, $\kappa$ is said to be {\em$\mu$-inaccessible} 
if $\kappa$ is regular and $\lambda^\mu<\kappa$ holds for all
$\lambda<\kappa$. Similarly, we say that $\kappa$ 
is {\em${<}\mu$-inaccessible} if $\kappa$ is regular and
$\lambda^{<\mu}<\kappa$ 
holds for all $\lambda<\kappa$. Thus, if $\mu$ is a successor cardinal, say
$\mu=\mu_0^+$, then $\kappa$ is ${<}\mu$-inaccessible if and only if 
$\kappa$ is $\mu_0$-inaccessible. In our context, ${<}\mu$-inaccessibility is 
relevant because of the following variant of the $\Delta$-System Lemma of 
Erd\H os and Rado. 
For cardinals $\mu<\kappa$, let 
\begin{xitemize}
\item[]
	$E^\kappa_{\geq\mu}=\setof{\alpha<\kappa}{\cf(\alpha)\geq\mu}$ 
\end{xitemize}
and let $E^\kappa_\mu$, $E^\kappa_{\leq\mu}$ etc.\ be defined analogously. 
\begin{Thm}\label{delta-lemma}
	{\rm (P.\ Erd\H os and R.\ Rado, see \cite{ju-so-sz})}
	Suppose that $\kappa$ is ${<}\mu$-inaccessible and 
$S\subseteq E^\kappa_{\geq\mu}$ is stationary in $\kappa$. For any sequence
	$\seqof{x_\alpha}{\alpha\in S}$ of sets of cardinality $<\mu$ there is a 
	stationary $S^*\subseteq S$ \st\ $\seqof{x_\alpha}{\alpha\in S^*}$ form a
	$\Delta$-system. \ifextended\else\qed\fi
\end{Thm}
\ifextended
{{\bf [\,The following proof will be omitted in the version for 
		publication.]} 
\smallskip\\
{\footnotesize
\prf
\Wolog, we may assume that $x_\alpha\subseteq\kappa$ for every 
$\alpha\in S$. 
By Fodor's lemma and by assumption on $\kappa$ and $\mu$, we may 
further assume that there 
is a stationary $S'\subseteq S$ \st\ $\sup(x_\alpha)\geq\alpha$ for 
every $\alpha\in S'$. 

Let $\chi$ be sufficiently large. By the 
assumptions on $\kappa$ and $\mu$ there is $M\prec\calH(\chi)$ \st\ 
$\seqof{x_\alpha}{\alpha\in S}$, $S'$,$\kappa$, $\mu\in M$; 
$\cardof{M}<\kappa$; $[M]^{<\mu}\subseteq M$ and 
$M\cap\kappa\in S'$. Let $\alpha_1=M\cap\kappa$ and 
$r=x_{\alpha_1}\cap M$. Then $r\in M$. Hence $\alpha_0=(\sup r)+1$ is an 
element of $M$. 

For any $\alpha_0\leq\alpha^*<\alpha_1$, we have 
$\calH(\chi)\models\alpha^*<\alpha_1\land x_{\alpha_1}\cap\alpha^*=r$. 
Hence 
$\calH(\chi)\models
	\exists\beta<\kappa(\alpha^*<\beta\land x_\beta\cap\alpha^*=r)$. 
By elementarity, 
\begin{xitemize}
\item[] $M\models
	\exists\beta<\kappa(\alpha^*<\beta\land x_\beta\cap\alpha^*=r)$. 
\end{xitemize}
Since $\alpha_1=\kappa\cap M$, it follows that 
\begin{xitemize}
\item[] 
	$M\models\forall\alpha<\kappa
	\exists\beta<\kappa\,(\alpha<\beta\land x_\beta\cap\alpha=r)$.
\end{xitemize}
Again by elementarity the same sentence also holds in $\calH(\chi)$ and 
hence in $V$ as $\chi$ was taken sufficiently large. 

Now let $\xi_\alpha$, $\alpha<\kappa$ be the increasing sequence of 
ordinal numbers defined inductively by: 
$\xi_0=\mbox{the minimal }\xi\in S'\mbox{ \st\ }x_\xi\cap\alpha_0=r$
and for $\alpha>0$
\begin{xitemize}
\item[] 
	$\begin{array}{r@{\;}l}
	\xi_\alpha=&\mbox{the minimal }\xi\in S'\mbox{ \st\ \ }
		\xi>\sup(\bigcup\setof{x_{\xi_\beta}}{\beta<\alpha})\mbox{\ \ and\ \ }
	x_\xi\cap\xi=r.    
	\end{array}$
\end{xitemize}
Let 
$S^*=\setof{\xi_\alpha}{\alpha<\kappa}$. Then $S^*\subseteq S'$. 
Since the definition involves only parameters from 
$M$, we have $S^*\in M$. Clearly 
$\setof{x_\alpha}{\alpha\in S^*}$ forms a $\Delta$-system with the root 
$r$ and $S^*\in M$. Hence it is enough to show that $S^*$ is stationary. 

First note that $\alpha_1=\xi_{\alpha_1}$ and hence $\alpha_1\in S^*$. If 
$C\in M$ is a club subset of 
$\kappa$, then $\alpha_1\in C$. Hence $C\cap S^*\not=\emptyset$. By 
elementarity it follows that $M\models C\cap S^*\not=\emptyset$. Thus 
we have $M\models\xmbox{``$S^*$ is stationary''}$. Again by elementarity 
and by the choice of $\chi$ being sufficiently large, it follows that 
$S^*$ is really stationary.
\qedofThm\qedskip

}
\fi
For a sequence $\poP_\alpha$, $\alpha<\delta$ of \pos\ and an ideal
$I\subseteq\psof{\delta}$, we consider the {\em$I$-support 
product\/} $\prod^I_{\alpha<\delta}\poP_\alpha$ of $\poP_\alpha$, 
$\alpha<\delta$ here as being defined as 
\begin{xitemize}
	\xitem[prod-0] 
	$
	\begin{array}[t]{r@{}l}
		\prod^I_{\alpha<\delta}\poP_\alpha=
		\setof{f}{{}&\mapping{f}{D}{{}
				\bigcup_{\alpha<\delta}\poP_\alpha}
			\mbox{ for some }D\in I\\[\jot]
			&\mbox{and }f(\alpha)\in\poP_\alpha\setminus\ssetof{\bbbone_{\poP_\alpha}}
			\mbox{ for all }\alpha\in D\quad}		
	\end{array}
	$
\end{xitemize}
with the ordering
\begin{xitemize}
	\xitem[] 
	$
	\begin{array}[t]{r@{}l}
		f\leq_{\prod^I_{\alpha<\delta}\poP_\alpha} g\ \ \Leftrightarrow\ \ 
		&\dom(f)\supseteq\dom(g)\mbox{ and}\\[\jot]
		&f(\alpha)\leq_{\poP_\alpha}g(\alpha)\mbox{ for all }\alpha\in\dom(g)
	\end{array}
$
\end{xitemize}
for all $f$, $g\in\prod^I_{\alpha<\delta}\poP_\alpha$. 
In particular, $\bbbone_{\prod^I_{\alpha<\delta}\poP_\alpha}=\emptyset$ is 
the largest element of $\prod^I_{\alpha<\delta}\poP_\alpha$ \wrt\
$\leq_{\prod^I_{\alpha<\delta}\poP_\alpha}$\,.

Though this definition of product of \pos\ is different from the standard 
one, it gives a \po\ forcing equivalent to the product given by the standard 
definition. The present definition is chosen here for the sake of smoother 
treatment of $p\restr X$, $P\restr X$, $G\restr X$ etc.\ (see 
\xitemof{prod-0-0}, \xitemof{prod-1} etc.)

As usual, the ideal $[\delta]^{<\aleph_0}$ is denoted by $fin$ and
$\prod^{fin}_{\alpha<\delta}\poP_\alpha$ is called the {\em finite support 
	product\/} of $\poP_\alpha$, $\alpha<\delta$. 

I.~Juh\'asz and K.~Kunen \cite{juhasz-kunen} proved the following theorem 
for $\mu=\aleph_1$ and $I=[\delta]^{<\aleph_0}$. Their proof also applies 
to the following slight generalization.  
\begin{Thm}{\rm (I.~Juh\'asz and K.~Kunen \cite{juhasz-kunen})} 
	Suppose that $\poP=\prod^I_{\alpha<\delta}\poP_\alpha$ for some ideal 
	$I\subseteq\psof{\delta}$, $\poP$ satisfies the $\mu$-c.c.\ and 
	$\cardof{\poP_\alpha}\leq2^{<\mu}$ for all $\alpha<\delta$. Then, for all
	${<}\mu$-inaccessible $\kappa$ we have $\forces{\poP}{\Cs(\kappa)}$. 
	\ifextended\else\qed\fi
\end{Thm}
\ifextended
\noindent
{\bf [\,The following proof will be omitted in the version for 
		publication.]} 
\smallskip\\
{\footnotesize
\prf Let $G$ be a $(V,\poP)$-generic filter. 
Suppose that 
$\calA=\seqof{a_{\alpha,n}}{\alpha<\kappa,n\in\omega}$ is a matrix of subsets of 
$\omega$ and $T\subseteq\fnsp{\omega>}{\omega}$, both in $V[G]$. For 
$\alpha<\kappa$ and $n\in\omega$, let $\dot{a}_{\alpha,n}$ be a 
nice $\poP$-name of $a_{\alpha,n}$ and $\dot{T}$ a nice $\poP$-name 
of $T$. Let $s_{\alpha,n}=\supp(\dot{a}_{\alpha,n})$ for $\alpha<\kappa$ 
and $n\in\omega$. Since $\poP$ satisfies the $\mu$-c.c., each $s_{\alpha,n}$ is of 
cardinality $<\mu$.
Hence by \Thmof{delta-lemma}, there is a stationary $S\subseteq\kappa$ \st\
$\setof{\bigcup_{n\in\omega}s_{\alpha,n}}{\alpha\in S}$ forms 
a $\Delta$-system with root $R$. In $V[G]$ we can thin out $S$ so that it 
remains stationary but satisfies additionally that:
\begin{xitemize}
	\xitem[j-k-0'] for all $x\in \psof{\omega}\cap V[G_R]$ and for all
	$n\in\omega$,  
	$\setof{\alpha\in S}{a_{\alpha,n}\subseteq x}$ is either stationary or empty
\end{xitemize}
where $G_R=G\cap\poP_R$ with $\poP_R$ as defined below. 
This is possible since $(2^{\aleph_0})^{V[G_R]}<\kappa$ by
$\cardof{\poP_R}\leq 2^{<\mu}<\kappa$.

If $\bigcap_{k<\ellof{t}}a_{\alpha_k,t(k)}\not=\emptyset$ 
for all $t\in T$ and
$\pairof{\alpha_0\ctentenc \alpha_{\ellof{t}-1}}\in\injprod{S}^{<\omega}$, 
then we have \assertof{c0} for these $\calA$ and $T$ in $V[G]$. Otherwise 
there are 
$t^*\in T$ and
$\pairof{\alpha^*_0\ctentenc \alpha^*_{\ellof{t^*}-1}}\in\injprod{S}^{<\omega}$ 
\st\ 
\begin{xitemize}
	\xitem[] $\bigcap_{k<\ellof{t^*}}a_{\alpha^*_k,t^*(k)}=\emptyset$. 
\end{xitemize}
By Lemma 3.4 in \cite{juhasz-kunen}, we can find $x_0$\ctentenc 
$x_{\ellof{t^*}-1}\in V[G_R]$ \st\ 
\begin{xitemize}
	\xitem[] $V[G]\models a_{\alpha^*_k, t^*(k)}\subseteq x_k$ for all 
$k<\ellof{t^*}$ and 
	\xitem[j-k-empty] $\bigcap_{k<\ellof{t^*}}x_k=\emptyset$. 
\end{xitemize}
Now let 
\begin{xitemize}
	\xitem[] $S_k=\setof{\alpha<\kappa}{a_{\alpha,t^*(k)}\subseteq x_k}$
\end{xitemize}
for $k<\ellof{t^*}$. Then $\alpha^*_k\in S_k$. By 
\xitemof{j-k-0'}, it follows that each $S_k$ is stationary. 
By \xitemof{j-k-empty} and by the definition of $S_k$,
$\bigcap_{k<\ellof{t^*}}a_{t^*(k),\alpha_k}=\emptyset$ for all
$\pairof{\alpha_0\ctentenc \alpha_{\ellof{t^*}-1}}
\in \injprod{S_0\ctentenc S_{\ellof{t^*}-1}}$. 
Thus \assertof{c1} 
holds in $V[G]$ for these $\calA$ and $T$. \qedofThm\qedskip}\fi

Suppose that $I\subseteq\poP(\delta)$ is an ideal and 
$\poP=\prod^I_{\alpha<\delta}\poP_\alpha$ is an $I$-support product of 
\pos\ $\poP_\alpha$, $\alpha<\delta$. 
For $p\in\poP$, the {\em support\/} $\supp(p)$ of $p$  is defined by
\begin{xitemize}
	\xitem[] $\supp(p)=\dom(p)$. 
\end{xitemize}
We assume in the following that $\poP$-names are constructed just as in 
\cite{kunen-book}. For a $\poP$-name $\dot{a}$, the 
support $\supp(\dot{a})$ is defined by 
\begin{xitemize}
	\xitem[prod-0-a] $\supp(\dot{a})=\bigcup\setof{\supp(p)}{
		\pairof{\dot{b},p}\in tcl(\dot{a})\mbox{ for some }\poP\mbox{-name }\dot{b}\/}$.
\end{xitemize}
For $X\in\psof{\delta}$ (not necessarily in $I$), let 
\begin{xitemize}
	\xitem[prod-0-0] $\poP\restr X=\setof{p\restr X}{p\in\poP}$. 
\end{xitemize}
By \xitemof{prod-0} and since $I$ is an ideal, we have
\begin{xitemize}
	\xitem[] $\poP\restr X=\setof{p\in\poP}{\supp(p)\subseteq X}$.
\end{xitemize}
In particular, 
\begin{xitemize}
	\xitem[prod-1] $\poP\restr X\subseteq\poP$. 
\end{xitemize}
Furthermore, it is easy to see that
$\poP\restr X\circleq\poP$. Thus, if $G$ is a $(V,\poP)$-generic filter  
then $G\cap(\poP\restr X)$ is a $(V,\poP\restr X)$-generic filter. We shall 
denote the generic filter $G\cap(\poP\restr X)$ by $G_X$. 
Note that a $\poP$-name $\dot{a}$ is a $\poP\restr X$-name if and only if
$\supp(\dot{a})\subseteq X$.

We shall call an $I$-support product
$\poP=\prod^I_{\alpha<\delta}\poP_\alpha$ {\em homogeneous} if 
$\poP_\alpha\cong\poP_\beta$ for all $\alpha$, $\beta<\delta$ and $I$ is 
translation invariant, that is, $I=\setof{j\imageof x}{x\in I}$ for 
all bijections $\mapping{j}{\delta}{\delta}$. 

Note that if $I$ is translation invariant then $I=[\delta]<\lambda$ for some
$\lambda$. 

For a homogeneous $\poP=\prod^I_{\alpha<\delta}\poP_\alpha$, we shall 
always assume that a commutative system
$\isom{i_{\alpha,\beta}}{\poP_\alpha}{\poP_\beta}$,
$\alpha$, $\beta<\delta$ of isomorphisms is fixed. With such a fixed system 
of isomorphisms, every bijection $\mapping{j}{\delta}{\delta}$ 
induces an isomorphism $\isom{\tilde{j}}{\poP}{\poP}$ defined by 
\begin{xitemize}
	\xitem[] $\dom(\tilde{j}(p))=j\imageof\dom(p)$\,;\\[\jot]
	for 
	$\alpha\in\dom(\tilde{j}(p))$,\ \   
	$\tilde{j}(p)(\alpha)=
	i_{j^{-1}(\alpha),\alpha}\circ p\circ j^{-1}(\alpha)$ 
\end{xitemize}
for all $p\in\poP$. 

%% \mbox{}\\
%% \mbox{}\hspace{-1cm}\input{coloring-x-fig3x}\hfill\\
%% \mbox{}\hfill {\bf fig.\ 6}\hfill\mbox{}\bigskip\bigskip\\

For notational simplicity we shall denote the isomorphism on $\poP$-names 
induced from $\tilde{j}$ also by $\tilde{j}$. 

Note that for $\poP$ and $j$ as above, $p\in\poP$, 
$\poP$-names $\dot{a}_0$\ctentenc\ $\dot{a}_{n-1}$ and a formula $\varphi$ 
in the language of set theory $\calL_\ZF$, we have
\begin{xitemize}
	\xitem[homogeneous-P] 
	$p\forces{\poP}{\varphi(\dot{a}_0\ctentenc\dot{a}_{n-1})}$  
	if and only if 
	$\tilde{j}(p)\forces{\poP}{
		\varphi(\tilde{j}(\dot{a}_0)\ctentenc\tilde{j}(\dot{a}_{n-1}))}$. 
\end{xitemize}

We are now ready to formulate the main result of the present section: 
\begin{Thm}\label{forcing-IP*}
	Suppose that 
	\begin{xitemize}
		\xitem[otimes-0] $\lambda$ is a regular uncountable cardinal with
		$2^{<\lambda}=\lambda$, $\mu\in\ssetof{\lambda,\lambda^+}$ 
		and $\kappa$ is a ${<}\lambda$-inaccessible cardinal. 
	\end{xitemize}
	Let $\poP=\prod^I_{\alpha<\delta}\poP_\alpha$ be a 
	homogeneous $I$-support product \st\
	\begin{xitemize}
		\xitem[otimes-1] $I\subseteq[\delta]^{<\lambda}$;
		\xitem[otimes-2] $\cardof{\poP_\alpha}\leq\lambda$ for all $\alpha<\delta$ 
		and $\poP$ satisfies the $\mu$-c.c.;
%%%		\xitem[otimes-3] $\poP$ satisfies the $\kappa$-c.c.;
		\xitem[otimes-4] $\poP$ is proper.
	\end{xitemize}
Then $\forces{\poP}{\IP(\kappa,\mu)}$ holds. 
\end{Thm}
\newcounter{thetheorem}
\setcounter{thetheorem}{\value{Thm}}

The proof of 
\Thmof{forcing-IP*} will be given after the following Lemmas   
\ref{proper-slim} to \ref{homogeneous}. 

As in \cite{kunen-book}, a $\poP$-name $\dot{x}$ of a subset of $\omega$ 
for a \po\ $\poP$ 
is called a {\em nice $\poP$-name} if there are antichains $A_{\dot{x},n}$,
$n\in\omega$ in $\poP$ \st\
$\dot{x}=\setof{\pairof{\check{n},p}}{p\in A_{\dot{x},n}}$. 
Note that, for such a name $\dot{x}$, we have 
$\supp(\dot{x})=\bigcup_{n\in\omega}A_{\dot{x},n}$. 
It is easy to see that, 
for all $\poP$-names $\dot{x}$ of subsets of $\omega$, there is a nice
$\poP$-name $\dot{x}'$ \st\ $\forces{\poP}{\dot{x}=\dot{x}'}$. 
We say that a nice $\poP$-name of a subset of $\omega$ 
with $A_{\dot{x},n}$, $n\in\omega$ as above is {\em slim} 
if $A_{\dot{x},n}$ is countable for all $n<\omega$. 

The following lemmas are well-known:
\begin{Lemma}\label{proper-slim}
	Suppose that $\poP$ is a proper \po\ and $p\in\poP$. For any $\poP$-name 
$\dot{x}$ of a subset of $\omega$, there are $q\leq_\poP p$ and a slim
	$\poP$-name $\dot{x}'$ \st\ $q\forces{\poP}{\dot{x}=\dot{x}'}$. 
\end{Lemma}
\prf By the remark above, we may assume \wolog\ that $\dot{x}$ is a nice 
$\poP$-name. Let $A_{\dot{x},n}$, $n\in\omega$ be as above and 
$\dot{y}$ be a $\poP$-name \st\ 
$\forces{\poP}{\dot{y}=\setof{s\in \poP}%
	{s\in(\bigcup_{n\in\omega}A_{\dot{x},n})\cap\dot{G}}}$. Then we have 
$\forces{\poP}{\dot{y}\xmbox{ is a countable subset of }\poP}$. As $\poP$ is proper 
there exist $q\leq_\poP p$ and countable $y\subseteq \poP$ \st\ 
$q\forces{\poP}{\dot{y}\subseteq y}$. Let 
$\dot{x}'=
	\setof{\pairof{\check n,s}}{n\in\omega,\ s\in A_{\dot{x},n}\cap y}$. These 
$q$ and $\dot{x}'$ are as desired. 
\qedofLemma
\begin{Lemma}\label{stat-in-G}
	Suppose that $\poP=\prod^I_{\alpha<\delta}\poP_\alpha$ is a $\kappa$-c.c.\
	$I$-support product for an ideal $I\subseteq[\delta]^{<\lambda}$ and $\kappa$ is
	${<}\lambda$-inaccessible.   
%% is a regular cardinal and 
%% 	$\poP$ is an $I$-supported product for an ideal $I$ over the 
%% 	index set $K=\bigcup I$, $I\subseteq[K]^{<\lambda}$, $\kappa$ is
%% 	${<}\lambda$-inaccessible and 
%% 	\begin{xitemize}
%% 		\xitem[fhp-7] $\cardof{\poP_X}<\kappa$ for all $X\in [K]^{{<}\lambda}$.
%% 	\end{xitemize}
	If $S\subseteq E^\kappa_{\geq\lambda}$ is stationary and $p_\alpha\in\poP$ for
	$\alpha\in S$ are \st\ $\supp(p_\alpha)$, $\alpha\in S$ form 
	a $\Delta$-system with the root $R$ and there is $p^*\in\poP\restr R$ 
	\st\ $p_\alpha\restr R=p^*$ for all
	$\alpha\in S$, then 
	\[ p^*\forces{\poP}{\setof{\alpha\in S}{p_\alpha\in \dot{G}}
		\mbox{ is stationary}}
	\]\noindent
	where $\dot{G}$ denotes the standard $\poP$-name of 
	a $(V,\poP)$-generic filter. 
\end{Lemma}
\prf By $\kappa$-c.c.\ of $\poP$, 
$\kappa$ remains a regular cardinal in 
$\poP$-generic extensions. 
Let $\dot{S}$ be a $\poP$-name of
$\setof{\alpha\in S}{p_\alpha\in \dot{G}}$. 
Suppose that $\dot{C}$ is a $\poP$-name of a 
club subset of $\kappa$ and $p\leq_\poP p^*$. It is enough to show that 
there is a $q\leq_\poP p$ \st\
$q\forces{\poP}{\dot{C}\cap\dot{S}\not=\emptyset}$. 

Let $\theta$ be sufficiently large and let $M\prec\calH(\theta)$  be \st\ 
\begin{xitemize}
	\xitem[fhp-3] $I$, $\poP$, $\kappa$, $\dot{C}$,
	$\seqof{p_\alpha}{\alpha\in S}$, $p\in M$;
	\xitem[fhp-5] $\cardof{M}<\kappa\cap M<\kappa$;
	\xitem[fhp-8] $[M]^{<\lambda}\subseteq M$ and 
	\xitem[fhp-6] $\alpha^*\in S$ where $\alpha^*=\kappa\cap M$.
\end{xitemize}
\xitemof{fhp-8} is possible since $\kappa$ is ${<}\lambda$-inaccessible. 
\xitemof{fhp-6} is possible since $S\subseteq E^\kappa_{\geq\lambda}$ and $S$ is 
stationary in $\kappa$.  

\begin{Claim}
	$\forces{\poP}{\alpha^*\in \dot{C}}$.
\end{Claim}
\prfofClaim  
Since 
$\poP$ satisfies the $\kappa$-c.c., we have 
\[\calH(\theta)\models\forall\alpha<\kappa\,\exists\beta\in\kappa\setminus\alpha\ 
(\forces{\poP}{\beta\in\dot{C}}).
\]\noindent
By \xitemof{fhp-3}, and elementarity of $M$ it follows that
\[ M\models\forall\alpha<\kappa\,\exists\beta\in\kappa\setminus\alpha\ 
(\forces{\poP}{\beta\in\dot{C}}).
\]\noindent
Thus $\forces{\poP}{\dot{C}\cap\alpha^*\mbox{ is unbounded in }\alpha^*}$. 
Since $\forces{\poP}{\dot{C}\xmbox{ is a club in }\kappa}$, it follows 
that $\forces{\poP}{\alpha^*\in \dot{C}}$. \qedofClaim
\begin{Claim}
	$\supp(p_{\alpha^*})\cap M=R$ and $\supp(p)\cap\supp(p_{\alpha^*})=R$\,. 
\end{Claim}
\prfofClaim
Suppose $u=(\supp(p_{\alpha^*})\cap M)\setminus R\not=\emptyset$. By 
\xitemof{fhp-8}, $u\in M$. Hence by elementarity 
$M\models\exists\alpha<\kappa\ (u\subseteq \supp(p_\alpha))$. Let
$\alpha\in\kappa\cap M$ be \st\ $u\subseteq\supp(p_\alpha)$. Then 
$\alpha<\alpha^*$ and $R\cup u\subseteq p_\alpha\cap p_{\alpha^*}$. This is 
a contradiction to the assumption 
that $R$ is the root of the $\Delta$-system
$\setof{\supp(p_\alpha)}{\alpha\in S}$. This shows
$\supp(p_{\alpha^*})\cap M=R$. 

By \xitemof{fhp-3} and \xitemof{fhp-8}, $\supp(p)\in M$. It follows that 
$\supp(p)\cap\sup(p_{\alpha^*})=\supp(p)\cap(\sup(p_{\alpha^*})\cap M)
=\supp(p)\cap R=R$. 

\qedofClaim\qedskip 
\\
Since 
$p\restr R\leq_\poP p^*\restr R=p^*=p_{\alpha^*}\restr R$, 
$q=p\cup p_{\alpha^*}\in\poP$. We have $q\leq_\poP p$. 
By
$p_{\alpha^*}\forces{\poP}{{\alpha^*}\in\dot{S}}$ and
$q\leq_\poP p_{\alpha^*}$, we have  
$q\forces{\poP}{\alpha^*\in\dot{C}\cap\dot{S}}$. 
In particular $q\forces{\poP}{\dot{C}\cap\dot{S}\not=\emptyset}$. 
\qedofLemma
\qedskip

The arguments of the following two lemmas are also well-known. For 
\Lemmaof{injectivity-lemma} 
see e.g.\ \cite{juhasz-kunen}. 
\begin{Lemma}\label{injectivity-lemma}
	Suppose that 
	$\poP=\prod^I_{\alpha<\delta}\poP_\alpha$ is an $I$-support 
	product and $G$ is a $(V,\poP)$-generic 
	filter. For $X$, $Y\subseteq \delta$, let $Z=X\cap Y$. Then, in $V[G]$, for 
	any $\kappa\in\Card^{V[G]}$, we have 
	\begin{xitemize}
	\item[] $[\On]^{<\kappa}\cap\left(V[G_X]\setminus V[G_Z]\right)\cap
	\left(V[G_Y]\setminus V[G_Z]\right)=\emptyset$.\qed
	\end{xitemize}
\end{Lemma}
\begin{Lemma}\label{homogeneous}
	Suppose that $\kappa\leq\delta$ and
	$\poP=\prod^I_{\alpha<\delta}\poP_\alpha$ is a $\kappa$-c.c.\ homogeneous
	$I$-support product, $p\in\poP$, $\dot{a}_0\ctentenc \dot{a}_{n-1}$ 
	are $\poP$-names with
	\begin{xitemize}
		\xitem[hh-0] $\supp(\dot{a}_0)\ctentenc \supp(\dot{a}_{n-1})\subseteq X$  
	\end{xitemize}
	for some $X\subseteq\delta$ and $\varphi=\varphi(x_0\ctentenc x_{n-1})$ 
	is a formula in $\calL_\ZF$ (possibly with some parameters from $V$). \smallskip 

	\assert{a} If 
	\begin{xitemize}
		\xitem[h-a-0] $p\forces{\poP}{\varphi(\dot{a}_0\ctentenc\dot{a}_{n-1})}$
		and 
		\xitem[h-a-1] $\delta\setminus X\not\in I$, 
	\end{xitemize}
	then $p\restr X\forces{\poP}{\varphi(\dot{a}_0\ctentenc\dot{a}_{n-1})}$. 
	\smallskip

	\assert{b} If 
	\begin{xitemize}
		\xitem[h-a-2] 	$p\forces{\poP}{
		(\exists x\in\fnsp{\omega}{\omega})\ 
		\varphi(x,\dot{a}_1\ctentenc\dot{a}_{n-1})}$, 
		\xitem[h-0] $\supp(p)\subseteq X$ and 
		\xitem[h-1] 
		$\cardof{X\setminus(\supp(p)\cup\supp(\dot{a}_0)\cup\cdots\cup
			\supp(\dot{a}_{n-1}))}\geq\kappa$, 
	\end{xitemize}
	then there is a $\poP_X$-name $\dot{a}$ \st\ 
	$p\forces{\poP}{\varphi(\dot{a},\dot{a}_1\ctentenc\dot{a}_{n-1})}$. 
\end{Lemma}
\prf \assertof{a}: Suppose that
$p\restr X\notforces{\poP}{\varphi(\dot{a}_0\ctentenc\dot{a}_{n-1})}$. Then 
there is $q\leq_\poP p\restr X$ \st\
$q\forces{\poP}{\neg\varphi(\dot{a}_0\ctentenc\dot{a}_{n-1})}$.
Let $\mapping{j}{\delta}{\delta}$ be a bijection \st\ 
\begin{xitemize}
	\xitem[h-1-0] $j\restr X=id_X$ and 
	\xitem[h-1-1] $(j\imageof\supp(q)\setminus X)\ \cap\ \supp(p)=\emptyset$. 
\end{xitemize}
Note that the last condition is possible by \xitemof{h-a-1}.
By \xitemof{h-1-0} and \xitemof{hh-0}, we have
\begin{xitemize}
	\xitem[h-1-2] $\tilde{j}(q)\restr X=\tilde{j}(q\restr X)=q\restr X$ and 
	\xitem[h-1-3] $\tilde{j}(\dot{a}_0)=\dot{a}_0$\ctentenc\,
	$\tilde{j}(\dot{a}_{n-1})=\dot{a}_{n-1}$. 
\end{xitemize}
By \xitemof{h-1-3} and by the choice of $q$, we have
$\tilde{j}(q)\forces{\poP}{\neg\varphi(\dot{a}_0\ctentenc\dot{a}_{n-1})}$. 
On the other hand, by \xitemof{h-1-2} $p$ and $\tilde{j}(q)$ are 
compatible. This is a contradiction to \xitemof{h-a-0}.\smallskip

\assertof{b}: 
By maximal principle, there is a nice $\poP$-name $\dot{a}'$ of a real \st\ 
\begin{xitemize}
\item[] $p\forces{\poP}{\varphi(\dot{a}',\dot{a}_1\ctentenc\dot{a}_{n-1})}$. 
\end{xitemize}
By the $\kappa$-c.c.\ of $\poP$, we have $\cardof{\supp(\dot{a}')}<\kappa$. 
By \xitemof{hh-0}, \xitemof{h-0} and 
\xitemof{h-1}, we can find a bijection $\mapping{j}{\delta}{\delta}$ \st\ 
\begin{xitemize}
	\xitem[h-2] $j$ on $\supp(p)\cup\supp(\dot{a}_1)\cup\cdots\cup
			\supp(\dot{a}_{n-1})$ is the identity mapping, and
	\xitem[h-3] $j\imageof\supp(\dot{a}')\subseteq X$.
\end{xitemize}
By \xitemof{h-2}, $\tilde{j}(p)=p$ 
and $\tilde{j}(\dot{a}_1)=\dot{a}_1$\ctentenc\  
$\tilde{j}(\dot{a}_{n-1})=\dot{a}_{n-1}$. 
Let $\dot{a}=\tilde{j}(\dot{a}')$. 
Then $p\forces{\poP}{\varphi(\dot{a},\dot{a}_1\ctentenc \dot{a}_{n-1})}$
and $\supp(\dot{a})\subseteq X$ by \xitemof{h-3}. \qedofLemma
\qedskip

%---------------------------------------------------------------------------
\noindent%
\newcounter{oldthm}%
\setcounter{oldthm}{\value{Thm}}%
\setcounter{Thm}{\value{thetheorem}}%
\prfof{\bfThmof{forcing-IP*}}
%% We consider only the case of $\mu=\lambda$. The proof for $\mu=\lambda^+$ is 
%% exactly the same except that $\lambda$ is to be replaced by $\lambda^+$ at several 
%% places. 
By \Propof{continuum-less-than-kappa}, we may assume that
$\forces{\poP}{\kappa\leq 2^{\aleph_0}}$. 
In particular, by \xitemof{otimes-1},
\xitemof{otimes-2} and \xitemof{otimes-4}, we may assume that $\delta\geq\kappa$. 
By  
the $\mu$-c.c.\ of $\poP$, $\mu$ and $\kappa$ remain regular cardinals in 
the generic extension by $\poP$. 

Let $G$ be a $(V,\poP)$-generic filter. In $V[G]$, let 
$\mapping{f}{\kappa}{\psetof{\omega}}$ and  
$\mapping{g}{\injprod{\psetof{\omega}}^{<\omega}}{\psetof{\omega}}$ be definable, 
say by a formula $\varphi$. We may 
assume that 
$\varphi$ has a real $a\in V[G]$ as  its unique parameter. 
Let $\dot{f}$, $\dot{a}$ and $\dot{g}$ be $\poP$-names of $f$, $a$ 
and $g$ respectively \st\  
$\forces{\poP}{\mapping{\dot{f}}{\kappa}{\psetof{\omega}}}$, 
$\forces{\poP}{
	\mapping{\dot{g}}{\injprod{\psetof{\omega}}^{<\omega}}{\psetof{\omega}}}$ 
and 
\begin{xitemize}
	\xitem[hIp-14] $\forces{\poP}{
		\forall\overline{x}\in\injprod{\psetof{\omega}}^{<\omega}\ 
		\forall x\in\psetof{\omega}\ \bigl({\dot{g}(\overline{x})=x\leftrightarrow
			\calH(\aleph_1)\models\varphi(\overline{x},x,\dot{a})}\bigr)}$.
\end{xitemize}
Suppose that, for a $p\in G$, 
\begin{xitemize}
	\xitem[hIp-14A]  
	$p\forces{\poP}{
		\mbox{\assertof{i0} for }\IP(\kappa,\mu)
		\mbox{ does not hold for }\dot{f}\mbox{ and }\dot{g}}$.
\end{xitemize}

In particular, we have
\begin{xitemize}
	\xitem[hIp-15] $ p\forces{\poP}{%
	\forall \alpha<\kappa\ (\setof{\beta\in\kappa}{\dot{f}(\beta)=\dot{f}(\alpha)}
	\mbox{ is non-stationary})}$.
\end{xitemize}

\begin{Claim}\label{hIp-25}
There is a stationary $S\subseteq E^\kappa_{\geq\lambda}$ \st\ 
\[ p\forces{\poP}{\dot{f}\restr S\mbox{\rm\ is 1-1}}.
\]\noindent
\end{Claim}
\prfofClaim
By the $\kappa$-c.c.\ of $\poP$ and by \xitemof{hIp-15}, there are club sets 
$C_\alpha\subseteq\kappa$ (in $V$) for each $\alpha<\kappa$ \st\
\[ p\forces{\poP}{%
	C_\alpha\cap\setof{\beta\in\kappa}{\dot{f}(\beta)=\dot{f}(\alpha)}=\emptyset}.
\]\noindent
Then 
$C=\Delta_{\alpha<\kappa}C_\alpha$ is club and $S=E^\kappa_{\geq\lambda}\cap C$ 
has the desired property. \\\qedofClaim\qedskip\\
We show that $p$ forces \assertof{i1} for $\dot{f}$ and 
$\dot{g}$. Let $p'\leq_\poP p$. It is enough to show that there is 
$p^*\leq_\poP p'$ forcing \assertof{i1}. 

By \Lemmaof{proper-slim}, \Thmof{delta-lemma}, \xitemof{otimes-0}, 
\xitemof{otimes-1} and \xitemof{otimes-4},  
there are $p''\leq_\poP p'$, a slim $\poP$-name $\dot{a}'$ of a real, 
a stationary $S^*\subseteq S$, a sequence 
$\seqof{\dot{x}'_\alpha}{\alpha\in S^*}$ of slim $\poP$-names and a 
sequence $\seqof{p_\alpha}{\alpha\in S^*}$ of conditions in $\poP$ \st\ 
\begin{xitemize}
	\xitem[hIp-16] \assert{i} $p''\forces{\poP}{\dot{a}=\dot{a}'}$,\\
	\assert{ii} $p_\alpha\leq_\poP p''$ and\\
	\assert{iii}
	$p_\alpha\forces{\poP}{\dot{f}(\alpha)=\dot{x}'_\alpha}$ for every
	$\alpha\in S^*$; 
	\xitem[hIp-17] 
	$d_\alpha=\supp(p_\alpha)\cup\supp(\dot{a}')\cup\supp(\dot{x}'_\alpha)$, 
	$\alpha\in S^*$ are all of the same cardinality and form 
	a $\Delta$-system with root $R$; 
	\xitem[hIp-18] for each $\alpha$, $\beta\in S^*$ there is a bijection 
	$\mapping{j_{\alpha,\beta}}{\delta}{\delta}$ \st\\ 
	\assert{i}
	$j_{\alpha,\beta}\restr\left(\delta\setminus(d_\alpha\Delta d_\beta)\right)
		=id_{\delta\setminus(d_\alpha\Delta d_\beta)}$\,,\\
	\assert{ii} $j_{\alpha,\beta}\imageof d_\alpha=d_\beta$, 
	$\tilde{j}_{\alpha,\beta}(p_\alpha)=p_\beta$ and\\
	\assert{iii}
	$\tilde{j}_{\alpha,\beta}(\dot{x}'_\alpha)=\dot{x}'_\beta$ for every 
	$\alpha$, $\beta\in S^*$. 
\end{xitemize}

Note that, by \xitemof{hIp-17}, we have
\begin{xitemize}
	\xitem[] 	$\supp(\dot{a}')=\supp(\dot{a}')\cap d_\alpha\subseteq R$ for every
	$\alpha\in S^*$. 
\end{xitemize}

By \xitemof{hIp-18}, $p_\alpha\restr R$ for $\alpha\in S^*$ are all the 
same.  
Let $q=p_\alpha\restr R$ for some/any $\alpha\in S^*$. $q\leq_\poP p''$ by 
\xitemof{hIp-16}, \assertof{ii}. Let 
$\dot{S}$ be a $\poP$-name \st\ 
\begin{xitemize}
	\xitem[hIp-19] 
	$\forces{\poP}{\dot{S}=\setof{\alpha\in S^*}{p_\alpha\in\dot{G}}}$.
\end{xitemize}
By \Lemmaof{stat-in-G}, 
$q \forces{\poP}{\dot{S}\xmbox{ is stationary}}$. 
Hence, by \xitemof{hIp-14A},
\begin{xitemize}
	\xitem[] $ q \forces{\poP}{
		\exists n\in\omega\,\forall\alpha<\kappa\ \big(\cardof{
			\dot{g}\imageof\injprod{\dot{f}\imageof(\dot{S}\setminus\alpha)}^n}
		\geq\mu\big)}$\,.
\end{xitemize}
Let $q'\leq_\poP q $ and $n^*\in\omega$ be \st\ 
\begin{xitemize}
	\xitem[hIp-20A] 
	 $q'\forces{\poP}{\forall\alpha<\kappa\ (\cardof{
			\dot{g}\imageof\injprod{\dot{f}\imageof(\dot{S}\setminus\alpha)}^{n^*}}
		\geq\mu)}$.
\end{xitemize}
%% We may assume that 
%% \begin{xitemize}
%% 	\xitem[hIp-20X] $
%% 	\begin{array}[t]{r@{}l}
%% 		q'\forces{\poP}{{}&\mbox{there are pairwise disjoint }
%% 		u_\alpha\in\injprod{\dot{f}\imageof\dot{S}}^{n^*},\,\alpha<\mu
%% 		\mbox{ \st}\\[\jot]
%% 		&g\imageof u_\alpha,\alpha<\mu\mbox{ are distinct}
%% 		},		
%% 	\end{array}
%% 	$
%% \end{xitemize}
%% (otherwise modify $\dot{g}$ so that some of the $n^*$ coordinates of 
%% $\dot{g}\restr\injprod{\dot{f}\imageof\dot{S}}^{n^*}$ are fixed as 
%% parameters). Here, we mean with $\pairof{x_0\ctentenc x_{n-1}}$, 
%% $\pairof{y_0\ctentenc y_{n-1}}\in\injprod{X}^n$ being disjoint that 
%% $\ssetof{x_0\ctentenc x_{n-1}}$ and $\ssetof{y_0\ctentenc y_{n-1}}$ are 
%% disjoint. 
Let 
\begin{xitemize}
	\xitem[hIp-21]
	$S^{**}=\setof{\alpha\in S^*}{\supp(q')\cap d_\alpha\subseteq R}$.  
\end{xitemize}
Since $\cardof{\supp(q')}<\lambda$ by \xitemof{otimes-1}, we have
\begin{xitemize}
	\xitem[star] $S^*\setminus S^{**}$ is of cardinality $<\lambda$. 
\end{xitemize}
In particular $S^{**}$ is still stationary and 
\begin{xitemize}
	\xitem[hIp-20B]
 	$q'\forces{\poP}{\cardof{
			\dot{g}\imageof\injprod{\dot{f}\imageof(\dot{S}\cap S^{**})}^{n^*}}\geq\mu}$
\end{xitemize}
by \xitemof{hIp-20A}. 
\begin{Claim}\label{not-in-V[Gr]}
There is $\pairof{\alpha_0\ctentenc\alpha_{n^*-1}}\in \injprod{S^{**}}^{n^*}$  
\st\ 
\begin{xitemize}
\item[] $q'\cup p_{\alpha_0}\cup\cdots\cup p_{\alpha_{n^*-1}}\,
	\notforces{\poP}{%
	\dot{g}(\pairof{\dot{x}'_{\alpha_0}\ctentenc\dot{x}'_{\alpha_{n^*-1}}})\in 
	V[\dot{G}_R]}$.   
\end{xitemize}
\end{Claim}
\prfofClaim
%% Let 
%% \begin{xitemize}
%% 	\xitem[hIp-26-0] $q''=q'\cup p_{\alpha_0}\cup\cdots\cup p_{\alpha_{n^*-1}}$. 
%% \end{xitemize}
Otherwise, we would have
\begin{xitemize}
\item[] $q'\cup p_{\beta_0}\cup\cdots\cup p_{\beta_{n^*-1}}\forces{\poP}{%
	\dot{g}(\pairof{\dot{x}'_{\beta_0}\ctentenc\dot{x}'_{\beta_{n^*-1}}})\in 
	V[\dot{G}_R]}$
\end{xitemize}
for all
$\pairof{\beta_0\ctentenc\beta_{n^*-1}}\in \injprod{S^{**}}^{n^*}$.  

%% Let $G$ be a $(V,\poP)$-generic filter \st\ $q'\in G$. 
%% Suppose that
%% $\pairof{\beta_0\ctentenc \beta_{n^*-1}}\in \injprod{\dot{S}^G}^{n^*}$. 
%% By \xitemof{hIp-19}, we have
%% $p_{\beta_0}\ctentenc p_{\beta_{n^*-1}}\in G$. 
%% Hence $q'\cup p_{\beta_0}\cup\cdots\cup p_{\beta_{n^*-1}}\in G$. 
%% By assumption it follows that
%% $\dot{g}^G(\dot{f}^G(\beta_0)\ctentenc\dot{f}^G(\beta_{n^*-1}))\in V[G_R]$. 

%% Since 
%% $\cardof{\fnsp{\omega}{\omega}^{V[G_R]}}<\lambda$ by 
%% \xitemof{otimes-2}, this is a contradiction to \xitemof{hIp-20A} and 
%% \xitemof{star}.\\
Fix $\pairof{\alpha_0\ctentenc\alpha_{n^*-1}}\in\injprod{S^{**}}^{n^*}$ and 
let
\begin{xitemize}
\item[] $
	\begin{array}{r@{}l}
		\calD=\setof{r\in\poP}{{}
			&r\leq_\poP q'\cup p_{\alpha_0}\cup\cdots\cup p_{\alpha_{n^*-1}},\\[\jot]
		&\supp(r)\subseteq R\cup\bigcup\setof{d_{\alpha_i}}{i<n^*}\cup\supp(q'),\\[\jot]
		&r\forces{\poP}{\dot{g}(\dot{x}'_{\alpha_0}\ctentenc \dot{x}'_{\alpha_{n^*-1}})
			=\dot{x}}
			\mbox{ for some }\poP_R\mbox{-name }\dot{x}\ }.
	\end{array}
	$
\end{xitemize}
Let $\calA$ be a maximal antichain in $\calD$. 
By the $\mu$-c.c.\ of $\poP$, $\cardof{\calA}<\mu$. 

For each $r\in\calA$, let $\dot{x}_r$ be a $\poP_R$-name \st\ 
\begin{xitemize}
\item[] $r\forces{\poP}{\dot{g}(\dot{x}'_{\alpha_0}\ctentenc \dot{x}'_{\alpha_{n^*-1}})
			=\dot{x}_r}$
\end{xitemize}
and $\dot{\calX}$ be a $\poP_R$-name \st\
\begin{xitemize}
\item[] $\forces{\poP}{\dot{\calX}=\setof{\dot{x}_r}{r\in\calA}}$. 
\end{xitemize}
Then, we have $\forces{\poP}{\cardof{\dot{\calX}}<\mu}$. 
By \Lemmaof{homogeneous},\,\assertof{a}
\begin{xitemize}
\item[] $q'\cup p_{\alpha_0}\cup\cdots\cup p_{\alpha_{n^*-1}}
	\forces{\poP}{\dot{g}(\pairof{\dot{x}'_{\alpha_0}\ctentenc
		\dot{x}'_{\alpha_{n^*-1}}})\in\dot{\calX}}$. 
\end{xitemize}
Hence by \xitemof{hIp-18} and \xitemof{homogeneous-P}, we have 
\begin{xitemize}
\item[] $q'\cup p_{\beta_0}\cup\cdots\cup p_{\beta_{n^*-1}}
	\forces{\poP}{\dot{g}(\pairof{\dot{x}'_{\beta_0}\ctentenc
		\dot{x}'_{\beta_{n^*-1}}})\in\dot{\calX}}$ 
\end{xitemize}
for all $\pairof{\beta_0\ctentenc \beta_{n^*-1}}\in\injprod{S^{**}}^{n^*}$. 
But this is a contradiction to \xitemof{hIp-20B}. 
\qedofClaim\qedskip

Let $\pairof{\alpha_0\ctentenc \alpha_{n^*-1}}\in \injprod{S^{**}}^{n^*}$ 
be as in \Claimabove\ and 
\begin{xitemize}
	\xitem[hIp-26-0] $q''=q'\cup p_{\alpha_0}\cup\cdots\cup p_{\alpha_{n^*-1}}$. 
\end{xitemize}
Note that 
\begin{xitemize}
	\xitem[hIp-22] $q''\forces{\poP}{f(\alpha_i)=\dot{x}'_{\alpha_i}}$ for
	$i<n^*$ by \xitemof{hIp-26-0}. 
\end{xitemize}

Let $p^*\leq_\poP q''$ be \st\ 
\begin{xitemize}
	\xitem[hIp-27] $p^*\forces{\poP}{\dot{g}(
	\dot{x}'_{\alpha_0}\ctentenc\dot{x}'_{\alpha_{n^*-1}})
	\not\in V[\dot{G}_R]}$.
\end{xitemize}

By thinning out $S^{**}$ further, if necessary, we may assume that 
$\supp(p^*)\cap\supp(p_\alpha)\subseteq R$ for all $\alpha\in S^{**}$. 
%% Let 
%% $X=\supp(p^*)\setminus\bigcup_{i<n^*}\supp(p_{\alpha_i})$. 
For 
$i<n^*$, let $\dot{S}_i$ be a $\poP$-name \st\ 
%\[ \forces{\poP}{\dot{D}_i=\setof{\alpha\in S^{**}}{%
%	\tilde{j}_{\alpha_i,\alpha}(p^*\restr (X\cup\supp(p_{\alpha_i})))
%		\in\dot{g} }}.
%\]\noindent
\begin{xitemize}
	\xitem[hIp-23] $ \forces{\poP}{\dot{S}_i=\setof{\alpha\in S^{**}}{%
	\tilde{j}_{\alpha_i,\alpha}(p^*)
		\in\dot{G} }}$.
\end{xitemize}
By \Lemmaof{stat-in-G}, we have 
$p^*\forces{\poP}{\dot{S}_i\xmbox{ is a stationary subset of }\kappa}$ for all 
$i<n^*$. Note that we have $\tilde{j}_{\alpha_i,\alpha}(p^*)\leq_\poP p_\alpha$ by 
\xitemof{hIp-26-0} and \xitemof{hIp-18},\,\assertof{ii}. 

%% Let $\dot{S}'_i$, $i<n^*$ be $\poP$-names \st\ 
%% \begin{xitemize}
%% 	\xitem[hIp-28] $
%% 	\begin{array}[t]{r@{}l}
%% 		p^*\forces{\poP}{{}&\dot{S}'_i\subseteq\dot{S}_i\mbox{ and }
%% 			\dot{S}'_i\mbox{ is stationary for all }i<n^*\mbox{, and}\\[\jot]		
%% 			&\dot{S}'_i,\,i<n^*\mbox{ are pairwise disjoint}}
%% 	\end{array}
%% $
%% \end{xitemize}
\begin{Claim}\label{Claim4.3.3}\mbox{}\vspace{-1.2\baselineskip}\\
	\begin{xitemize}
	\item[] 
		$\begin{array}[t]{r@{}l}
		p^*\forces{\poP}{\forall\beta_0\cdots\forall\beta_{n^*-1}\ 
			\Big(\,&\pairof{\beta_0\ctentenc\beta_{n^*-1}}\in
			{\injprod{\dot{S}_0\ctentenc \dot{S}_{n^*-1}}}
			\\[\jot]
				&\,\rightarrow
			\dot{g}(\pairof{\dot{f}(\beta_0)\ctentenc\dot{f}(\beta_{n^*-1})})
			\not\in V[\dot{G}_R]\Big)}. 
	\end{array}
	$
	\end{xitemize}
	
\end{Claim}
\prfofClaim
Suppose that $q\leq_\poP p^*$ and 
$q\forces{\poP}{\pairof{\beta_0\ctentenc \beta_{n^*-1}}
	\in\injprod{\dot{S}_0\ctentenc \dot{S}_{n^*-1}}}$. Then, by 
\xitemof{hIp-23}, $q\forces{\poP}{\tilde{j}_{\alpha_i,\beta_i}(p^*)\in \dot{G}}$ 
for $i<n^*$. It follows that
\begin{xitemize}
	\xitem[hIp-24]
	$q\forces{\poP}{\tilde{j}_{\alpha_i,\beta_i}(p^*)\restr d_{\beta_i}\in\dot{G}}$
	for $i<n^*$.  
\end{xitemize}
Let 
\[\tilde{j}=\tilde{j}_{\alpha_0,\beta_0}\circ \tilde{j}_{\alpha_1,\beta_1}\circ\cdots\circ
	\tilde{j}_{\alpha_{n^*-1},\beta_{n^*-1}}.
\]\noindent
Then 
\begin{xitemize}
	\xitem[oplus] $\tilde{j}(p^*)=
	p^*\restr (\delta\setminus \bigcup_{i<n^*}d_{\alpha_i})
	\cup \tilde{j}_{\alpha_0,\beta_0}(p^*)\restr d_{\beta_0}
	\cup\cdots\cup \tilde{j}_{\alpha_{n^*-1},\beta_{n^*-1}}(p^*)\restr d_{\beta_{n^*-1}}$
\end{xitemize}
by \xitemof{hIp-18}. 
Hence
\begin{xitemize}
	\xitem[@-1] $q\forces{\poP}{\tilde{j}(p^*)\in\dot{G}}$ 
\end{xitemize}
by $q\leq_\poP p^*$ and \xitemof{hIp-24} and \xitemof{oplus}. 
By definition of $\tilde{j}$ and $q''$, and by \xitemof{hIp-18}, we have 
\begin{xitemize}
	\xitem[@-2]
	$\tilde{j}(p^*)\leq_\poP \tilde{j}(q'')\leq_\poP \tilde{j}_{\alpha_i,\beta_i}(p_{\alpha_i})=p_{\beta_i}$ 
	for $i<n^*$ and
	\xitem[@-3] $\tilde{j}(\dot{x}'_{\alpha_i})=\dot{x}'_{\beta_i}$ for $i<n^*$ 
	by \xitemof{hIp-22}.
\end{xitemize}
Hence by \xitemof{hIp-27}
\[ q\forces{\poP}{\dot{g}(\pairof{\dot{x}'_{\beta_0}\ctentenc 
		\dot{x}'_{\beta_{n^*-1}}})\not\in V[\dot{G}_R]}. 
\]\noindent
By \xitemof{hIp-16}, \xitemof{@-1} and \xitemof{@-2}, it follows that 
\[ q\forces{\poP}{\dot{f}(\beta_i)=\dot{x}'_{\beta_i}}
\]\noindent
for $i<n^*$.
Hence $q\forces{\poP}{
	\dot{g}(\pairof{\dot{f}(\beta_0)\ctentenc \dot{f}(\beta_{n^*-1})})
	\not\in V[G_R]}$.\nopagebreak
\qedofClaim\qedskip\\

To show that $p^*\forces{\poP}{\mbox{\assertof{i1} holds}}$, suppose that 
$q\leq_\poP p^*$ and 
$\pairof{\beta_0\ctentenc\beta_{n^*-1}}$,
$\pairof{\gamma_0\ctentenc\gamma_{n^*-1}}\in\injprod{S^{**}}^{n^*}$ are \st
\begin{xitemize}
	\xitem[@-4] $\ssetof{\beta_0\ctentenc\beta_{n^*-1}}
	\ \cap\ \ssetof{\gamma_0\ctentenc\gamma_{n^*-1}}=\emptyset$ and
	\xitem[@-5] $q\forces{\poP}{\pairof{\beta_0\ctentenc\beta_{n^*-1}},\,
	\pairof{\gamma_0\ctentenc\gamma_{n^*-1}}
	\in\injprod{\dot{S}_0\ctentenc\dot{S}_{n^*-1}}}$. 
\end{xitemize}
Note that it is enough to consider $\pairof{\beta_0\ctentenc\beta_{n^*-1}}$,
$\pairof{\gamma_0\ctentenc\gamma_{n^*-1}}\in\injprod{S^{**}}^{n^*}$ with 
\xitemof{@-4} since we can thin out $\dot{S}^G_0$\ctentenc\, 
$\dot{S}^G_{n^*-1}$ afterwards if necessary so that they are pairwise disjoint. 

By the remark after \xitemof{hIp-23}, 
we may assume that 
\begin{xitemize}
\item[] $q\leq_\poP p^*\cup p_{\beta_0}\cup\cdots\cup p_{\beta_{n^*-1}}
	\cup p_{\gamma_0}\cup\cdots\cup p_{\gamma_{n^*-1}}$. 
\end{xitemize}
By \Lemmaof{homogeneous},\,\assertof{b}, there are $\poP$-names $\dot{y}$,  
$\dot{z}$ \st 
\begin{xitemize}
	\item[] $\supp(\dot{y})\cap\supp(\dot{z})\subseteq R$ and 
	\item[] $\begin{array}[t]{r@{}l}
		p^*\cup p_{\beta_0}\cup\Ctenten\cup p_{\beta_{n^*-1}}
		\cup p_{\gamma_0}\cup\Ctenten\cup p_{\gamma_{n^*-1}}\forces{\poP}{{}&
			\dot{g}(\pairof{\dot{f}(\beta_0)\ctentenc \dot{f}(\beta_{n^*-1})})
			=\dot{y}\\[\jot]
			&\land\ 
			\dot{g}(\pairof{\dot{f}(\gamma_0)\ctentenc \dot{f}(\gamma_{n^*-1})})=\dot{z}}.		
	\end{array}
	$ 
\end{xitemize}
By \Claimof{Claim4.3.3} and \Lemmaof{injectivity-lemma}, it follows that 
\begin{xitemize}
\item[] $
	\begin{array}[t]{r@{}l}
		q\leq_\poP p^*\cup p_{\beta_0}\cup\Ctenten\cup p_{\beta_{n^*-1}}
	\cup p_{\gamma_0}\cup\Ctenten\cup p_{\gamma_{n^*-1}}\forces{\poP}{{}&
	\dot{g}(\pairof{\dot{f}(\beta_0}\ctentenc \dot{f}(\beta_{n^*-1}))\not=\\[\jot]
	&\dot{g}(\pairof{\dot{f}(\gamma_0)\ctentenc \dot{f}{\gamma_{n^*-1})})}}. 
	\end{array}$
\end{xitemize}
Since $q$ as above may be chosen below arbitrary $r\leq_\poP p^*$, it 
follows that 
\begin{xitemize}
\item[] $p^*\forces{\poP}{\mbox{\assertof{i1} holds}}$.
\end{xitemize}
\qedof{\Thmof{forcing-IP*}}
\setcounter{Thm}{\value{oldthm}}

\begin{Cor}
	\label{fs-prod-IP}
	\assert{a} Assume \CH\ and $\poP=\Fn(\mu,2)$ for some cardinal $\mu$. Then 
	$\forces{\poP}{\IP(\aleph_2,\aleph_1)}$ holds. \smallskip

	\assert{b} Assume \GCH\ and $\poP=\Fn(\mu,2)$ for some cardinal $\mu$. Then 
	$\forces{\poP}{\IP(\kappa^+,\aleph_1)}$ holds for every uncountable $\kappa$ of 
	uncountable cofinality and $\forces{\poP}{\IP(\lambda,\aleph_1)}$ for every 
	inaccessible $\lambda$. \smallskip

	\assert{c} Assume \CH\ and $\poP$ is a finite support product of copies 
	of a productively c.c.c.\ \po\ of cardinality $\aleph_1$. Then 
	$\forces{\poP}{\IP(\aleph_2,\aleph_1)}$ holds. In particular, we have 
	$\forces{\poP}{\HP(\aleph_2)}$. \smallskip

	\assert{d} Assume \GCH\ and $\poP$ is a finite support product of copies 
	of a productively c.c.c.\ \po\ of cardinality $\aleph_1$. Then 
	$\forces{\poP}{\IP(\kappa^+,\aleph_1)}$ holds for every uncountable $\kappa$ of 
	uncountable cofinality and $\forces{\poP}{\IP(\lambda,\aleph_1)}$ for every
	inaccessible $\lambda$. \smallskip

	\assert{e} Assume \CH\ and $\poP$ is a countable support product of 
	copies of a proper \po\ of cardinality $\aleph_1$ \st\ its product is 
	also proper. Then $\forces{\poP}{\IP(\aleph_2,\aleph_2)}$ holds. In 
	particular, we have  
	$\forces{\poP}{\HP(\aleph_2)}$. \smallskip

	\assert{f} Assume \GCH\ and $\poP$ is a countable support product of copies 
	of a proper \po\ of cardinality $\aleph_1$ \st\ its product is 
	also proper. Then 
	$\forces{\poP}{\IP(\kappa^+,\aleph_2)}$ holds for every uncountable $\kappa$ of 
	uncountable cofinality and $\forces{\poP}{\IP(\lambda,\aleph_2)}$ for every
	inaccessible $\lambda$. 
\end{Cor}
Note that countable support products of Sacks or Prikry-Silver forcing are 
instances of  
\assertof{e} and \assertof{f} above.\qedskip\\
\prf Under \CH, $\omega_1=2^{<\omega_1}$ and 
$\omega_2$ is ${<}\omega_1$-inaccessible. 
In \assertof{a} and \assertof{b}, 
$\poP$ is forcing equivalent to a finite support product of copies of the 
countable \po\ $\Fn(\omega,2)$.
Clearly $\poP$'s in all of 
\assertof{a} $\sim$ \assertof{f} are homogeneous; $\poP$'s in \assertof{a}
$\sim$ \assertof{d} satisfy the c.c.c.\ and  
hence they are proper. Thus we can apply \Thmof{forcing-IP*}. The second parts of 
\assertof{c} and  \assertof{e} follow from 
\Thmof{IP-implies-HP}.  
\qedofCor
\qedskip
\\
Results similar to 
\Thmof{fircing-IP*} and \Corabove\ also hold for partial orderings with 
product-like structure as those considered in \cite{FuShSo}. Thus, we can 
prove e.g.\ that $\IP(\aleph_2,\aleph_2)$ together with clubsuit principle 
is consistent.

In \cite{fu-ge-sh-so} it is shown that, if we start from a model $V$ which is 
obtained by adding a dominating real to a model of \GCH\ $+$ Chang's conjecture for
$\aleph_\omega$, i.e.\
$(\aleph_{\omega+1},\aleph_\omega)\pfeil(\aleph_1,\aleph_0)$, then adding 
more than $\aleph_{\omega+1}$ Cohen reals forces $\neg\WFN$. Since $V$ 
satisfies \GCH, $\IP(\kappa,\aleph_1)$ is forced for every $\kappa\geq\aleph_2$ which is 
not a successor of a singular cardinal of cofinality $\omega$ by adding 
any number of Cohen reals by \Corof{fs-prod-IP}. In 
particular:

\begin{Cor} \label{chang-conj-then-hechler}
	Suppose that Chang's conjecture for $\aleph_\omega$ is 
	consistent. Then so is $\IP(\aleph_2,\aleph_1)$ $\land$
	$\boundingno^*=\aleph_1$ $\land$ $\neg\WFN$. \qed
\end{Cor}
\section{Models of\ \ $\IP(\aleph_2,\aleph_2)$ $\land$ $\neg\IP(\aleph_2,\aleph_1)$}
\label{IP-aleph2-aleph1}
Recall that Prikry-Silver forcing $\poS$ is the forcing with partial 
functions with co-infinite domain, that is
\begin{xitemize}
\item[] $\poS=\setof{f}{\mapping{f}{D}{2},\ D\subseteq \omega,\ 
	\cardof{\omega\setminus D}=\aleph_0 }$
\end{xitemize}
with the ordering
\begin{xitemize}
\item[] $f\leq_\poS g$\ \ $\Leftrightarrow$\ \ $f\supseteq g$
\end{xitemize}
for $f$, $g\in\poS$. 

A $(V,\poS)$-generic filter $G$ gives rise to the function
$\mapping{s_G=\bigcup G}{\omega}{2}$ which is often called a Prikry-Silver real. 

For $f\in\poS$ let $\codom(f)=\omega\setminus\dom(f)$. 

It is easy to check that Prikry-Silver forcing $\poS$ as well as its 
countable support products $\poS^I$ over any index set $I$ satisfy the 
Axiom A. Hence they are all proper. 

Note that, by definition of $\leq_\poS$, we have:
\begin{xitemize}
	\xitem[silver-0] $f$, $g\in\poS$ are 
	incompatible if $\cardof{\codom(f)\cap\codom(g)}<\aleph_0$.
	\xitem[silver-1] For any $\pairof{f_0,f_1}\in\poS^2$, there is
	$\pairof{g_0,g_1}\leq_{\poS^2}\pairof{f_0,f_1}$ \st\
	$\cardof{\codom(g_0)\cap\codom(g_1)}<\aleph_0$.  
\end{xitemize}
\begin{Lemma}
	\label{poS}
	For any $f\in\poS$ and $\pairof{g^n_0,g^n_1}\in\poS^2$, $n\in\omega$ \st\ 
	$\cardof{\codom(g^n_0)\cap\codom(g^n_1)}<\aleph_0$ there is
	$g\leq_\poS f$ \st\ $\pairof{g,g}$ is incompatible with all
	$\pairof{g^n_0,g^n_1}$, $n\in\omega$. 
\end{Lemma}
\prf
Construct $i_n\in 2$, $n\in\omega$ and $A\subseteq\codom(f)$ recursively so 
that 
\begin{xitemize}
\item[] $\cardof{\codom(f)\cap\bigcap_{k\leq n}\dom(g^k_{i_k})}=\aleph_0$ and 
\item[] $\cardof{A\cap\codom(g^n_{i_n})}<\aleph_0$ for all $n\in\omega$. 
\end{xitemize}
Then any extension $g$ of $f$ on $\omega\setminus A$ 
will do.
\qedofLemma
\qedskip

Working in $V=L$, we can construct recursively a maximal antichain 
$\setof{\pairof{g^\alpha_0,g^\alpha_1}}{\alpha<\omega_1}$ in $\poS^2$ \st\ 
\begin{xitemize}
	\xitem[] $\cardof{\codom(g^\alpha_0)\cap\codom(g^\alpha_1)}<\aleph_0$ for 
	all $\alpha<\omega_1$.
\end{xitemize}
Note that each step of the recursive construction is possible by 
\xitemof{silver-1} and \xitemof{silver-1}. Furthermore by 
choosing $\pairof{g^\alpha_0,g^\alpha_1}$ in each step of the construction 
according to the 
$\Sigma^1_2$-well ordering of the reals (which exists because of $V=L$), we can make 
$\setof{\pairof{g^\alpha_0,g^\alpha_1}}{\alpha<\omega_1}$ 
a $\Sigma^1_2$-set 
(actually we can even choose such a maximal antichain as a $\Pi^1_1$-set 
arguing similarly to \cite{miller}). 

Let $\mapping{\varphi}{\poS^2}{\fnsp{\omega}{2}}$ be a Borel bijection and 
let $\mapping{g}{\injprod{\fnsp{\omega}{2}}^{<\omega}}{\fnsp{\omega}{2}}$ be 
defined by 
\begin{xitemize}
	\xitem[silver-2] $g(\pairof{x_0\ctentenc x_{n-1}})= {}\left\{\  
	\begin{array}[c]{@{}l@{\,\ }l}
		\varphi(g^{\alpha^*}_0,g^{\alpha^*}_1) 
		&\mbox{; if }n=2\mbox{, there is }\alpha<\omega_2\mbox{ \st}\\[-\jot]
		&\mbox{\phantom{; }}
		x_0\supseteq g^\alpha_0,\ x_1\supseteq g^\alpha_1\mbox{ and }
		\alpha^*\mbox{ is minimal }\\[-\jot]
		&\mbox{\phantom{; }}\mbox{among such }\alpha\mbox{'s}\\
		0 &\mbox{; otherwise.}
	\end{array}\right.
	$
\end{xitemize}
It is easy to check that $g$ is a $\Delta^1_3$-set. 
\begin{Thm}
	\label{prikry-silver}
Assume $V=L$. Then we have 
\begin{xitemize}
\item[] $\forces{\poS^{\omega_2}}{\IP(\aleph_2,\aleph_2)\mbox{ and }
	\neg\IP(\aleph_2,\aleph_1)}$.
\end{xitemize}
\end{Thm}
\prf $\forces{\poS^{\omega_2}}{\IP(\aleph_2,\aleph_2)}$ follows from 
\Corof{fs-prod-IP},\,\assertof{e}. 

To show that $\forces{\poS^{\omega_2}}{\neg\IP(\aleph_2,\aleph_1)}$, 
let $G$ be a $(V,\poS^{\omega_2})$-generic filter. Working in $L[G]$, let
$s_\beta$ be the $\beta$'th Prikry-Silver real added by $G$ for
$\beta<\omega_2$. Let
$\mapping{f}{\omega_2}{\fnsp{\omega}{2}}$ be defined by 
\begin{xitemize}
	\xitem[] $f(\beta)=s_\beta$ for $\beta<\omega_2$
\end{xitemize}
and let $\mapping{g}{\injprod{\fnsp{\omega}{2}}^{<\omega}}{\fnsp{\omega}{2}}$ be 
the mapping as in
\xitemof{silver-2}, or more precisely, let $g$ be the mapping (in $L[G]$) defined by 
the $\Delta^1_3$ definition corresponding to \xitemof{silver-2}.

We show that $f$ and $g$ build a counter-example 
to $\IP(\aleph_2,\aleph_1)$. 

Since $\cardof{\range(g)}\leq\aleph_1$, \assertof{i1} clearly fails for 
these $f$ and $g$. 
Hence we are done by showing that $f$ and $g$ do not satisfy 
\assertof{i0}. 

Assume, for a contradiction, that $f$ and $g$ satisfy \assertof{i0}. 
Returning to $L$, 
let
$\dot{f}$, $\dot{g}$, $\dot{s}_\beta$, $\beta<\omega_2$ etc.\ be 
$\poS^{\omega_2}$-names of $f$, $g$, $s_\beta$, $\beta<\omega_2$ etc.\ 
respectively. In particular, we can choose $\dot{f}$ \st
\begin{xitemize}
	\xitem[silver-2-a] $\forces{\poS^{\omega_2}}{\dot{f}(\beta)=\dot{s}_\beta}$ for all
	$\beta<\omega_2$. 
\end{xitemize}

Since $\poS^{\omega_2}$ is proper, there are $p\in G$,
$\poS^{\omega_2}$-name $\dot{S}$ and a countable set $Z$ (in $L$) \st\ 
\begin{xitemize}
	\xitem[silver-2-0] $p\forces{\poS^{\omega_2}}{\dot{S}\subseteq\omega_2
	\mbox{ is stationary and }
	\dot{g}\imageof\injprod{\dot{f}\imageof\dot{S}}^2\subseteq Z}$. 
\end{xitemize}

Let $U=\setof{\beta<\omega_2}{\mbox{there is }p'\leq_{\poS^{\omega_2}}p
	\mbox{ \st\ }p'\forces{\poS^{\omega_2}}{\beta\in\dot{S}}}$. 
Then $U$ is a stationary subset of $\omega_2$. For each $\beta\in U$, let 
$p_\beta\leq_{\poS^{\omega_2}}p$ be \st\
$p_\beta\forces{\poS^{\omega_2}}{\beta\in\dot{S}}$  and $\beta\in\supp(p_\beta)$. 

By $\Delta$-System Lemma and \CH, there is $U^*\in[U]^{\aleph_2}$ \st\ 
\begin{xitemize}
	\xitem[silver-3] $\supp(p_\beta)$, $\beta\in U^*$ form a $\Delta$-system 
	with root $R$ which is an initial segment of all of $\supp(p_\beta)$,
	$\beta\in U^*$; 
	\xitem[silver-3-0] $\sup R<\min U^*$; 
	\xitem[silver-4] $p_\beta\restr R$, $\beta\in U^*$ are all the same;  
	and 
	\xitem[silver-5] $p_\beta(\beta)$, $\beta\in U^*$ are all the same, say
	$h\in\poS$.  
\end{xitemize}
Note that $p_\beta$, $\beta\in U^*$ are compatible by \xitemof{silver-3} 
and \xitemof{silver-4}. 

Let 
\begin{xitemize}
	\xitem[] $X=\varphi^{-1}(Z)$.
\end{xitemize}
By \Lemmaof{poS}, there is a $k\leq_\poS h$ \st\ $\pairof{k,k}$ is 
incompatible with all $\pairof{g^\alpha_0,g^\alpha_1}$ from the countable 
set $X$. 

Fix two distinct $\beta$, $\gamma\in U^*$ and let
$q\leq_{\poS^{\omega_2}}p_\beta$, $p_\gamma$ be defined by 
$\dom(q)=\dom(p_\beta)\cup\dom(p_\gamma)$ and
\begin{xitemize}
	\xitem[] $q(\delta)=\left\{\,
	\begin{array}[c]{@{}ll}
		p_\beta(\delta) 
		&\mbox{; if }\delta\in\supp(p_\beta)\setminus\ssetof{\beta}\\[\jot]
		p_\gamma(\delta)
		&\mbox{; else if }\delta\in\supp(p_\gamma)\setminus\ssetof{\gamma}\\[\jot]
		k
		&\mbox{; else if }\delta=\beta\mbox{ or }\delta=\gamma
	\end{array}\right.
	$
\end{xitemize}
for $\delta\in\dom(q)$. 

By $q\leq_{\poS^{\omega_2}}p_\beta$, $p_\gamma$, we have
$q\forces{\poS^{\omega_2}}{\beta,\,\gamma\in\dot{S}}$. Thus the following 
claim yields a contradiction to \xitemof{silver-2-0}: 
\begin{Claim}
	$q\forces{\poS^{\omega_2}}{
		\dot{g}(\pairof{\dot{f}(\beta),\dot{f}(\gamma)})\not\in Z}$.
\end{Claim}
\prfofClaim 
By \xitemof{silver-2-a}, we have to show
$q\forces{\poS^{\omega_2}}{\dot{g}(\pairof{\dot{s}_\beta,\dot{s}_\gamma})\not\in Z}$. 

First, we show that
$q\forces{\poS^{\omega_2}}{\dot{g}(\pairof{\dot{s}_\beta,\dot{s}_\gamma})\not=0}$. 
Note that, 
by the complete embedding 
$\poS^2\ni\pairof{g_0,g_1}\mapsto
\ssetof{\pairof{\beta,g_0},\pairof{\gamma,g_1}}\in\poS^{\ssetof{\beta,\gamma}}
\circleq\poS^{\omega_2}$ we have:  
\begin{xitemize}
	\item[] $\setof{\ssetof{\pairof{\beta,g^\alpha_0},\pairof{\gamma,g^\alpha_1}}}{
		\alpha<\omega_1}$ is a maximal antichain in $\poS^{\omega_2}$. 
\end{xitemize}

For any $r\leq_{\poS^{\omega_2}}q$, let $\alpha^*<\omega_1$ be \st\ $r$ 
and $\ssetof{\pairof{\beta,g^{\alpha^*}_0},\pairof{\gamma,g^{\alpha^*}_1}}$ 
are compatible. Let $s\leq_{\poS^{\omega_2}}r$,
$\ssetof{\pairof{\beta,g^{\alpha^*}_0},\pairof{\gamma,g^{\alpha^*}_1}}$. 
Then we have 
\begin{xitemize}
	\item[] $s\forces{\poS^{\omega_2}}{
	\dot{s}_\beta\supseteq g^{\alpha^*}_0,\ 
	\dot{s}_\gamma\supseteq g^{\alpha^*}_1}$. 
\end{xitemize}
Hence, by 
\xitemof{silver-2}, it follows that
$s\forces{\poS^{\omega_2}}{\dot{g}(\pairof{\dot{s}_\beta,\dot{s}_\gamma})\not=0}$. 

Now, suppose, for contradiction, that there is $r\leq_{\poQ^{\omega_2}}q$ 
\st\ 
\begin{xitemize}
\item[] 
	$r\forces{\poS^{\omega_2}}{\dot{g}(\pairof{\dot{s}_\beta,\dot{s}_\gamma})\in Z}$. 
\end{xitemize}
Then, by the first part of the proof,  there are $s\leq_{\poS^{\omega_1}}r$ and
$\pairof{g^\alpha_0,g^\alpha_1}\in X$ \st\ 
$s\forces{
	\poS^{\omega_2}}{\dot{s}_\beta\supseteq g^\alpha_0\mbox{ and }
	\dot{s}_\gamma\supseteq g^\alpha_1}$. In particular $s(\beta)$ and $s(\gamma)$ 
are compatible with $g^\alpha_0$ and $g^\alpha_1$, respectively. 
Since
$r\leq_{\poS^{\omega_2}}q\leq_{\poS^{\omega_2}}
	\ssetof{\pairof{\beta,k},\pairof{\gamma,k}}$, 
it follows that $k$ is compatible with both of $g^\alpha_0$ and $g^\alpha_1$. This 
is a contradiction to the choice of $k$.\qedofClaim\\
\qedofThm
\qedskip
\\
We can prove a Lemma similar to \Lemmaof{poS} for $omega$ product of Sacks 
forcing. Thus, by a similar argument as above, we can also prove that
$\IP(\aleph_2,\aleph_1)$ fails in 
a generic extension by countable support side-by-side product of Sacks 
forcing. 

\section{The Consistency of 
	$\boundingno^*=\aleph_2$ $\land$ $\donum=\aleph_1$}
%% \subsection{Forcing $\HP(\aleph_2)$ $+$ $\boundingno^*=\aleph_2$ from the 
%% 	assumption (A)}
\label{Omega}
In the following we shall refer by (A) the assertion that there is a 
structure $\pairof{(\omega_2)^2,A,\calF}$ with the properties 
\xitemof{Omega-0} $\sim$ \xitemof{Omega-4} below. Recall that a mapping 
$\mapping{f}{X}{X}$ is called an involution if it is a bijection exchanging 
(some) pairs of elements of $X$, that is, if $f\circ f=id_X$ holds.  
\begin{xitemize}
	\xitem[Omega-0] 
	$\omega_2\times\omega_2\supseteq A
		\supseteq\setof{\pairof{\alpha,\beta}\in\omega_2\times\omega_2}{\beta<\alpha}$; 
	\xitem[Omega-1]
	For any $C\in[\omega_2]^{\aleph_0}$ there is 
	an $X\in[\omega_2]^{\aleph_2}$ \st\ 
	$(C\times X)\cap A=\emptyset$;
	\xitem[Omega-2]
	For all $\pairof{\phi,\psi}\in\calF$, $\phi$ and $\psi$ are involutions 
	on $\omega_2$;
%% 	\xitem[Omega-2A] If $\pairof{\phi,\psi}\in\calF$ and
%% 	$\pairof{\phi',\psi'}\in\calF$ then
%% 	$\pairof{\phi\circ\phi',\psi\circ\phi'}\in\calF$; 
	\xitem[Omega-3]
	For each $\pairof{\phi,\psi}\in\calF$ and for all
	$\pairof{\alpha,\beta}\in \omega_2\times\omega_2$,  
	we have $\pairof{\alpha,\beta}\in A$ if and only if
	$\pairof{\phi(\alpha),\psi(\beta)}\in A$; 
	\xitem[Omega-4] For any stationary $S\subseteq E^{\omega_2}_{\omega_1}$ 
	and any $A_\zeta$, $B_\zeta\in[\omega_2]^{\aleph_0}$ for $\zeta\in S$, 
	there is a stationary $T\subseteq S$ \st, for any $n\in\omega$, if 
	$\zeta_i$, $\eta_i\in T$ for $i\in n$ are pairwise distinct ($2n$ 
	elements) then there is $\pairof{\phi,\psi}\in\calF$ \st\ 
	\begin{xxitemize}
	\item[] $\phi\imageof A_{\zeta_i}=A_{\eta_i}$, 
		$\psi\imageof B_{\zeta_i}=B_{\eta_i}$\,; and
	\item[] $\mapping{\phi\restr A_{\zeta_i}}{A_{\zeta_i}}{A_{\eta_i}}$, 
		$\mapping{\psi\restr B_{\zeta_i}}{B_{\zeta_i}}{B_{\eta_i}}$ are order 
		isomorphisms
	\end{xxitemize}
	for all $i\in n$. 
%% 	\xitem[Omega-4a] 
%% 	For any $X_0$, $Y_0\in[\omega_2]^{\aleph_0}$, there are $X$, 
%% $Y\in[\omega_2]^{\aleph_0}$ with $X_0\subseteq X$ and $Y_0\subseteq Y$ \st, 
%% 	for any $X'$, $Y'\in[\omega_2]^{\aleph_0}$ with $X\cap X'=\emptyset$ and
%% 	$Y\cap Y'=\emptyset$, there are $X''$, $Y''\in[\omega_2]^{\aleph_0}$ and 
%% $\pairof{\phi,\psi}\in\calF$ \st\ 
%% 	\begin{xxitemize}
%% 		\item[] $X''\cap(X\cup X')=\emptyset$, $Y''\cap(Y\cup Y')=\emptyset$\,;
%% 		\item[] $\phi\restr X=id_X$, $\psi\restr Y=id_Y$\,;  
%% 		\item[] $\phi\imageof X'=X''$, $\psi\imageof Y'= Y''$\,; 
%% 		\item[] $\mapping{\phi\restr X'}{X'}{X''}$ and 
%% 			$\mapping{\psi\restr Y'}{Y'}{Y''}$ are order isomorphisms.
%% 	\end{xxitemize}
\end{xitemize}

The consistency of (A) together with \CH\ over \ZFC\ is proved in the next section. 
Below, we will prove the consistency of
$\continuum=\boundingno^*=\aleph_2$ $\land$ $\donum=\aleph_1$ $\land$
$\neg\Cs(\aleph_2)$ by constructing a model of this combination of 
assertions starting  
from a model of (A) and \CH. 

Let us begin with introducing some notation for the forcing construction we 
use in the proof. 

For a cardinal $\kappa$, a 
sequence $\bar{f}=\seqof{f_\xi}{\xi<\kappa}$ in $\fnsp{\omega}{\omega}$ 
and $X\subseteq\kappa$, let $\poD_{\bar{f},X}$ be the canonical \po\ adding 
an element of  $\fnsp{\omega}{\omega}$ dominating $\setof{f_\xi}{\xi\in X}$. That is
\begin{xitemize}
	\xitem[poD-0] $\poD_{\bar{f},X}
	=\setof{\pairof{s,F}}{s\in\fnsp{\omega{>}}{\omega},\,F\in[\kappa]^{<\aleph_0}}$
\end{xitemize}
and, for $\pairof{s,F}$, $\pairof{s',F'}\in \poD_{\bar{f},X}$\,,
\begin{xitemize}
	\xitem[poD-1] $\pairof{s',F'}\leq_{\poD_{\bar{f},X}}\pairof{s,F}\ \Leftrightarrow\ 
	\begin{array}[t]{@{}l}
		s'\supseteq s,\,F'\supseteq F,\\[\jot]
		\forall \alpha\in F\cap X\ \forall n\in\dom(s')\setminus\dom(s)\ 
		(f_\alpha(n)\leq s'(n)).
	\end{array}$
\end{xitemize}
Since any $\pairof{s,F}$, $\pairof{s',F'}\in\poD_{\bar{f},X}$ with $s=s'$ 
are compatible, we have:
\begin{Lemma}
	$\poD_{\bar{f},X}$ is $\sigma$-centered. \ifextended\else\qed\fi
\end{Lemma}
\ifextended
\prf For any $s\in\fnsp{\omega{>}}{\omega}$ and $F$, 
$F'\in[X]^{<\aleph_0}$, $\pairof{s,F}$ and $\pairof{s,F'}$ are compatible 
in $\poD_{\bar{f},X}$ by \xitemof{poD-1}. 
Hence
$\poD_{\bar{f},X}
=\bigcup_{s\in\fnsp{\omega{>}}{\omega}}
\setof{\pairof{s,F}}{F\in[X]^{<\aleph_0}}$ is a countable union of 
compatible sets. \qedofLemma
\qedskip
\fi

Note that the underlying set of $\poD_{\bar{f},X}$ does not depend on the 
sequence $\bar{f}$. 
So we shall denote this set with $\poD_X$. 
Actually $\poD_X$ as a set does not depend on $X$ either. Nevertheless we 
shall add the suffix $X$ so that we can distinguish $\poD$'s by their 
intended function. 
%% For convenience, we shall regard $\poD_X$ as a partial ordering with the 
%% largest element $\pairof{\emptyset,\emptyset}$ but, except the order 
%% relation making $\pairof{\emptyset,\emptyset}$ the largest element, without 
%% any other ordering  
%% among elements of $\poD_X\setminus\ssetof{\pairof{\emptyset,\emptyset}}$. 
%% With this partial ordering, we can talk about the finite support side by 
%% side product 
%% $\prod^{fin}_{\alpha<\kappa}\poD_{X_\alpha}$ in the sense of 
%% \xitemof{prod-0} 
%% for some fixed sequence
%% $\seqof{X_\alpha}{\alpha<\kappa}$ of subsets of $\kappa$.  %% Here, a finite 
%% %% support side by side product $\prod^{fin}_{\alpha<\kappa}\poP_\alpha$ for 
%% %% \pos\ $\poP_\alpha$, $\alpha<\kappa$ is considered as consisting of finite 
%% %% partial function $p$ from $\kappa$ to $\bigcup_{\alpha<\kappa}\poP_\alpha$ \st\ 
%% %% for each $\alpha\in\dom(p)$ $p(\alpha)\in\poP_\alpha$; and equiped with the 
%% %% canonical ordering induced from the orderings of $\poP_\alpha$,
%% %% $\alpha<\kappa$. 

Note also that, as a set, 
$\prod^{fin}_{\alpha<\kappa}\poD_{\bar{f},X_\alpha}$ for 
any $\kappa$-sequence $\bar{f}$ of reals is the same: we shall denote this 
set by $\prod^{fin}_{\alpha<\kappa}\poD_{X_\alpha}$.

If $d\in\poD_X$ and $d=\pairof{s,F}$ then we shall write $s^d$ 
and $F^d$  to denote these $s$ and $F$ respectively. 

In the following we assume that a 
sequence $\bar{X}=\seqof{X_\alpha}{\alpha<\kappa}$ of nonempty subsets of $\kappa$ 
is fixed. 
Let 
\begin{xitemize}
	\xitem[] $\poQ_{\bar{X}}
	=\poC_\kappa\ast\prod^{fin}_{\alpha<\kappa}\poD_{\dot{\bar{f}},X_\alpha}$
\end{xitemize}
where $\poC_\kappa=\Fn(\kappa\times\omega,\omega)$ and $\dot{\bar{f}}$ 
denotes the $\poC_\kappa$-name of the sequence of Cohen reals
($\in\fnsp{\omega}{\omega}$) 
of length $\kappa$ added by $\poC_\kappa$. Thus, if $G$ is a 
$(V,\poC_\kappa)$-generic set and $c_\alpha$ is the $\alpha$'th element of
$\dot{\bar{f}}^G$, then $c_\alpha(n)=m$ if and only if there is a condition 
$c\in G$ \st\ $\pairof{\alpha,n}\in\dom(c)$ and $c(\alpha,n)=m$. 

Let 
\begin{xitemize}
	\xitem[poQd-0] $\textstyle
\poQ^\dagger_{\bar{X}}=
\setof{\pairof{c,d}}{
	{}\begin{array}[t]{@{}l}
			c\in\poC_\kappa,\,d\in\prod^{fin}_{\alpha\in \kappa}\poD_{X_\alpha},\,
%%			\dom(d)=\pi_0(\dom(c)),
			\\[\jot]
			\bigcup_{\xi\in\dom(d)} 
			F^{d(\xi)}\times\dom(s^{d(\xi)})\ \ \subseteq\  \dom(c) }
	\end{array}
$
\end{xitemize}

%% where $\pi_0$ denotes the projection to the $0$th coordinate. 

For $\pairof{c,d}$, $\pairof{c',d'}\in\poQ^\dagger_{\bar{X}}$, 
\begin{xitemize}
	\xitem[poQd-1] 
	$
	\begin{array}[t]{@{}l}
		\pairof{c',d'}\leq_{\poQ^\dagger_{\bar{X}}}\pairof{c,d}\ \Leftrightarrow\\[\jot] 
		\begin{array}{@{\mbox{}\quad\ \ }l}
			c'\leq_{\poC_\kappa}c,\, \dom(d')\supseteq\dom(d),\\[\jot]
			\forall\alpha\in\dom(d)\ \Big(\,s^{d'(\alpha)}\supseteq s^{d(\alpha)}\ \land
			\ F^{d'(\alpha)}\supseteq F^{d(\alpha)} \land\\[\jot]
			\forall \xi\in %%s^{d'(\alpha)}\setminus 
			F^{d(\alpha)}\cap X_\alpha\ 
			\forall n\in\dom(s^{d'(\alpha)})\setminus\dom(s^{d(\alpha)})\ 
			\left(c'(\xi,n)\leq s^{d'(\alpha)}(n)\right)\Big)\,. 
		\end{array}
	\end{array}
$
\end{xitemize}
The following can be shown easily by standard arguments: 

\begin{Lemma}\label{poQdager}
%% 	\assert{b} For $X\subseteq Y\subseteq\kappa$, we have
%% 	$\poQ^\dagger_{\kappa,X}\circleq\poQ^\dagger_{\kappa,Y}$\,. \smallskip

%% 	\assert{c} 
	$\mapping{\Phi}{\poQ^\dagger_{\bar{X}}}{\poQ_{\bar{X}}}\,;$
	$\pairof{c,d}\mapsto \pairof{c,\check{d}}$ is a dense embedding of 
$\poQ^\dagger_{\bar{X}}$ into $\poQ_{\bar{X}}$.\ifextended\else\\\qed\fi
\end{Lemma}
\ifextended
\noindent
{\bf [\,The following proof will be omitted in the version for 
		publication.]} 
\smallskip\\
{\footnotesize
\prf Suppose that $\pairof{c,\dot{d}}$ is an element of $\poQ_{\bar{X}}$. We show 
that there is $\pairof{c',d'}\in\poQ^\dagger_{\bar{X}}$ \st\
$\Phi(\pairof{c',d'})\leq_{\poQ_{\bar{X}}}\pairof{c,\dot{d}}$. 
Let $c''\leq_{\poC_\kappa}c$ be \st\
$c''\forces{\poC_\kappa}{\dot{d}=\check{d}'}$ for some
$d'\in\prod^{fin}_{\alpha<\kappa}\poD_{X_\alpha}$. 
Let $c'\leq_{\poC_\kappa}c''$ be \st\ 
\[\bigcup_{\xi\in\dom(d')} 
			F^{d'(\xi)}\times\dom(s^{d'(\xi)})\ \ \subseteq\  \dom(c')\,.
\]\noindent
%% and let 
%% \[ q=p\cup\setof{\pairof{\alpha,\pairof{\emptyset,\ssetof{\xi_\alpha}}}}{
%% 		\alpha\in\pi_0(\dom(c))\setminus\dom(p)}.
%% \]\noindent
Then $\pairof{c',d'}\in\poQ^\dagger_{\bar{X}}$ 
and $\Phi(\pairof{c',d'})\leq_{\poQ_{\bar{X}}}\pairof{c,\dot{d}}$. 

Suppose now $\pairof{c,d}$, $\pairof{c',d'}\in\poQ^\dagger_{\bar{X}}$ and  
$\pairof{c,\check{d}}\leq_{\poQ_{\bar{X}}}\pairof{c',\check{d}'}$. 
By definition of the two step iteration, this means that
$c\leq_{\poC_\kappa}c'$ and
\begin{xitemize}
	\xitem[] $c'\forces{\poC_\kappa}{%
	\check{d}\leq_{\prod^{fin}_{\alpha<\kappa}\poD_{\dot{\bar{f}},X_\alpha}}
	\check{d}'}$.
\end{xitemize}
The condition \xitemabove\ is equivalent to 
\begin{xitemize}
	\xitem[xxxxx-0]
	$c'\forces{\poC_\kappa}{\dom(\check{d})\supseteq\dom(\check{d}')\ \land\ 
		\forall\alpha\in\dom(\check{d}')\ \Big(
		\check{d}(\alpha)\leq_{\poQ_{\dot{\bar{f}},X_\alpha}}\check{d}'(\alpha)\Big)}$.
\end{xitemize}
Since $\pairof{c,d}\in\poQ^\dagger_{\bar{X}}$, we have 
$\bigcup_{\xi\in\dom(d)} 
			F^{d(\xi)}\times\dom(s^{d(\xi)})\ \subseteq\  \dom(c)$. Hence,
by \xitemof{poD-1},
\xitemof{xxxxx-0} is equivalent in turn to 
\begin{xitemize}
	\xitem[] $\dom(d)\supseteq\dom(d)'\ \land\ \forall\alpha\in\dom(d)\ 
	\Big(s^{d(\alpha)}\supseteq s^{d'(\alpha)}
	\ \land\ F^{d(\alpha)}\supseteq F^{d'(\alpha)}\ \land$\\
	$\forall \xi\in F^{d(\alpha)}\,\forall n\in\dom(s^{d(\alpha)})
	\setminus\dom(s^{d'(\alpha)})\ \left(c(\xi,n)\leq s^{d(\alpha)}(n)\right)\Big)$.
\end{xitemize}
Thus we have $\pairof{c,d}\leq_{\poQ^\dagger_{\bar{X}}}\pairof{c',d'}$. 

Similarly, $\pairof{c,\check{d}}\leq_{\poQ_{\bar{X}}}\pairof{c,\check{d}'}$ 
follows from $\pairof{c,d}\leq_{\poQ^\dagger_{\bar{X}}}\pairof{c',d'}$. \qedofLemma
\qedskip}\\\fi
$\poQ_{\bar{X}}$ and $\poQ^\dagger_{\bar{X}}$ are thus forcing equivalent. 

For $p\in\poQ^\dagger_{\bar{X}}$ with $p=\pairof{c,d}$, let 
\begin{xitemize}
\item[] $\supp_0(p)=\setof{\alpha<\kappa}{\pairof{\alpha,n}\in\dom(c)
	\mbox{ for some }n\in\omega}$ and
\item[] $\supp_1(p)=\dom(d)$.
\end{xitemize}
For a $\poQ^\dagger_{\bar{X}}$-name $\dot{a}$, $\supp_0(\dot{a})$ and 
$\supp_1(\dot{a})$ are also defined in analogy to \xitemof{prod-0-a}.

In \Thmof{new-model}, we assume \CH\ $+$ (A) and let, for a structure 
$\pairof{(\omega_2)^2,A,\calF}$ as in (A), $\kappa=\omega_2$ and
$\bar{X}=\seqof{X_\alpha}{\alpha<\omega_2}$ where
$X_\alpha=\setof{\beta\in\omega_2}{\pairof{\alpha,\beta}\in A}$ for
$\alpha<\omega_2$. For such $\bar{X}$, the next lemma follows immediately 
from \xitemof{Omega-2} and \xitemof{Omega-3}.
\begin{Lemma}
	Suppose that $\pairof{(\omega_2)^2,A,\calF}$ and $\bar{X}$ are as above. 
	If $\pairof{\phi,\psi}\in\calF$, then the mapping 
	$\mapping{j_{\pairof{\phi,\psi}}}{\poQ^\dagger_{\bar{X}}}{\poQ^\dagger_{\bar{X}}}$
	defined by 
	\begin{xitemize}
	\item[] $j_{\pairof{\phi,\psi}}(\pairof{c,d})=\pairof{c',d'}$
	\end{xitemize}
	for $\pairof{c,d}\in\poQ^\dagger_{\bar{X}}$ where $c'$ and $d'$ are \st\ 
	\begin{xitemize}
	\item[] $\dom(c')=\setof{\pairof{\phi(\alpha),n}}{\pairof{\alpha,n}\in\dom(c)}$;
	\item[] $c'(\pairof{\phi(\alpha),n})=c(\pairof{\alpha,n})$ for
		$\pairof{\alpha,n}\in\dom(c)$;
	\item[] $\dom(d')=\psi\imageof\dom(d)$;
	\item[] $F^{d'(\psi(\xi))}=F^{d(\xi)}$ and $s^{d'(\psi(\xi))}=s^{d(\xi)}$ 
		for $\xi\in\dom(d)$
	\end{xitemize}
	is an automorphism on the \po\ $\poQ^\dagger_{\bar{X}}$. \qed
\end{Lemma}
Similarly to \sectionof{forcing-construction}, we shall also denote with 
$j_{\pairof{\phi,\psi}}$ the corresponding mapping on
$\poQ^\dagger_{\bar{X}}$-names. 

The following theorem together with the consistency result in 
\sectionof{forcing(A)} gives the consistency of the conjunction of 
the assertions $\continuum=\boundingno^*=\aleph_2$, $\donum=\aleph_1$ and 
$\neg\Cs(\aleph_2)$ over \ZFC. 

\begin{Thm}
	\label{new-model}
	Assume 
	\CH\ and {\rm(A)}. 
	Let $\pairof{(\omega_2)^2,A,\calF}$ be a 
	structure satisfying \xitemof{Omega-0} $\sim$ \xitemof{Omega-4} and let 
	$\bar{X}=\seqof{X_\alpha}{\alpha<\omega_2}$ where
	$X_\alpha=\setof{\beta\in\omega_2}{\pairof{\alpha,\beta}\in A}$. 
	Then 
	$\forces{\poQ^\dagger_{\bar{X}}}{
		\continuum=\boundingno^*=\aleph_2\ \land\ \donum=\aleph_1\ 
		\land\ \neg\Cs(\aleph_2)}$. 
\end{Thm}
\prf
First, we show that
$\forces{\poQ^\dagger_{\bar{X}}}{\continuum=\boundingno^*=\aleph_2}$.  
Let $G$ be a $(V,\poQ^\dagger_{\bar{X}})$-generic filter. 
Working in $V[G]$, let $\bar{f}=\seqof{c_\alpha}{\alpha<\omega_2}$ be 
the sequence of Cohen reals added by the $\poC_{\omega_2}$ part of 
$\poQ_{\bar{X}}$ and $d_\alpha$ be the Hechler type real added 
by $\poD_{\bar{f},X_\alpha}$ for $\alpha<\omega_2$. 
By \xitemof{Omega-0}, $\setof{c_\alpha}{\alpha<\gamma}$ is bounded by
$d_\gamma$ for all $\gamma<\omega_2$. On the other hand,
$\setof{c_\alpha}{\alpha<\omega_2}$  
is unbounded by \xitemof{Omega-1} and the c.c.c.\ 
of $\prod_{\alpha<\omega_2}\poD_{\bar{f}, X_\alpha}$ (in $V[\bar{f}]$).  
This shows that
$V[G]\models\aleph_2\leq\boundingno^*$. Since
$\cardof{\poQ^\dagger_{\bar{X}}}=\aleph_2$ by \CH, 
we have 
$V[G]\models\continuum\leq\aleph_2$. 

To show that $\poQ^\dagger_{\bar{X}}$ forces $\donum=\aleph_1$, suppose 
that $\dot{f}_\alpha$, $\alpha<\omega_2$ are $\poQ^\dagger_{\bar{X}}$-names 
of elements of $\fnsp{\omega}{\omega}$, $\varphi(x,y,z)$ a formula in 
$\calL_\ZF$ and $\dot{a}$ a $\poQ^\dagger_{\bar{X}}$-name of an element of 
$\fnsp{\omega}{\omega}$ \st
\begin{xitemize}
	\xitem[5.4-a-a-0] $\forces{\poQ^\dagger_{\bar{X}}}{\calH(\aleph_1)
	\models\varphi(\dot{f}_\alpha,\dot{f}_\beta,\dot{a})}$ for all
	$\alpha<\beta<\omega_2$. 
\end{xitemize}
By Maximal Principle, 
it is enough to show that there are 
$\eta_1<\eta_0<\omega_2$ \st\ 
\begin{xitemize}
\item[] $\forces{\poQ^\dagger_{\bar{X}}}{\calH(\aleph_1)
	\models\varphi(\dot{f}_{\eta_0},\dot{f}_{\eta_1},\dot{a})}$. 
\end{xitemize}

For $\xi<\omega_2$, let 
\begin{xitemize}
\item[] $A_\xi=\supp_1(\dot{f}_\xi)\cup\supp_1(\dot{a})$ and
\item[] $B_\xi=\supp_0(\dot{f}_\xi)\cup\supp_0(\dot{a})$. 
\end{xitemize}
By \CH, $\Delta$-System Lemma and \xitemof{Omega-4}, we can find a stationary 
$S\subseteq E^{\omega_2}_{\omega_1}$ \st\ 
\begin{xitemize}
	\xitem[5.4-a-a-1] $A_\xi$, $\xi\in S$ form a $\Delta$-system \st\ its 
	root is an initial segment of each of $A_\xi$, $\xi\in S$\,;\ \ 
	$B_\xi$, $\xi\in S$ form a $\Delta$-system \st\ its 
	root is an initial segment of each of $B_\xi$, $\xi\in S$\,;
	\xitem[5.4-a-0] for any distinct $\zeta_0$, $\zeta_1$, $\eta_0$, $\eta_1\in S$,
	there is $\pairof{\phi,\psi}\in\calF$ \st\
	\begin{xxitemize}
	\item[\xxitemof{5.4-a-0}{a}] $\phi\imageof A_{\zeta_i}=A_{\eta_i}$, 
		$\psi\imageof B_{\zeta_i}=B_{\eta_i}$\,; and
	\item[\xxitemof{5.4-a-0}{b}]
		$\mapping{\phi\restr A_{\zeta_i}}{A_{\zeta_i}}{A_{\eta_i}}$,  
		$\mapping{\psi\restr B_{\zeta_i}}{B_{\zeta_i}}{B_{\eta_i}}$ are order 
		isomorphisms
	\end{xxitemize}
	for $i\in 2$;
	\xitem[5.4-a-1] $j_{\pairof{\phi,\psi}}(\dot{f}_\zeta)=\dot{f}_\eta$ for 
	any distinct $\zeta$, $\eta\in S$ and $\pairof{\phi,\psi}\in\calF$ as in 
	\xitemof{5.4-a-0} with $\zeta_0=\zeta$ and $\eta_0=\eta$. 
\end{xitemize}
Note that, by \xitemof{5.4-a-a-1} and \xxitemof{5.4-a-0}{b}, we have 
\begin{xitemize}
	\xitem[] $j_{\pairof{\phi,\psi}}(\dot{a})=\dot{a}$
\end{xitemize}
for any $\pairof{\phi,\psi}$ as in \xitemof{5.4-a-0}. 

Now,  let $\zeta_0$, $\zeta_1$, $\eta_0$, $\eta_1\in S$ be four distinct 
elements of $S$ \st\ $\zeta_0<\zeta_1$ and $\eta_1<\eta_0$. 
By \xitemof{5.4-a-a-0}, we have 
\begin{xitemize}
\item[] $\forces{\poQ^\dagger_{\bar{X}}}{\calH(\aleph_1)
	\models\varphi(\dot{f}_{\zeta_0},\dot{f}_{\zeta_1},\dot{a})}$\,. 
\end{xitemize}
Hence, by mapping this situation by $j_{\pairof{\phi,\psi}}$ for 
$\pairof{\phi,\psi}\in\calF$ as in \xitemof{5.4-a-0} for these $\zeta_0$,
$\zeta_1$, $\eta_0$, $\eta_1$\,, we obtain
\begin{xitemize}
\item[] $\forces{\poQ^\dagger_{\bar{X}}}{\calH(\aleph_1)
	\models\varphi(\dot{f}_{\eta_0},\dot{f}_{\eta_1},\dot{a})}$\,. 
\end{xitemize}
Thus, $\eta_0$, $\eta_1$ above are as desired.

Finally, we show that $\poQ^\dagger_{\bar{X}}$ forces the negation of
$\Cs(\aleph_2)$. 

Let $\pairof{r^0_n,s^0_n,r^1_n,s^1_n}$, $n\in\omega$ list all quadruples of 
finite sequences
$r^0$, $s^0$, $r^1$, $s^1\in\fnsp{\omega{>}}{\omega}$ \st\ 
\begin{xitemize}
	\xitem[5.4-0] $\cardof{r^0}=\cardof{s^0}=\cardof{r^1}=\cardof{s^1}$ and 
	\xitem[5.4-1] $\pairof{r^0,s^0}\not=\pairof{r^1,s^1}$ if $\cardof{r^0}>0$. 
\end{xitemize}
We further assume that the enumeration 
$\seqof{\pairof{r^0_n,s^0_n,r^1_n,s^1_n}}{n\in\omega}$ is arranged so that 
\begin{xitemize}
%% 	\xitem[5.4-1-0] the enumeration
%% 	$\seqof{\pairof{r^0_n,s^0_n,r^1_n,s^1_n}}{n\in\omega}$ is injective;
	\xitem[5.4-2] $\cardof{r^0_n}\leq n$ for all
	$n\in\omega$.%%  and
%% 	\xitem[5.4-3] for $n\in\omega$, if $\ell<\cardof{r^0_n}$ and 
%% 	$\pairof{r^0_n\restr\ell,s^0_n\restr\ell,r^1_n\restr\ell,s^1_n\restr\ell}$
%% 	satisfies \xitemof{5.4-1}, then there is $m<n$ \st\ 
%% 	$\pairof{r^0_n\restr\ell,s^0_n\restr\ell,r^1_n\restr\ell,s^1_n\restr\ell}
%% 	=\pairof{r^0_m,s^0_m,r^1_m,s^1_m}$. 
\end{xitemize}

Now, working in $V[G]$, let $a_\alpha$, $\alpha<\omega_2$ be the subsets  
of $\omega$ defined by 
\begin{xitemize}
\item[] $n\in a_\alpha$ $\Leftrightarrow$ one of the following 
	\xitemof{5.4-4} and \xitemof{5.4-5} holds:
\end{xitemize}
\begin{xitemize}
	\xitem[5.4-4] $r^0_n\subseteq c_\alpha$,\ \
	$s^0_n\subseteq d_{\alpha+1}$,\ \ \\
	$c_\alpha(n)=0$,\ \ $d_{\alpha+1}(n)=1$,\ \ $c_\alpha(n+1)=2$\ \ and\ \ 
	$d_{\alpha+1}(n+1)=3$\,;
	\xitem[5.4-5] $r^1_n\subseteq c_\alpha$,\ \
	$s^1_n\subseteq d_{\alpha+1}$,\ \ \\
	$c_\alpha(n)=2$,\ \  $d_{\alpha+1}(n)=3$,\ \ $c_\alpha(n+1)=0$\ \ and\ \ 
	$d_{\alpha+1}(n+1)=1$\,.
\end{xitemize}
Let 
\begin{xitemize}
	\xitem[5.4-5-0] $a_{\alpha,n}=a_\alpha\setminus 
	\setof{k}{\cardof{r^0_k}<n}$\ \ 
for $\alpha<\omega_2$ and $n\in\omega$. 
\end{xitemize}

We show that the 
matrix $\seqof{a_{\alpha,n}}{\alpha<\omega_2,n\in\omega}$ together with
$T=\fnsp{2}{\omega}$ is a 
counter-example to $\Cs(\aleph_2)$. For this, it is enough to prove the 
following:
\begin{Claim}
	\label{main-claim}%%
	If $S_0$, $S_1$ are cofinal subsets of $\omega_2$, then\smallskip

	\assert{1} there exist $n<\omega$, $\alpha\in S_0$ and $\beta\in S_1$ \st\ 
	$a_{\alpha,n}\cap a_{\beta,n}=\emptyset$; and\smallskip

	\assert{2} for any $t\in\fnsp{2}{\omega}$, there are $\alpha\in S_0$ and
	$\beta\in S_1$ \st\ $a_{\alpha, t(0)}\cap a_{\beta,t(1)}\not=\emptyset$. 
\end{Claim}
\prfofClaim Working in the ground model, let $\dot{S}_0$ and $\dot{S}_1$ be
$\poQ^\dagger_{\bar{X}}$-names for the cofinal subsets of $\omega_2$. 

Let $p\in\poQ^\dagger_{\bar{X}}$. For $\alpha<\omega_2$, let
$p_\alpha\in\poQ^\dagger_{\bar{X}}$ and $\gamma_\alpha$, 
$\delta_\alpha\in\omega_2$ be \st\
\begin{xitemize}
	\xitem[] $\gamma_\alpha<\delta_\alpha<\gamma_\beta<\delta_\beta$ for all
	$\alpha<\beta<\omega_2$\,; 
	\xitem[] $p_\alpha\leq_{\poQ^\dagger_{\bar{X}}}p$\,, 
	$p_\alpha=\pairof{c^\alpha,d^\alpha}$ for all $\alpha<\omega_2$\,; and 
	\xitem[] $p_\alpha\forces{\poQ^\dagger_{\bar{X}}}{\gamma_\alpha\in\dot{S}_0\,,
	\delta_\alpha\in\dot{S}_1}$. 
\end{xitemize}
By $\Delta$-System Lemma, we find a stationary $U\subseteq E^{\omega_2}_{\omega_1}$
and $A_\alpha$, $B_\alpha\in[\omega_2]^{<\aleph_0}$ for $\alpha\in U$ \st\ 
\begin{xitemize}
	\xitem[5.4-6] $\supp_0(p_\alpha)\subseteq A_\alpha$,
	$\supp_1(p_\alpha)\subseteq B_\alpha$\,;
	\xitem[5.4-7] $A_\alpha$, $\alpha\in U$ form a $\Delta$-system with root
	$A$;\ \ 
	$B_\alpha$, $\alpha\in U$ form a $\Delta$-system with root $B$; 
	\xitem[] $\gamma_\alpha$, $\gamma_\alpha+1$, $\delta_\alpha$,
	$\delta_\alpha+1\in (A_\alpha\cap B_\alpha)\setminus(A\cup B)$\,.
\end{xitemize}
By thinning out $U$ further, if necessary, 
we may also assume that there are some $k^*$, $n^*\in\omega$ \st\ 
\begin{xitemize}
	\xitem[5.4-8] $\dom(c^\alpha)=\supp_0(p_\alpha)\times k^*$ and 
	$\dom(s^{d^\alpha(\xi)})=k^*$ for all $\xi\in\supp_1(p_\alpha)$\,;
	\xitem[5.4-8-0] $c^\alpha(\gamma_\alpha,\cdot)=r^0_{n^*}$,
	$s^{d^\alpha(\gamma_\alpha+1)}=s^0_{n^*}$\,;
	\xitem[5.4-8-1]  $c^\alpha(\delta_\alpha,\cdot)=r^1_{n^*}$,
	$s^{d^\alpha(\delta_\alpha+1)}=s^1_{n^*}$\,.
\end{xitemize}
\Wolog, we may also assume that, for some fixed
$c^*$, $d^*$,
\begin{xitemize}
	\xitem[5.4-9] $c^\alpha\restr A\times k^*=c^*$ and
	$\seqof{s^{d^\alpha(\eta)}}{\eta\in B}=d^*$ for all $\alpha\in U$. 
\end{xitemize}
Note that $p_\alpha$, $\alpha\in U$ are compatible by \xitemof{5.4-6}, 
\xitemof{5.4-7} and \xitemof{5.4-9}. 

Now, since $\dot{S}_0$, $\dot{S}_1$, $p$ were arbitrary, 
\Claimof{main-claim},\,\assertof{1} is proved by the following subclaim:
\begin{Subclaim}
	For any $\alpha$, $\beta\in U$ with $\alpha<\beta$, there is 
$q\leq_{\poQ^\dagger_{\bar{X}}}p$ \st\ 
	\begin{xitemize}
	\item[] $q\forces{\poQ^\dagger_{\bar{X}}}{
	\gamma_\alpha\in\dot{S}_0,\ \delta_\beta\in\dot{S}_1,\ 
	\dot{a}_{\gamma_\alpha,n^*}\cap\dot{a}_{\delta_\beta,n^*}=\emptyset}$. 
	\end{xitemize}
\end{Subclaim}
\prfofClaim
Let $q=\pairof{c^q,d^q}$ be the common extension of $p_\alpha$ and $p_\beta$ 
\st\ 
\begin{xitemize}
	\xitem[5.4-10] $\gamma_\alpha\in F^{d^q(\delta_\beta+1)}$ 
	\xitem[5.4-10-0] $\dom(s^{d^q(\xi)})=k^*$ for all $\xi\in\dom(d^q)$. 
\end{xitemize}

Let $G$ be a $(V,\poQ^\dagger_{\bar{X}})$-generic filter with $q\in G$. 
In $V[G]$, we have 
\begin{xitemize}
	\xitem[5.4-11] $c_{\gamma_\alpha}(m)\leq d_{\delta_\beta+1}(m)$ 
	for all $m\geq k^*$ 
\end{xitemize}
by \xitemof{5.4-8}, \xitemof{5.4-10} and \xitemof{5.4-10-0}. 

Now, toward a contradiction, assume that
$a_{\gamma_\alpha,n^*}\cap a _{\delta_\beta,n^*}\not=\emptyset$ and let 
$m\in a_{\gamma_\alpha,n^*}\cap a_{\delta_\beta,n^*}$. 
By the definition of $a_\alpha$'s it follows that, for some $i$, $j\in 2$, 
we have 
\begin{xitemize}
\item[] $r^i_m\subseteq c_{\gamma_\alpha}$,
	$s^i_m\subseteq d_{\gamma_\alpha+1}$\,;
\item[] $r^j_m\subseteq c_{\delta_\beta}$,
	$s^j_m\subseteq d_{\delta_\beta+1}$. 
\end{xitemize}
On the other hand, since $q\in G$, we have $p_\alpha$, $p_\beta\in G$. It 
follows that 
\begin{xitemize}
\item[] $r^0_{n^*}\subseteq c_{\gamma_\alpha}$,
	$s^0_{n^*}\subseteq d_{\gamma_\alpha+1}$\,;
\item[] $r^1_{n^*}\subseteq c_{\delta_\beta}$, 
	$s^1_{n^*}\subseteq d_{\delta_\beta+1}$
\end{xitemize}
by \xitemof{5.4-8-0} and \xitemof{5.4-8-1}. 
By the definition \xitemof{5.4-5-0}
of $a_{\gamma_\alpha,n}$'s, we have
$\cardof{r^0_m}\geq n^*$. 
Thus we have, either
\begin{xitemize}
\item[] $r^0_{n^*}\subseteq r^0_m\subseteq c_{\gamma_\alpha}$\,, 
	$s^0_{n^*}\subseteq s^0_m\subseteq d_{\gamma_\alpha+1}$\,;
\item[] $r^1_{n^*}\subseteq r^1_m\subseteq c_{\gamma_\beta}$\,, 
	$s^1_{n^*}\subseteq s^1_m\subseteq d_{\gamma_\beta+1}$\,;
\end{xitemize}
or 
\begin{xitemize}
\item[] $r^0_{n^*}\subseteq r^1_m\subseteq c_{\gamma_\alpha}$\,, 
	$s^0_{n^*}\subseteq s^1_m\subseteq d_{\gamma_\alpha+1}$\,;
\item[] $r^1_{n^*}\subseteq r^0_m\subseteq c_{\gamma_\beta}$\,, 
	$s^1_{n^*}\subseteq s^0_m\subseteq d_{\gamma_\beta+1}$\,.
\end{xitemize}
%%% 
In the first case, we must have $c_{\gamma_\alpha}(m+1)=2$ and
$d_{\delta_\beta+1}(m+1)=1$ by \xitemof{5.4-4} and \xitemof{5.4-5}. This is 
a contradiction to \xitemof{5.4-11}. Similarly, in the second case, we have 
$c_{\gamma_\alpha}(m)=2$ and $d_{\delta_\beta+1}(m)=1$. This is again a 
contradiction to \xitemof{5.4-11}. \\
\qedofSubclaim\qedskip

\assertof{2} of \Claimof{main-claim} follows from the next subclaim: 
\begin{Subclaim}
	For any $t\in\fnsp{2}{\omega}$ and $\alpha$, $\beta\in U$ with
	$\alpha<\beta$, there is $q\leq_{\poQ^\dagger_{\bar{X}}}p$ \st\ 
	\begin{xitemize}
	\item[] $q\forces{\poQ^\dagger_{\bar{X}}}{
	\gamma_\alpha\in\dot{S}_0,\ \delta_\beta\in\dot{S}_1,\ 
	\dot{a}_{\gamma(\alpha),t(0)}\cap\dot{a}_{\delta(\beta),t(1)}
	\not=\emptyset}$.
	\end{xitemize}
\end{Subclaim}
\prfofClaim
For each $\xi\in\ssetof{\alpha,\beta}$, let 
$\tilde{p}_\xi\leq_{\poQ^\dagger_{\bar{X}}}p_\xi$ with 
$\tilde{p}_\xi=\pairof{\tilde{c}^\xi,\tilde{d}^\xi}$ and $m\in\omega$ 
be \st
\begin{xitemize}
\xitem[5.4-11-0] $\tilde{c}^\xi(\gamma_\xi,\cdot)=r^0_m$\,, 
	$s^{\tilde{d}^\xi(\gamma_\xi+1)}=s^0_m$\,;\ \ 
	$\tilde{c}^\xi(\delta_\xi,\cdot)=r^1_m$\,, 
	$s^{\tilde{d}^\xi(\delta_\xi+1)}=s^1_m$\,; 
\xitem[5.4-11-a] $\cardof{r^0_m}\geq t(0), t(1)$\,;
\xitem[5.4-11-1] $\supp_0(\tilde{p}_\xi)=\supp_0(p_\xi)$\,; 
	$\supp_1(\tilde{p}_\xi)=\supp_1(p_\xi)$\,;  
\xitem[5.4-11-2] $\tilde{c}^\xi\restr A\times\omega=c^\xi\restr A\times\omega$ and
	$\seqof{s^{\tilde{d}^\xi(\eta)}}{\eta\in B}
	=\seqof{s^{d^\xi(\eta)}}{\eta\in B}$.
\end{xitemize}
Let $q^0=\pairof{c^{q^0},d^{q^0}}$ be the maximal (\wrt\
$\leq_{\poQ^\dagger_{\bar{X}}}$) common extension of $\tilde{p}_\alpha$ and
$\tilde{p}_\beta$ which exists because of \xitemof{5.4-11-1} and 
\xitemof{5.4-11-2}. Extend $q^0$ further to $q=\pairof{c^q,d^q}$ \st\  
\begin{xitemize}
	\xitem[] $\cardof{c^q(\gamma_\alpha,\cdot)}
	=\cardof{c^q(\delta_\beta,\cdot)}
	=\cardof{s^{d^q(\gamma_\alpha+1)}}
	=\cardof{s^{d^q(\delta_\beta+1)}}=m+2$\,;
	\xitem[5.4-12] $c^q(\gamma_\alpha,m)=0$, $s^{d^q(\gamma_\alpha+1)}(m)=1$, 
	$c^q(\gamma_\alpha,m+1)=2$, $s^{d^q(\gamma_\alpha+1)}(m+1)=3$\,; 
	\xitem[5.4-13] $c^q(\delta_\beta,m)=2$, $s^{d^q(\delta_\beta+1)}(m)=3$, 
	$c^q(\delta_\beta,m+1)=0$, $s^{d^q(\delta_\beta+1)}(m+1)=1$\,.  
\end{xitemize}
This is possible because $\gamma_\alpha\not\in F^{d^{q^0}(\delta_\beta+1)}$ 
and $\delta_\beta\not\in F^{d^{q^0}(\gamma_\alpha+1)}$ by the maximality of
$q^0$ and \xitemof{5.4-11-1}. 

By \xitemof{5.4-11-0}, \xitemof{5.4-12}, \xitemof{5.4-13}, by the 
definition \xitemof{5.4-5-0}
of $a_{\alpha,n}$'s, and since
$\cardof{r^0_m}\geq t(0)$, $t(1)$,  we have  
$q\forces{\poQ^\dagger_{\bar{X}}}{
	m\in \dot{a}_{\gamma_\alpha, t(0)}\cap\dot{a}_{\delta_\beta, t(1)}}$.
Since $q\leq_{\poQ^\dagger_{\bar{X}}}p_\alpha$, $p_\beta$, we also have 
$q\forces{\poQ^\dagger_{\bar{X}}}{\gamma_\alpha\in\dot{S}_0,\ 
	\delta_\beta\in\dot{S}_1}$. 

Thus, $q$ as above is as desired. 
\qedofSubclaim\\
\qedofClaim\\
\qedofThm
\qedskip
\\
Note that in the proof of $\neg\Cs(\aleph_2)$ in \Thmabove, we used only 
\xitemof{Omega-0} from the assumption (A). Note also that this proof 
actually shows that in the generic extension the negation of
${\rm C}(\aleph_2)$ from \cite{ju-so-sz} 
holds 
which is a weakening of $\Cs(\aleph_2)$ obtained by 
replacing the condition ``stationary'' in the formulation of $\Cs(\aleph_2)$ 
by ``cofinal''.  

\section{Forcing CH $+$ (A)}
\label{forcing(A)}
In this section, we define under \CH\ a $\sigma$-closed $\aleph_2$-c.c.\ 
\po\ $\poP_0$ which forces the combinatorial assertion (A) of the previous 
section. 

The \po\ $\poP_0$ is defined as follows:
\newcommand{\thecond}[1]{\pairof{X^{#1},Y^{#1},\,%
		\tau^{#1},\,%
	\seqof{\phi^{#1}_\xi,\psi^{#1}_\xi}{\xi\in D^{#1}}%
		}}
\begin{xitemize}
	\item[] $p\in\poP_0$\ \ $\Leftrightarrow$\ \ $p=\thecond{p}$ 
\end{xitemize}
where
\begin{xitemize}
	\xitem[(A)-1] $X^p$, $Y^p\in[\omega_2]^{\aleph_0}$\,;
	\xitem[(A)-2] $D^{p}\in[\omega_2]^{\aleph_0}$\,;
	\xitem[(A)-3] for all $\xi\in D^p$, $\mapping{\phi^p_\xi}{X^p}{X^p}$ 
	and $\mapping{\psi^p_\xi}{Y^p}{Y^p}$ are involutions (that is, bijections 
	$\phi$ 
	\st\ $\phi^{-1}=\phi$);
	\xitem[(A)-4] for all $\xi\in D^p$, $\alpha\in X^p$ and $\beta\in Y^p$,
	\begin{xxitemize}
	\xxitem[(A)-4][a] $\phi^p_\xi(\alpha)<\alpha+\xi+\omega_1$ and
	\xxitem[(A)-4][b] $\psi^p_\xi(\beta)<\beta+\xi+\omega_1$\,;
	\end{xxitemize}
\end{xitemize}
%% We shall refer the situation like in \xxitemof{(A)-4}{a} as {\em the distance 
%% of $\phi^p_\xi(\alpha)$ from $\alpha$ 
%% being strictly less than $\xi+\omega_1$}. 
Note that we have also 
$\alpha<\phi^p_\xi(\alpha)+\xi+\omega_1$ and
$\alpha<\psi^p_\xi(\beta)+\xi+\omega_1$ for all $\xi\in D^p$, 
$\alpha\in X^p$ and $\beta\in Y^p$ since $\phi^p_\xi$ and $\psi^p_\xi$ are 
involutions by \xitemof{(A)-3}. 
\begin{xitemize}
	\xitem[(A)-5] $\mapping{\tau^p}{X^p\times Y^p}{2}$\,;
	\xitem[(A)-6] for all $\xi\in D^p$, $\alpha\in X^p$ and $\beta\in Y^p$, 
	we have
	$\tau^p(\alpha,\beta)=\tau^p(\phi^p_\xi(\alpha),\psi^p_\xi(\beta))$\,;
	\xitem[(A)-7] $\tau^p(\alpha,\beta)=1$ for all
	$\pairof{\alpha,\beta}\in X^p\times Y^p$ with $\beta<\alpha$. 
\end{xitemize}

The ordering on $\poP_0$ is defined by the following: 
For $p$, $q\in\poP_0$ 
with 
\begin{xitemize}
\item[] $p=\thecond{p}$ and\\[\jot]
	$q=\thecond{q}$,
\end{xitemize}
\begin{xitemize}
	\xitem[(A)-8] $p\leq_{\poP_0}q$\ \ $\Leftrightarrow$\ \ 
	\begin{xxitemize}
	\xxitem[(A)-8][a] $X^p\supseteq X^q$, $Y^p\supseteq Y^q$\,;
	\xxitem[(A)-8][b] $D^p\supseteq D^q$\,;
	\xxitem[(A)-8][c] $\phi^p_\xi\supseteq \phi^q_\xi$ and 
	$\psi^p_\xi\supseteq \psi^q_\xi$ for all $\xi\in D^q$\,;
	\xxitem[(A)-8][d] $\tau^p\supseteq \tau^q$  and 
	\xxitem[(A)-8][e] $\tau^p\restr(X^p\setminus X^q)\times Y^q\equiv 1$. 
	\end{xxitemize}
\end{xitemize}

For $p\in\poP_0$ with $p=\thecond{p}$, we intend to approximate the 
characteristic function of the set $A$ in the assertion (A) 
by $\tau^p$. More precisely, in a generic 
extension $V[G]$ for a $(V,\poP_0)$-generic $G$, letting 
\begin{xitemize}
\xitem[p0-0] $\tau=\bigcup_{p\in G}\tau^p$\,;\quad
	$\phi_\xi=\bigcup_{p\in G}\phi_\xi^p$\ \ and\ \  
	$\psi_\xi=\bigcup_{p\in G}\psi_\xi^p$\ \  for $\xi\in\omega_2$\,;\vspace{\jot}
\xitem[p0-1]
	$A=\tau^{-1}\imageof\ssetof{1}$\ \ 
	and\ \ 
	$\calF=\setof{\pairof{\phi_\xi,\psi_\xi}}{\xi\in\omega_2}$, 
\end{xitemize}
we are aiming to force $\pairof{(\omega_2)^2,A,\calF}$ to satisfy 
\xitemof{Omega-0} $\sim$ \xitemof{Omega-4} in (A). 

Of the conditions in the definition of $\poP_0$, \xitemof{(A)-5} and 
\xxitemof{(A)-8}{d} 
force $\tau$ to be a function. Furthermore, 
$\mapping{\tau}{\omega_2\times\omega_2}{2}$ by density 
argument and the following \Lemmaof{density-argument},\,\assertof{a}. 

\xitemof{(A)-3} and 
\xxitemof{(A)-8}{c} make 
$\phi_\xi$ and 
$\psi_\xi$ mappings for all $\xi\in\omega_2$;   
they are forced to be involutions on $\omega_2$ by \xitemof{(A)-3} and the following 
\Lemmaof{density-argument},\,\assertof{a}. Thus 
$\pairof{(\omega_2)^2,A,\calF}$ is forced to satisfy \xitemof{Omega-2}. 

By 
\xitemof{(A)-7} (and by the following 
\Lemmaof{density-argument},\,\assertof{a}), $\pairof{(\omega_2)^2,A,\calF}$ 
is forced to satisfy the second inclusion of \xitemof{Omega-0}. 

By \xitemof{(A)-6}, 
$\pairof{(\omega_2)^2,A,\calF}$ is forced to satisfy \xitemof{Omega-3}.  

\xitemof{(A)-4} and \xxitemof{(A)-8}{e} are technical conditions
whose role will be clear later in the course of the proof. 

By the definition of $\poP_0$, it is 
clear that $\poP_0$ is $\sigma$-closed.  
Thus, we are done by showing that $\poP_0$ satisfies the $\aleph_2$-c.c.\ 
and it 
forces that 
$\pairof{(\omega_2)^2,A,\calF}$ as above satisfies the conditions 
\xitemof{Omega-1} and \xitemof{Omega-4}. 

The next Lemma follows readily from the definition of $\poP_0$. 

\begin{Lemma}\label{density-argument} \assert{a} For any $\alpha$,
	$\beta<\omega_2$, the set 
	\begin{xitemize}
	\item[] $
		\begin{array}{r@{}l}
			D_{\alpha,\beta}=\setof{\,p\in\poP_0}{{}&p=\thecond{p},\\[\jot]
				&\alpha\in X^p\mbox{ and }\beta\in Y^p\quad}
		\end{array}
		$
	\end{xitemize}
	is dense in $\poP_0$. 
	\smallskip

	\assert{b} For any $C\in[\omega_2]^{\aleph_0}$ and any $\beta<\omega_2$, 
	\begin{xitemize}
	\item[]
		$
		\begin{array}{r@{}l}
			E_{C,\beta}=\setof{\, p\in\poP_0}{{}&p=\thecond{p},\\[\jot]
				&C\subseteq X^p 
				\mbox{ and for some }\delta\in Y^p\mbox{ with }\delta\geq\beta\\[\jot] 
				&\tau^p(\gamma,\delta)=0\mbox{ for all }\gamma\in C\quad}			
		\end{array}
		$
	\end{xitemize}
	is dense in $\poP_0$.\qed
\end{Lemma}

In the rest of the section, we are going to work mainly in the ground model (where 
\CH\ holds). 
Let $\dot{\tau}$, $\dot{\phi}_\xi$, $\dot{\psi}_\xi$ for $\xi\in\omega_2$,  
$\dot{A}$ and $\dot\calF$ be $\poP_0$-names 
of $\tau$, $\phi_\xi$, $\psi_\xi$ for $\xi\in\omega_2$, $A$ and $\calF$ as 
above, respectively. 
\begin{Lemma}
	$\forces{\poP_0}{\pairof{(\omega_2)^2,\dot{A},\dot\calF}
		\models\mbox{\xitemof{Omega-1}}}$.
\end{Lemma}
\prf By density argument with \Lemmaof{density-argument},\,\assertof{b}.\qedofLemma
\qedskip 

For $\xi<\omega_2$, $X$, $X'$, $Y$, $Y'\in[\omega_2]^{\aleph_0}$ with
$X'\subseteq X$ and $Y'\subseteq Y$, $\mapping{\tau}{X\times Y}{2}$ and 
involutions $\mapping{\phi'}{X'}{X'}$, $\mapping{\psi}{Y'}{Y'}$, 
let us call the quintuple
$\pairof{X,Y,\tau,\phi',\psi'}$ a {\em$\xi$-extendable semi-condition} if 
\begin{xitemize}
	\xitem[7-6] $\phi'(\alpha)<\alpha+\xi+\omega_1$\ \ and\ \ 
			$\psi'(\beta)<\beta+\xi+\omega_1$\ \ for all $\alpha\in X'$ and $\beta\in Y'$\,;
	\xitem[7-7] $\tau(\alpha,\beta)=\tau(\phi(\alpha),\psi(\beta))$\ \  for 
	all $\alpha\in X'$ and $\beta\in Y'$\,; 
	\xitem[7-8] $\tau(\alpha,\beta)=1$\ \ for all $\alpha\in X$ and $\beta\in Y$ 
	with $\beta<\alpha$\,; 
	\xitem[7-9] $\tau\restr((X\setminus X')\times Y')\equiv 1$. 
\end{xitemize}
\mbox{}\\
\mbox{}\hfill
%%\input{coloring-x-fig-7-9}
%WinTpicVersion3.08
\unitlength 1cm
\begin{picture}(  3.5875,  3.6875)(  0.4125, -5.1625)
% STR 2 0 3 0
% 3 1000 1720 1000 1820 5 0
% $X'$
\put(2.5000,-4.5500){\makebox(0,0){$X'$}}%
% STR 2 0 3 0
% 3 1200 2050 1200 2150 5 0
% $X$
\put(3.0000,-5.3750){\makebox(0,0){$X$}}%
% LINE 2 0 3 0
% 6 800 1640 800 1670 800 1670 1200 1670 1200 1670 1200 1640
% 
\special{pn 8}%
\special{pa 788 1615}%
\special{pa 788 1644}%
\special{fp}%
\special{pa 788 1644}%
\special{pa 1182 1644}%
\special{fp}%
\special{pa 1182 1644}%
\special{pa 1182 1615}%
\special{fp}%
% LINE 2 0 3 0
% 6 760 1600 730 1600 730 1600 730 1000 730 1000 760 1000
% 
\special{pn 8}%
\special{pa 749 1575}%
\special{pa 719 1575}%
\special{fp}%
\special{pa 719 1575}%
\special{pa 719 985}%
\special{fp}%
\special{pa 719 985}%
\special{pa 749 985}%
\special{fp}%
% LINE 2 0 3 0
% 10 790 1940 790 1940 790 1970 790 1970 790 1970 1590 1970 1590 1970 1590 1940 790 1940 790 1970
% 
\special{pn 8}%
\special{pa 778 1910}%
\special{pa 778 1910}%
\special{fp}%
\special{pa 778 1939}%
\special{pa 778 1939}%
\special{fp}%
\special{pa 778 1939}%
\special{pa 1565 1939}%
\special{fp}%
\special{pa 1565 1939}%
\special{pa 1565 1910}%
\special{fp}%
\special{pa 778 1910}%
\special{pa 778 1939}%
\special{fp}%
% STR 2 0 3 0
% 3 580 1300 580 1400 5 0
% $Y'$
\put(1.4500,-3.5000){\makebox(0,0){$Y'$}}%
% LINE 2 0 3 0
% 8 470 1600 470 1600 470 1600 440 1600 440 1600 440 600 440 600 470 600
% 
\special{pn 8}%
\special{pa 463 1575}%
\special{pa 463 1575}%
\special{fp}%
\special{pa 463 1575}%
\special{pa 434 1575}%
\special{fp}%
\special{pa 434 1575}%
\special{pa 434 591}%
\special{fp}%
\special{pa 434 591}%
\special{pa 463 591}%
\special{fp}%
% STR 2 0 3 0
% 3 300 1110 300 1210 5 0
% $Y$
\put(0.7500,-3.0250){\makebox(0,0){$Y$}}%
% LINE 1 0 3 0
% 10 800 600 800 600 800 600 800 1600 800 1600 1600 1600 1600 1600 1600 590 1600 590 800 590
% 
\special{pn 13}%
\special{pa 788 591}%
\special{pa 788 591}%
\special{fp}%
\special{pa 788 591}%
\special{pa 788 1575}%
\special{fp}%
\special{pa 788 1575}%
\special{pa 1575 1575}%
\special{fp}%
\special{pa 1575 1575}%
\special{pa 1575 581}%
\special{fp}%
\special{pa 1575 581}%
\special{pa 788 581}%
\special{fp}%
% LINE 1 1 3 0
% 6 800 1190 800 1190 800 990 1600 990 1200 590 1200 1590
% 
\special{pn 13}%
\special{pa 788 1172}%
\special{pa 788 1172}%
\special{da 0.050}%
\special{pa 788 975}%
\special{pa 1575 975}%
\special{da 0.050}%
\special{pa 1182 581}%
\special{pa 1182 1565}%
\special{da 0.050}%
% LINE 1 0 3 0
% 2 800 1590 1600 1090
% 
\special{pn 13}%
\special{pa 788 1565}%
\special{pa 1575 1073}%
\special{fp}%
% LINE 2 0 3 0
% 10 1190 1390 990 1590 1140 1380 930 1590 980 1480 870 1590 1190 1450 1050 1590 1190 1510 1110 1590
% 
\special{pn 8}%
\special{pa 1172 1369}%
\special{pa 975 1565}%
\special{fp}%
\special{pa 1123 1359}%
\special{pa 916 1565}%
\special{fp}%
\special{pa 965 1457}%
\special{pa 857 1565}%
\special{fp}%
\special{pa 1172 1428}%
\special{pa 1034 1565}%
\special{fp}%
\special{pa 1172 1487}%
\special{pa 1093 1565}%
\special{fp}%
% LINE 2 0 3 0
% 20 1590 1410 1410 1590 1590 1350 1350 1590 1590 1290 1290 1590 1590 1230 1230 1590 1590 1170 1200 1560 1590 1110 1200 1500 1460 1180 1200 1440 1300 1280 1200 1380 1590 1470 1470 1590 1590 1530 1530 1590
% 
\special{pn 8}%
\special{pa 1565 1388}%
\special{pa 1388 1565}%
\special{fp}%
\special{pa 1565 1329}%
\special{pa 1329 1565}%
\special{fp}%
\special{pa 1565 1270}%
\special{pa 1270 1565}%
\special{fp}%
\special{pa 1565 1211}%
\special{pa 1211 1565}%
\special{fp}%
\special{pa 1565 1152}%
\special{pa 1182 1536}%
\special{fp}%
\special{pa 1565 1093}%
\special{pa 1182 1477}%
\special{fp}%
\special{pa 1438 1162}%
\special{pa 1182 1418}%
\special{fp}%
\special{pa 1280 1260}%
\special{pa 1182 1359}%
\special{fp}%
\special{pa 1565 1447}%
\special{pa 1447 1565}%
\special{fp}%
\special{pa 1565 1506}%
\special{pa 1506 1565}%
\special{fp}%
% LINE 2 0 3 0
% 16 1580 1000 1320 1260 1530 990 1200 1320 1470 990 1200 1260 1410 990 1200 1200 1350 990 1200 1140 1290 990 1200 1080 1230 990 1200 1020 1590 1050 1480 1160
% 
\special{pn 8}%
\special{pa 1556 985}%
\special{pa 1300 1241}%
\special{fp}%
\special{pa 1506 975}%
\special{pa 1182 1300}%
\special{fp}%
\special{pa 1447 975}%
\special{pa 1182 1241}%
\special{fp}%
\special{pa 1388 975}%
\special{pa 1182 1182}%
\special{fp}%
\special{pa 1329 975}%
\special{pa 1182 1123}%
\special{fp}%
\special{pa 1270 975}%
\special{pa 1182 1063}%
\special{fp}%
\special{pa 1211 975}%
\special{pa 1182 1004}%
\special{fp}%
\special{pa 1565 1034}%
\special{pa 1457 1142}%
\special{fp}%
\end{picture}%
\mbox{}\hfill\mbox{}\bigskip\bigskip\\
\mbox{}\hfill {\bf fig.\ 6}\hfill\mbox{}\bigskip\bigskip\\

Note that the sets $X'$ and
$Y'$, though not mentioned explicitly in the definition 
of $\xi$-extendable semi-condition, 
can be recovered from $\phi'$ and $\psi'$. 

Note also that \xitemof{7-9} holds vacuously if $X=X'$. Hence, if
$p\in\poP_0$ with $p=\thecond{p}$,  
the quintuple $\pairof{X^p,Y^p,\tau^p,\phi_\xi,\psi_\xi}$ is a $\xi$-extendable 
semi-condition %% with 
%% $(X^p)'=X^p$ and $(Y^p)'=Y^p$ 
for all $\xi\in D^p$. 

The following two lemmas explain the choice of the naming of $\xi$-extendable 
semi-conditions.  
\begin{Lemma}\label{xi-extendable}
	For any $\xi<\omega_2$, $X$, $X'$, $Y$, $Y'\in[\omega_2]^{\aleph_0}$ with
	$X'\subseteq X$ and $Y'\subseteq Y$, $\mapping{\tau}{X\times Y}{2}$ as well 
	as involutions $\mapping{\phi'}{X'}{X'}$, $\mapping{\psi'}{Y'}{Y'}$, 
	if $\pairof{X,Y,\tau,\phi',\psi'}$ is a $\xi$-extendable 
	semi-condition then there are 
	$\tilde{X}\supseteq X$, $\tilde{Y}\supseteq Y$,
	$\mapping{\tilde{\tau}}{\tilde{X}\times\tilde{Y}}{2}$ with
	$\tilde{\tau}\supseteq\tau$ and involutions
	$\mapping{\tilde{\phi}}{\tilde{X}}{\tilde{X}}$, 
	$\mapping{\tilde{\psi}}{\tilde{Y}}{\tilde{Y}}$ extending $\phi'$ and
	$\psi'$ respectively \st\
	$\pairof{\tilde{X},\tilde{Y},\tilde{\tau},\tilde{\phi},\tilde{\psi}}$ 
	is a $\xi$-extendable semi-condition and 
	\begin{xitemize}
		\xitem[7-9a] 
		$\tilde{\tau}\restr((\tilde{X}\setminus X)\times Y)\equiv 1$. 
	\end{xitemize}
\mbox{}\hfill
%%\input{coloring-x-fig-7-9a}
%WinTpicVersion3.08
\unitlength 1cm
\begin{picture}(  6.2500,  4.9125)( -0.7375, -5.9125)
% LINE 1 0 3 0
% 12 1005 400 1005 1600 1005 1600 2205 1600 2205 1600 2205 400 2205 400 1005 400 1005 800 2205 800 1805 400 1805 1600
% 
\special{pn 13}%
\special{pa 990 394}%
\special{pa 990 1575}%
\special{fp}%
\special{pa 990 1575}%
\special{pa 2171 1575}%
\special{fp}%
\special{pa 2171 1575}%
\special{pa 2171 394}%
\special{fp}%
\special{pa 2171 394}%
\special{pa 990 394}%
\special{fp}%
\special{pa 990 788}%
\special{pa 2171 788}%
\special{fp}%
\special{pa 1777 394}%
\special{pa 1777 1575}%
\special{fp}%
% STR 2 0 3 0
% 3 1205 1720 1205 1820 5 0
% $X'$
\put(3.0125,-4.5500){\makebox(0,0){$X'$}}%
% STR 2 0 3 0
% 3 1410 2010 1410 2110 5 0
% $X$
\put(3.5250,-5.2750){\makebox(0,0){$X$}}%
% LINE 2 0 3 0
% 6 1005 1640 1005 1670 1005 1670 1405 1670 1405 1670 1405 1640
% 
\special{pn 8}%
\special{pa 990 1615}%
\special{pa 990 1644}%
\special{fp}%
\special{pa 990 1644}%
\special{pa 1383 1644}%
\special{fp}%
\special{pa 1383 1644}%
\special{pa 1383 1615}%
\special{fp}%
% LINE 2 0 3 0
% 6 965 1600 935 1600 935 1600 935 1200 935 1200 965 1200
% 
\special{pn 8}%
\special{pa 950 1575}%
\special{pa 921 1575}%
\special{fp}%
\special{pa 921 1575}%
\special{pa 921 1182}%
\special{fp}%
\special{pa 921 1182}%
\special{pa 950 1182}%
\special{fp}%
% LINE 2 0 3 0
% 10 995 1940 995 1940 995 1970 995 1970 995 1970 1795 1970 1795 1970 1795 1940 995 1940 995 1970
% 
\special{pn 8}%
\special{pa 980 1910}%
\special{pa 980 1910}%
\special{fp}%
\special{pa 980 1939}%
\special{pa 980 1939}%
\special{fp}%
\special{pa 980 1939}%
\special{pa 1767 1939}%
\special{fp}%
\special{pa 1767 1939}%
\special{pa 1767 1910}%
\special{fp}%
\special{pa 980 1910}%
\special{pa 980 1939}%
\special{fp}%
% STR 2 0 3 0
% 3 785 1300 785 1400 5 0
% $Y'$
\put(1.9625,-3.5000){\makebox(0,0){$Y'$}}%
% LINE 2 0 3 0
% 8 675 1600 675 1600 675 1600 645 1600 645 1600 645 800 645 800 675 800
% 
\special{pn 8}%
\special{pa 665 1575}%
\special{pa 665 1575}%
\special{fp}%
\special{pa 665 1575}%
\special{pa 635 1575}%
\special{fp}%
\special{pa 635 1575}%
\special{pa 635 788}%
\special{fp}%
\special{pa 635 788}%
\special{pa 665 788}%
\special{fp}%
% STR 2 0 3 0
% 3 505 1110 505 1210 5 0
% $Y$
\put(1.2625,-3.0250){\makebox(0,0){$Y$}}%
% LINE 1 1 3 0
% 4 1405 400 1405 1600 1005 1200 2205 1200
% 
\special{pn 13}%
\special{pa 1383 394}%
\special{pa 1383 1575}%
\special{da 0.050}%
\special{pa 990 1182}%
\special{pa 2171 1182}%
\special{da 0.050}%
% LINE 2 0 3 0
% 24 2195 890 1895 1190 2195 830 1835 1190 2165 800 1805 1160 2105 800 1805 1100 2045 800 1805 1040 1985 800 1805 980 1925 800 1805 920 1865 800 1805 860 2195 950 1955 1190 2195 1010 2015 1190 2195 1070 2075 1190 2195 1130 2135 1190
% 
\special{pn 8}%
\special{pa 2161 876}%
\special{pa 1866 1172}%
\special{fp}%
\special{pa 2161 817}%
\special{pa 1807 1172}%
\special{fp}%
\special{pa 2131 788}%
\special{pa 1777 1142}%
\special{fp}%
\special{pa 2072 788}%
\special{pa 1777 1083}%
\special{fp}%
\special{pa 2013 788}%
\special{pa 1777 1024}%
\special{fp}%
\special{pa 1954 788}%
\special{pa 1777 965}%
\special{fp}%
\special{pa 1895 788}%
\special{pa 1777 906}%
\special{fp}%
\special{pa 1836 788}%
\special{pa 1777 847}%
\special{fp}%
\special{pa 2161 936}%
\special{pa 1925 1172}%
\special{fp}%
\special{pa 2161 995}%
\special{pa 1984 1172}%
\special{fp}%
\special{pa 2161 1054}%
\special{pa 2043 1172}%
\special{fp}%
\special{pa 2161 1113}%
\special{pa 2102 1172}%
\special{fp}%
% LINE 2 0 3 0
% 6 370 400 340 400 340 400 340 1600 340 1600 370 1600
% 
\special{pn 8}%
\special{pa 365 394}%
\special{pa 335 394}%
\special{fp}%
\special{pa 335 394}%
\special{pa 335 1575}%
\special{fp}%
\special{pa 335 1575}%
\special{pa 365 1575}%
\special{fp}%
% STR 2 0 3 0
% 3 200 880 200 980 5 0
% $\tilde{Y}$
\put(0.5000,-2.4500){\makebox(0,0){$\tilde{Y}$}}%
% LINE 2 0 3 0
% 6 1000 2220 1000 2250 1000 2250 2200 2250 2200 2250 2200 2220
% 
\special{pn 8}%
\special{pa 985 2186}%
\special{pa 985 2215}%
\special{fp}%
\special{pa 985 2215}%
\special{pa 2166 2215}%
\special{fp}%
\special{pa 2166 2215}%
\special{pa 2166 2186}%
\special{fp}%
% STR 2 0 3 0
% 3 1560 2350 1560 2450 5 0
% $\tilde{X}$
\put(3.9000,-6.1250){\makebox(0,0){$\tilde{X}$}}%
% LINE 2 0 3 0
% 26 1800 1260 1470 1590 1790 1210 1420 1580 1740 1200 1410 1530 1680 1200 1410 1470 1620 1200 1410 1410 1560 1200 1410 1350 1500 1200 1410 1290 1440 1200 1410 1230 1800 1320 1530 1590 1800 1380 1590 1590 1800 1440 1650 1590 1800 1500 1710 1590 1800 1560 1770 1590
% 
\special{pn 8}%
\special{pa 1772 1241}%
\special{pa 1447 1565}%
\special{fp}%
\special{pa 1762 1191}%
\special{pa 1398 1556}%
\special{fp}%
\special{pa 1713 1182}%
\special{pa 1388 1506}%
\special{fp}%
\special{pa 1654 1182}%
\special{pa 1388 1447}%
\special{fp}%
\special{pa 1595 1182}%
\special{pa 1388 1388}%
\special{fp}%
\special{pa 1536 1182}%
\special{pa 1388 1329}%
\special{fp}%
\special{pa 1477 1182}%
\special{pa 1388 1270}%
\special{fp}%
\special{pa 1418 1182}%
\special{pa 1388 1211}%
\special{fp}%
\special{pa 1772 1300}%
\special{pa 1506 1565}%
\special{fp}%
\special{pa 1772 1359}%
\special{pa 1565 1565}%
\special{fp}%
\special{pa 1772 1418}%
\special{pa 1625 1565}%
\special{fp}%
\special{pa 1772 1477}%
\special{pa 1684 1565}%
\special{fp}%
\special{pa 1772 1536}%
\special{pa 1743 1565}%
\special{fp}%
% LINE 2 0 3 0
% 24 2160 1200 1810 1550 2200 1220 1830 1590 2200 1280 1890 1590 2200 1340 1950 1590 2200 1400 2010 1590 2200 1460 2070 1590 2200 1520 2130 1590 2100 1200 1810 1490 2040 1200 1810 1430 1980 1200 1810 1370 1920 1200 1810 1310 1860 1200 1810 1250
% 
\special{pn 8}%
\special{pa 2126 1182}%
\special{pa 1782 1526}%
\special{fp}%
\special{pa 2166 1201}%
\special{pa 1802 1565}%
\special{fp}%
\special{pa 2166 1260}%
\special{pa 1861 1565}%
\special{fp}%
\special{pa 2166 1319}%
\special{pa 1920 1565}%
\special{fp}%
\special{pa 2166 1378}%
\special{pa 1979 1565}%
\special{fp}%
\special{pa 2166 1438}%
\special{pa 2038 1565}%
\special{fp}%
\special{pa 2166 1497}%
\special{pa 2097 1565}%
\special{fp}%
\special{pa 2067 1182}%
\special{pa 1782 1467}%
\special{fp}%
\special{pa 2008 1182}%
\special{pa 1782 1408}%
\special{fp}%
\special{pa 1949 1182}%
\special{pa 1782 1349}%
\special{fp}%
\special{pa 1890 1182}%
\special{pa 1782 1290}%
\special{fp}%
\special{pa 1831 1182}%
\special{pa 1782 1231}%
\special{fp}%
\end{picture}%
\mbox{}\hfill\mbox{}\bigskip\bigskip\\
\mbox{}\hfill {\bf fig.\ 7}\hfill\mbox{}\\
\end{Lemma}
\prf Let $X_0\in[\omega_2\setminus X]^{\leq\aleph_0}$ be \st\ 
\begin{xitemize}
	\xitem[7star] $X_0$ is order-isomorphic to $X\setminus X'$ and the 
	order-isomorphism identifies points of distance less than $\omega_1$ 
	(that is, if $\alpha\in X\setminus X'$ and $\alpha_0\in X_0$ are 
	identified then we have $\alpha<\alpha_0+\omega_1$ and
	$\alpha_0<\alpha+\omega_1$). 
\end{xitemize}
Since $X$ (and hence also $X\setminus X'$) is countable, we can easily 
choose the elements of $X_0$ recursively in $otp(X\setminus X')$ steps in 
accordance with \xitemof{7star}. 
Put $\tilde{X}=X\dotcup X_0$ and let $\tilde{\phi}$ be the extension of 
$\phi$ which maps $X_0$ order-isomorphically to $X\setminus X'$ and 
vice versa. Fix $\theta<\omega_1$ \st\ 
\begin{xitemize}
	\xitem[boxplus] $\tilde{\phi}(\alpha)\leq\alpha+\xi+\theta$ for all
	$\alpha\in\tilde{X}$.
\end{xitemize}
There is such $\theta$ by \xitemof{7-6}, \xitemof{7star} and since
$\cardof{\tilde{X}}\leq\aleph_0$. 

Let $Y_0\in[\omega_2\setminus Y]^{\leq\aleph_0}$ be \st\ 
\begin{xitemize}
	\xitem[dstar-1] $Y_0$ is order-isomorphic to $Y\setminus Y'$ and the 
	order-isomorphism identifies points of distance less than $\xi+\omega_1$; 
	and 
	\xitem[dstar-2] if $\beta\in Y\setminus Y'$ and $\beta_0\in Y_0$ are 
	identified then $\beta_0>\beta+\xi+\theta$.
\end{xitemize}
It is easy to see that the elements of such $Y_0$ can be chosen recursively 
in $otp(Y\setminus Y')$ steps.

Now let $\tilde{Y}=Y\dotcup Y_0$ and let $\tilde{\psi}$ be the extension of 
$\psi'$ which maps $Y_0$ order-isomorphically to $Y\setminus Y'$ and vice 
versa. 

Finally define $\mapping{\tilde{\tau}}{\tilde{X}\times\tilde{Y}}{2}$ by
\begin{xitemize}
	\xitem[7-10] 
	$\tilde{\tau}(\alpha,\beta)=\left\{\,
	\begin{array}[c]{@{}ll}
		\tau(\alpha,\beta),\ 
		&\mbox{if }\pairof{\alpha,\beta}\in X\times Y\\[\jot]
		\tau(\tilde{\phi}(\alpha),\tilde{\psi}(\beta)),\ 
		&\mbox{if }
		\pairof{\alpha,\beta}\in (X'\times Y_0)\cup(X_0\times Y')\cup(X_0\times 
		Y_0)\\[\jot]
		1, &\mbox{otherwise}
	\end{array}\right.
	$
\end{xitemize}
for every $\alpha\in\tilde{X}$ and $\beta\in\tilde{Y}$. 

\mbox{}\\
\mbox{}\hfill
%%\input{coloring-x-fig-7-10}
%WinTpicVersion3.08
\unitlength 1cm
\begin{picture}(  4.5875,  4.1625)(  0.4125, -5.1625)
% LINE 1 0 3 0
% 12 800 400 800 1600 800 1600 2000 1600 2000 1600 2000 400 2000 400 800 400 800 800 2000 800 1600 400 1600 1600
% 
\special{pn 13}%
\special{pa 788 394}%
\special{pa 788 1575}%
\special{fp}%
\special{pa 788 1575}%
\special{pa 1969 1575}%
\special{fp}%
\special{pa 1969 1575}%
\special{pa 1969 394}%
\special{fp}%
\special{pa 1969 394}%
\special{pa 788 394}%
\special{fp}%
\special{pa 788 788}%
\special{pa 1969 788}%
\special{fp}%
\special{pa 1575 394}%
\special{pa 1575 1575}%
\special{fp}%
% STR 2 0 3 0
% 3 1000 1720 1000 1820 5 0
% $X'$
\put(2.5000,-4.5500){\makebox(0,0){$X'$}}%
% STR 2 0 3 0
% 3 1200 2050 1200 2150 5 0
% $X$
\put(3.0000,-5.3750){\makebox(0,0){$X$}}%
% STR 2 0 3 0
% 3 1730 1650 1730 1750 1 0
% $X_0=\tilde{X}\setminus X$
\put(4.3250,-4.3750){\makebox(0,0)[lt]{$X_0=\tilde{X}\setminus X$}}%
% LINE 2 0 3 0
% 6 800 1640 800 1670 800 1670 1200 1670 1200 1670 1200 1640
% 
\special{pn 8}%
\special{pa 788 1615}%
\special{pa 788 1644}%
\special{fp}%
\special{pa 788 1644}%
\special{pa 1182 1644}%
\special{fp}%
\special{pa 1182 1644}%
\special{pa 1182 1615}%
\special{fp}%
% LINE 2 0 3 0
% 6 1600 1640 1600 1670 1600 1670 2000 1670 2000 1670 2000 1640
% 
\special{pn 8}%
\special{pa 1575 1615}%
\special{pa 1575 1644}%
\special{fp}%
\special{pa 1575 1644}%
\special{pa 1969 1644}%
\special{fp}%
\special{pa 1969 1644}%
\special{pa 1969 1615}%
\special{fp}%
% LINE 2 0 3 0
% 6 760 1600 730 1600 730 1600 730 1200 730 1200 760 1200
% 
\special{pn 8}%
\special{pa 749 1575}%
\special{pa 719 1575}%
\special{fp}%
\special{pa 719 1575}%
\special{pa 719 1182}%
\special{fp}%
\special{pa 719 1182}%
\special{pa 749 1182}%
\special{fp}%
% LINE 2 0 3 0
% 6 760 800 730 800 730 800 730 400 730 400 760 400
% 
\special{pn 8}%
\special{pa 749 788}%
\special{pa 719 788}%
\special{fp}%
\special{pa 719 788}%
\special{pa 719 394}%
\special{fp}%
\special{pa 719 394}%
\special{pa 749 394}%
\special{fp}%
% LINE 2 0 3 0
% 10 790 1940 790 1940 790 1970 790 1970 790 1970 1590 1970 1590 1970 1590 1940 790 1940 790 1970
% 
\special{pn 8}%
\special{pa 778 1910}%
\special{pa 778 1910}%
\special{fp}%
\special{pa 778 1939}%
\special{pa 778 1939}%
\special{fp}%
\special{pa 778 1939}%
\special{pa 1565 1939}%
\special{fp}%
\special{pa 1565 1939}%
\special{pa 1565 1910}%
\special{fp}%
\special{pa 778 1910}%
\special{pa 778 1939}%
\special{fp}%
% STR 2 0 3 0
% 3 580 1300 580 1400 5 0
% $Y'$
\put(1.4500,-3.5000){\makebox(0,0){$Y'$}}%
% LINE 2 0 3 0
% 8 470 1600 470 1600 470 1600 440 1600 440 1600 440 800 440 800 470 800
% 
\special{pn 8}%
\special{pa 463 1575}%
\special{pa 463 1575}%
\special{fp}%
\special{pa 463 1575}%
\special{pa 434 1575}%
\special{fp}%
\special{pa 434 1575}%
\special{pa 434 788}%
\special{fp}%
\special{pa 434 788}%
\special{pa 463 788}%
\special{fp}%
% STR 2 0 3 0
% 3 580 500 580 600 5 0
% $Y_0$
\put(1.4500,-1.5000){\makebox(0,0){$Y_0$}}%
% STR 2 0 3 0
% 3 300 1110 300 1210 5 0
% $Y$
\put(0.7500,-3.0250){\makebox(0,0){$Y$}}%
% SPLINE 1 0 3 0
% 3 1790 1410 1560 1330 1330 1410
% 
\special{pn 13}%
\special{pa 1762 1388}%
\special{pa 1733 1373}%
\special{pa 1703 1358}%
\special{pa 1673 1344}%
\special{pa 1643 1332}%
\special{pa 1614 1322}%
\special{pa 1584 1314}%
\special{pa 1554 1311}%
\special{pa 1524 1310}%
\special{pa 1495 1313}%
\special{pa 1465 1320}%
\special{pa 1436 1329}%
\special{pa 1405 1341}%
\special{pa 1376 1355}%
\special{pa 1346 1370}%
\special{pa 1316 1385}%
\special{pa 1310 1388}%
\special{sp}%
% SARROW 1 0 3 1
% 2 1337 1407 1330 1410
% 
\special{pn 13}%
\special{pa 1316 1385}%
\special{pa 1310 1388}%
\special{fp}%
\special{sh 1}%
\special{pa 1310 1388}%
\special{pa 1377 1380}%
\special{pa 1358 1368}%
\special{pa 1362 1344}%
\special{pa 1310 1388}%
\special{fp}%
% SPLINE 1 0 3 0
% 4 1800 610 1520 700 1390 960 1390 960
% 
\special{pn 13}%
\special{pa 1772 601}%
\special{pa 1739 606}%
\special{pa 1704 612}%
\special{pa 1672 618}%
\special{pa 1639 625}%
\special{pa 1608 632}%
\special{pa 1578 642}%
\special{pa 1551 655}%
\special{pa 1524 669}%
\special{pa 1501 686}%
\special{pa 1480 706}%
\special{pa 1461 729}%
\special{pa 1444 753}%
\special{pa 1429 781}%
\special{pa 1416 810}%
\special{pa 1403 840}%
\special{pa 1392 872}%
\special{pa 1381 904}%
\special{pa 1372 937}%
\special{pa 1369 945}%
\special{sp}%
% SARROW 1 0 3 1
% 2 1393 951 1390 960
% 
\special{pn 13}%
\special{pa 1372 937}%
\special{pa 1369 945}%
\special{fp}%
\special{sh 1}%
\special{pa 1369 945}%
\special{pa 1408 889}%
\special{pa 1385 895}%
\special{pa 1371 876}%
\special{pa 1369 945}%
\special{fp}%
% SPLINE 1 0 3 0
% 4 1000 600 940 870 980 1010 980 1010
% 
\special{pn 13}%
\special{pa 985 591}%
\special{pa 973 622}%
\special{pa 961 653}%
\special{pa 951 684}%
\special{pa 941 714}%
\special{pa 934 746}%
\special{pa 928 776}%
\special{pa 925 807}%
\special{pa 925 837}%
\special{pa 927 868}%
\special{pa 934 898}%
\special{pa 941 929}%
\special{pa 951 958}%
\special{pa 963 989}%
\special{pa 965 995}%
\special{sp}%
% SARROW 1 0 3 1
% 2 978 1004 980 1010
% 
\special{pn 13}%
\special{pa 963 989}%
\special{pa 965 995}%
\special{fp}%
\special{sh 1}%
\special{pa 965 995}%
\special{pa 963 926}%
\special{pa 948 944}%
\special{pa 926 938}%
\special{pa 965 995}%
\special{fp}%
% LINE 1 1 3 0
% 4 1200 400 1200 1600 800 1200 2000 1200
% 
\special{pn 13}%
\special{pa 1182 394}%
\special{pa 1182 1575}%
\special{da 0.050}%
\special{pa 788 1182}%
\special{pa 1969 1182}%
\special{da 0.050}%
% LINE 2 0 3 0
% 32 1500 1200 1200 1500 1340 1420 1200 1560 1380 1440 1230 1590 1550 1330 1290 1590 1590 1350 1350 1590 1590 1410 1410 1590 1590 1470 1470 1590 1590 1530 1530 1590 1470 1350 1390 1430 1440 1200 1200 1440 1380 1200 1200 1380 1320 1200 1200 1320 1260 1200 1200 1260 1560 1200 1390 1370 1590 1230 1490 1330 1590 1290 1560 1320
% 
\special{pn 8}%
\special{pa 1477 1182}%
\special{pa 1182 1477}%
\special{fp}%
\special{pa 1319 1398}%
\special{pa 1182 1536}%
\special{fp}%
\special{pa 1359 1418}%
\special{pa 1211 1565}%
\special{fp}%
\special{pa 1526 1310}%
\special{pa 1270 1565}%
\special{fp}%
\special{pa 1565 1329}%
\special{pa 1329 1565}%
\special{fp}%
\special{pa 1565 1388}%
\special{pa 1388 1565}%
\special{fp}%
\special{pa 1565 1447}%
\special{pa 1447 1565}%
\special{fp}%
\special{pa 1565 1506}%
\special{pa 1506 1565}%
\special{fp}%
\special{pa 1447 1329}%
\special{pa 1369 1408}%
\special{fp}%
\special{pa 1418 1182}%
\special{pa 1182 1418}%
\special{fp}%
\special{pa 1359 1182}%
\special{pa 1182 1359}%
\special{fp}%
\special{pa 1300 1182}%
\special{pa 1182 1300}%
\special{fp}%
\special{pa 1241 1182}%
\special{pa 1182 1241}%
\special{fp}%
\special{pa 1536 1182}%
\special{pa 1369 1349}%
\special{fp}%
\special{pa 1565 1211}%
\special{pa 1467 1310}%
\special{fp}%
\special{pa 1565 1270}%
\special{pa 1536 1300}%
\special{fp}%
% LINE 2 0 3 0
% 32 1920 1200 1740 1380 1980 1200 1790 1390 1990 1250 1650 1590 1770 1410 1600 1580 1730 1390 1600 1520 1690 1370 1600 1460 1650 1350 1600 1400 1990 1310 1710 1590 1990 1370 1770 1590 1990 1430 1830 1590 1990 1490 1890 1590 1990 1550 1950 1590 1860 1200 1710 1350 1800 1200 1660 1340 1740 1200 1620 1320 1680 1200 1600 1280
% 
\special{pn 8}%
\special{pa 1890 1182}%
\special{pa 1713 1359}%
\special{fp}%
\special{pa 1949 1182}%
\special{pa 1762 1369}%
\special{fp}%
\special{pa 1959 1231}%
\special{pa 1625 1565}%
\special{fp}%
\special{pa 1743 1388}%
\special{pa 1575 1556}%
\special{fp}%
\special{pa 1703 1369}%
\special{pa 1575 1497}%
\special{fp}%
\special{pa 1664 1349}%
\special{pa 1575 1438}%
\special{fp}%
\special{pa 1625 1329}%
\special{pa 1575 1378}%
\special{fp}%
\special{pa 1959 1290}%
\special{pa 1684 1565}%
\special{fp}%
\special{pa 1959 1349}%
\special{pa 1743 1565}%
\special{fp}%
\special{pa 1959 1408}%
\special{pa 1802 1565}%
\special{fp}%
\special{pa 1959 1467}%
\special{pa 1861 1565}%
\special{fp}%
\special{pa 1959 1526}%
\special{pa 1920 1565}%
\special{fp}%
\special{pa 1831 1182}%
\special{pa 1684 1329}%
\special{fp}%
\special{pa 1772 1182}%
\special{pa 1634 1319}%
\special{fp}%
\special{pa 1713 1182}%
\special{pa 1595 1300}%
\special{fp}%
\special{pa 1654 1182}%
\special{pa 1575 1260}%
\special{fp}%
% LINE 2 0 3 0
% 24 1990 890 1690 1190 1990 830 1630 1190 1960 800 1600 1160 1900 800 1600 1100 1840 800 1600 1040 1780 800 1600 980 1720 800 1600 920 1660 800 1600 860 1990 950 1750 1190 1990 1010 1810 1190 1990 1070 1870 1190 1990 1130 1930 1190
% 
\special{pn 8}%
\special{pa 1959 876}%
\special{pa 1664 1172}%
\special{fp}%
\special{pa 1959 817}%
\special{pa 1605 1172}%
\special{fp}%
\special{pa 1930 788}%
\special{pa 1575 1142}%
\special{fp}%
\special{pa 1871 788}%
\special{pa 1575 1083}%
\special{fp}%
\special{pa 1812 788}%
\special{pa 1575 1024}%
\special{fp}%
\special{pa 1752 788}%
\special{pa 1575 965}%
\special{fp}%
\special{pa 1693 788}%
\special{pa 1575 906}%
\special{fp}%
\special{pa 1634 788}%
\special{pa 1575 847}%
\special{fp}%
\special{pa 1959 936}%
\special{pa 1723 1172}%
\special{fp}%
\special{pa 1959 995}%
\special{pa 1782 1172}%
\special{fp}%
\special{pa 1959 1054}%
\special{pa 1841 1172}%
\special{fp}%
\special{pa 1959 1113}%
\special{pa 1900 1172}%
\special{fp}%
% LINE 2 0 3 0
% 18 1520 400 1200 720 1580 400 1200 780 1590 450 1250 790 1590 510 1310 790 1590 570 1370 790 1460 400 1200 660 1400 400 1200 600 1340 400 1200 540 1280 400 1200 480
% 
\special{pn 8}%
\special{pa 1497 394}%
\special{pa 1182 709}%
\special{fp}%
\special{pa 1556 394}%
\special{pa 1182 768}%
\special{fp}%
\special{pa 1565 443}%
\special{pa 1231 778}%
\special{fp}%
\special{pa 1565 502}%
\special{pa 1290 778}%
\special{fp}%
\special{pa 1565 562}%
\special{pa 1349 778}%
\special{fp}%
\special{pa 1438 394}%
\special{pa 1182 650}%
\special{fp}%
\special{pa 1378 394}%
\special{pa 1182 591}%
\special{fp}%
\special{pa 1319 394}%
\special{pa 1182 532}%
\special{fp}%
\special{pa 1260 394}%
\special{pa 1182 473}%
\special{fp}%
% LINE 2 0 3 0
% 4 1590 690 1490 790 1590 750 1550 790
% 
\special{pn 8}%
\special{pa 1565 680}%
\special{pa 1467 778}%
\special{fp}%
\special{pa 1565 739}%
\special{pa 1526 778}%
\special{fp}%
\end{picture}%
\mbox{}\hfill\mbox{}\bigskip\bigskip\\
\mbox{}\hfill {\bf fig.\ 8}\hfill\mbox{}\bigskip\bigskip\\
\indent
$\pairof{\tilde{X},\tilde{Y},\tilde{\tau},\tilde{\phi},\tilde{\psi}}$ 
satisfies \xitemof{7-6} by \xitemof{boxplus} and \xitemof{dstar-1}. It 
satisfies \xitemof{7-7} by \xitemof{7-10}. Thus we are done by checking 
$\pairof{\tilde{X},\tilde{Y},\tilde{\tau},\tilde{\phi},\tilde{\psi}}$ 
also satisfies \xitemof{7-9a} and \xitemof{7-8}. 

For \xitemof{7-9a}, suppose that
$\pairof{\alpha,\beta}\in(\tilde{X}\setminus X)\times Y$ ($=X_0\times Y$). 
If $\pairof{\alpha,\beta}\in X_0\times Y'$ then
$\tilde{\tau}(\alpha,\beta)=\tau(\tilde{\phi}(\alpha),\tilde{\psi}(\beta))$
by \xitemof{7-10}. But
$\pairof{\tilde{\phi}(\alpha),\tilde{\psi}(\beta)}
\in (X\setminus X')\times Y'$ by definition of $\tilde{\phi}$ and
$\tilde{\psi}$. Hence, by \xitemof{7-9}, we have
$\tilde{\tau}(\alpha,\beta)=\tau(\tilde{\phi}(\alpha),\tilde{\psi}(\beta))=1$. 
If $\pairof{\alpha,\beta}\in X_0\times(Y\setminus Y')$ then 
$\tilde{\tau}(\alpha,\beta)=1$ by the ``otherwise'' clause of 
\xitemof{7-10}. 

For \xitemof{7-8}, it is enough to check that $\tilde{\tau}(\alpha,\beta)=1$ 
for all $\pairof{\alpha,\beta}\in (X'\cup X_0)\times Y_0$ with $\beta<\alpha$
by \xitemof{7-10} and \xitemof{7-9a}.   
For such $\pairof{\alpha,\beta}$, we have 
$\tilde{\tau}(\alpha,\beta)=\tau(\tilde{\phi}(\alpha),\tilde{\psi}(\beta))$ 
by \xitemof{7-10}. Suppose that
$\tilde{\tau}(\alpha,\beta)=0$. Then, since $\tau$ satisfies \xitemof{7-8}, 
we should have $\tilde{\phi}(\alpha)\leq\tilde{\psi}(\beta)$. By 
\xitemof{dstar-2}, we have $\beta>\tilde{\psi}(\beta)+\xi+\theta$. 
On the other hand, by \xitemof{boxplus}, we have
$\alpha=\tilde{\phi}^2(\alpha)\leq\tilde{\phi}(\alpha)+\xi+\theta$. 
It follows that 
\begin{xitemize}
\item[] 
	$\alpha\leq\tilde{\phi}(\alpha)+\xi+\theta
	\leq\tilde{\psi}(\beta)+\xi+\theta<\beta$. 
\end{xitemize}
This is a contradiction.
\qedofLemma
\qedskip

A quartet $p=\thecond{p}$ 
(not necessarily an element of $\poP_0$) with 
$D^p\in[\omega_2]^{\aleph_0}$ is said to be an {\em extendable condition} 
if $\pairof{X^p,Y^p,\tau^p,\phi^p_\xi,\psi^p_\xi}$ is a $\xi$-extendable 
semi-condition for all $\xi\in D^p$. 

For extendable conditions $p$, $q$ with 
\begin{xitemize}
	\item[] $p=\thecond{p}$,\ $q=\thecond{q}$, 
\end{xitemize}
we denote $p\leq_1 q$ if 
\begin{xitemize}
	\xitem[7-11] $X^p\supseteq X^q$, $Y^p\supseteq Y^q$,
	$\tau^p\supseteq\tau^q$, $D^p\supseteq D^q$, $\phi^p_\xi\supseteq \phi^q_\xi$ 
	and $\psi^p_\xi\supseteq\psi^q_\xi$ for all $\xi\in D^p$\,; 
\end{xitemize}	
and 
\begin{xitemize}
	\xitem[7-9b] $\tau^p\restr (X^p\setminus X^q)\times Y^q\equiv 1$. 
\end{xitemize}
Note that for $p$, $q\in\poP_0$, we have $p\leq_1 q$ if and only 
if $p\leq_{\poP_0}q$. 
\begin{Lemma}{\rm (Extension Lemma)} \label{extension-lemma}
	Suppose that 
	\begin{xitemize}
	\item[] $p=\thecond{p}$ 
	\end{xitemize}
	is an extendable condition for some
	$D^p\in[\omega_2]^{\aleph_0}$. Then there is a $q\in\poP_0$ with
	$q=\thecond{q}$ \st\ $D^q=D^p$ and $q\leq_1 p$. 

	Furthermore , if $p_0\in\poP_0$ is \st\ $p\leq_1 p_0$ 
	then we have $q\leq_{\poP_0}p_0$. 
\end{Lemma}
\prf
The second part of the lemma is clear once the condition $q$ as in the 
claim of the lemma is found since \xxitemof{(A)-8}{e} holds for such $q$ 
and $p_0$ since the relation $\leq_1$ is easily seen to be transitive. 

To construct the desired $q\in\poP_0$, let $\seqof{\xi_n}{n\in\omega}$ be 
an enumeration of $D^p$ \st\ each $\xi\in D^p$ appears infinitely often in 
the enumeration. 

First, construct 
$\pairof{X_n,Y_n,\tau_n,
	\seqof{\phi_{\xi,n},\,\psi_{\xi,n}}{\xi\in D^p}}$, $n\in\omega$ recursively \st\
\begin{xitemize}
	\xitem[ext-0] $\pairof{X_0,Y_0,\tau_0,
		\seqof{\phi_{\xi,0},\,\psi_{\xi,0}}{\xi\in D^p}}=p$,
	\xitem[ext-1] 
	$\pairof{X_{n+1},Y_{n+1},\tau_{n+1},\phi_{\xi_n,n+1},
	\psi_{\xi_n,n+1}}$ is the $\xi_n$-extendable semi-condition which is 
	constructed 	just as in \Lemmaof{xi-extendable} from the 
	$\xi_n$-extendable semi-condition
	$\pairof{X_{n},Y_{n},\tau_{n},\phi_{\xi_n,n},\psi_{\xi_n,n}}$. 
	\xitem[ext-2] $\phi_{\xi,n+1}=\phi_{\xi,n}$ and 
	$\psi_{\xi,n+1}=\psi_{\xi,n}$ for all $\xi\in D^p$ with $\xi\not=\xi_n$. 
\end{xitemize}
Along with the recursive construction above, it can be shown easily that 
$\pairof{X_{n},Y_{n},\tau_{n},\linebreak[2]\phi_{\xi,n},\psi_{\xi,n}}$ is 
a $\xi$-extendable semi-condition for all $n\in\omega$ and $\xi\in D^p$. Hence 
the construction in \xitemof{ext-1} is actually possible at each step. 

Let
\begin{xitemize}
\item[] $X^q=\bigcup_{n\in\omega}X_n$\,, $Y^q=\bigcup_{n\in\omega}X_n$\,,   
$\tau^q=\bigcup_{n\in\omega}\tau_n$\,,\\[\jot]    
$D^q=D^p$\ \ and\ \ $\phi^q_\xi=\bigcup_{n\in\omega}\phi_{\xi,n}$\,,
$\psi^q_\xi=\bigcup_{n\in\omega}\psi_{\xi,n}$\ \ for all $\xi\in D^q$.
\end{xitemize}
For all $\xi\in D^q$, there are infinitely many $n\in\omega$ 
\st\ $\xi_n=\xi$. For such $n$, $\phi_{\xi,n}$ is an involution on $X_n$ and 
$\psi_{\xi,n}$ is an involution on $Y_n$. It follows that $\phi^q_\xi$ is 
an involution on $X^q$ and $\psi^q_\xi$ is an involution on $Y^q$. 
Hence 
\begin{xitemize}
	\item[] $q=\thecond{q}$ 
\end{xitemize}
is a condition in $\poP_0$. 
Also we have 
\begin{xitemize}
\item[] $\tau^q\restr(X^q\setminus X^p)\times Y^p
	=\bigcup_{n\in\omega}\tau^q\restr(X_{n+1}\setminus X_n)\times Y^p\equiv 
	1$. 
\end{xitemize}
Thus, this $q$ is as desired. \qedofLemma
\begin{Lemma}{\rm (\CH)} \label{aleph-2-c.c.}
	$\poP_0$ satisfies the $\aleph_2$-c.c.\ 
\end{Lemma}
\prf Actually we shall show that $\poP_0$ satisfies a strong form 
of $\aleph_2$-Knaster property. 

Suppose that $p^\zeta\in\poP_0$ with $p^\zeta=\thecond{\zeta}$ for $\zeta\in\omega_2$. 
By the $\Delta$-System Lemma (\Thmof{delta-lemma}) and the Pigeon Hole 
Principle, there are a stationary $S\subseteq\omega_2$, $X$, $Y$,
$D\in[\omega_2]^{\aleph_0}$, $\mapping{\tau}{X\times Y}{2}$ and
$\mapping{\phi_\xi}{X}{X}$, $\mapping{\psi_\xi}{Y}{Y}$ for $\xi\in D$ \st\ 
\begin{xitemize}
	\xitem[cc-0] $X^\zeta$, $\zeta\in S$ form a $\Delta$-system with root $X$ and 
	$Y^\zeta$, $\zeta\in S$ form a $\Delta$-system with root $Y$\,;
	\xitem[cc-1] $\tau^\zeta\restr X\times Y=\tau$ for all $\zeta\in S$\,;
	\xitem[cc-2] $D^\zeta$, $\zeta\in S$ form a $\Delta$-system with root $D$\,;
	\xitem[cc-3] $\phi^\zeta_\xi\restr X=\phi_\xi$ and
	$\psi^\zeta_\xi\restr Y=\psi_\xi$ for all $\zeta\in S$ and $\xi\in D$\,;
	\xitem[cc-4] $\tau^\zeta\restr(X^\zeta\setminus X)\times Y\equiv 1$ for all
	$\zeta\in S$\,.
\end{xitemize}
Note that \xitemof{cc-1} is possible since, by \CH, there are at most 
$\cardof{\fnsp{X\times Y}{2}}\leq 2^{\aleph_0}=\aleph_1<\aleph_2$ many 
possible values of $\tau^\zeta\restr X\times Y$. \xitemof{cc-3} is possible 
since, by \CH, \xitemof{(A)-4} and countability of $D$, there are at most
$\aleph_1$ possible values of $\seqof{\phi^\zeta_\xi}{\xi\in D}$ and
$\seqof{\psi^\zeta_\xi}{\xi\in D}$. \xitemof{cc-4} is possible by 
\xitemof{(A)-7} and since we can choose $S$ \st\
$\min(X^\zeta\setminus X)>\sup(Y)$ for all $\zeta\in S$. 

Now suppose $\zeta$, $\eta\in S$ with $\zeta<\eta$. We show that $p^\zeta$ 
and $p^\eta$ are compatible. Let $X^p=X^\zeta\cup X^\eta$,
$Y^p=Y^\zeta\cup Y^\eta$ and $D^p=D^\zeta\cup D^\eta$. For $\xi\in D^p$, let
$\mapping{\phi^p_\xi}{X^p}{X^p}$ and $\mapping{\psi^p_\xi}{Y^p}{Y^p}$ be 
defined by 
\begin{xitemize}
\xitem[] $\phi^p_\xi=\left\{ 
	\begin{array}[c]{@{}ll}
		\phi^\zeta_\xi, &\mbox{if }\xi\in D^\zeta\setminus D,\\[\jot]
		\phi^\zeta_\xi\cup \phi^\eta_\xi, &\mbox{if }\xi\in D,\\[\jot]
		\phi^\eta_\xi, &\mbox{otherwise}
	\end{array}\right.$
	\quad\ and\quad\  
	$\psi^p_\xi=\left\{ 
	\begin{array}[c]{@{\ }ll}
		\psi^\zeta_\xi, &\mbox{if }\xi\in D^\zeta\setminus D,\\[\jot]
		\psi^\zeta_\xi\cup \psi^\eta_\xi, &\mbox{if }\xi\in D,\\[\jot]
		\psi^\eta_\xi, &\mbox{otherwise.}
	\end{array}\right.$
\end{xitemize}
Finally let $\mapping{\tau^p}{X^p\times Y^p}{2}$ be \st\ 
\begin{xitemize}
\xitem[cc-5] $\tau^p(\alpha,\beta)=\left\{ 
	\begin{array}[c]{@{\ }ll}
		\tau^\zeta(\alpha,\beta), &\mbox{if }\pairof{\alpha,\beta}\in 
		X^\zeta\times Y^\zeta\\[\jot]
		\tau^\eta(\alpha,\beta), &\mbox{else if }\pairof{\alpha,\beta}\in 
		X^\eta\times Y^\eta\\[\jot]
		1, &\mbox{otherwise.}
	\end{array}\right.$
\end{xitemize}
for all $\pairof{\alpha,\beta}\in X^p\times Y^p$. 

It is easy to see that $p=\thecond{p}$ is an extendable condition and 
$p\leq_1 p^\zeta$, $p^\eta$. In particular, \xitemof{7-9b} for
$p\leq_1 p^\zeta$ and $p\leq_1 p^\eta$ holds because of \xitemof{cc-4} 
and ``otherwise'' clause of \xitemof{cc-5}. 

By Extension Lemma (\Lemmaof{extension-lemma}), there is a $q\in\poP_0$ 
with $q\leq_1 p$. Hence, by the second half of the lemma, it follows that 
$q\leq_{\poP_0}p^\zeta$, $p^\eta$. \qedofLemma
\qedskip

A modification of the $\Delta$-system argument in the proof of 
\Lemmaabove\ is also used to prove the following:
\begin{Lemma}{\rm (CH)}
	\label{homogeneity-theorem}
	$\poP_0$ forces \xitemof{Omega-4}. 
\end{Lemma}
\prf 	
We show that $\poP_0$ forces the following: 
\begin{xitemize}
	\xitem[homo-0] For any stationary $S\subseteq E^{\omega_2}_{\omega_1}$ 
	and $A_\zeta$, $B_\zeta\in[\omega_2]^{\aleph_0}$ for $\zeta\in S$, 
	there is a stationary $T\subseteq S$ \st\ for any $n\in\omega$ and 
	pairwise distinct $\zeta_i$, $\eta_i\in T$, $i\in n$, there is 
	$\xi<\omega_2$ \st\ $\dot\phi_\xi\imageof A_{\zeta_i}= A_{\eta_i}$ and
	$\dot\psi_\xi\imageof B_{\zeta_i}= B_{\eta_i}$ for all $i\in n$. 
\end{xitemize}

Note that, 
by $\sigma$-closedness and $\aleph_2$-c.c.\  of $\poP_0$ (proved in 
\Lemmaof{aleph-2-c.c.}), $\omega_1$ 
and $\omega_2$ in generic extensions by $\poP_0$ remain $\omega_1$ and
$\omega_2$. 

Suppose that $\dot{S}$ is a
$\poP_0$-name of a stationary subset of $E^{\omega_2}_{\omega_1}$. Let 
$\seqof{\dot{A}_\zeta}{\zeta\in\dot{S}}$ and 
$\seqof{\dot{B}_\zeta}{\zeta\in\dot{S}}$ be $\poP_0$-names of sequences 
of countable subsets of $\omega_2$. Let 
\begin{xitemize}
\item[] $\tilde{S}
	=\setof{\zeta\in E^{\omega_2}_{\omega_1}}{
			\notforces{\poP_0}{\zeta\not\in\dot{S}}}$.
\end{xitemize}
Then we have $\forces{\poP_0}{\dot{S}\subseteq\tilde{S}}$ and hence 
$\tilde{S}$ is a stationary subset of $E^{\omega_2}_{\omega_1}$. 

Since $\poP_0$ is $\sigma$-closed, we can find $p_\zeta\in\poP_0$ and 
$A_\zeta$, $B_\zeta\in[\omega_2]^{\aleph_0}$ \st\ 
\begin{xitemize}
	\xitem[po-0] $p_\zeta=\thecond{\zeta}$ and 
	\xitem[po-1] $p_\zeta\forces{\poP_0}{\zeta\in\dot{S}\,,\ 
	\dot{A}_\zeta=A_\zeta \mbox{ and }\dot{B}_\zeta=B_\zeta}$ 
\end{xitemize}
for all $\zeta\in\tilde{S}$. 

\Wolog, we may assume that
\begin{xitemize}
	\xitem[po-2] $A_\zeta\subseteq X^\zeta$ and
	$B_\zeta\subseteq Y^\zeta$. 
\end{xitemize}
By $\Delta$-System Lemma (\Thmof{delta-lemma}) and the Pigeon Hole 
Principle, there are a stationary
$\tilde{S}_0\subseteq \tilde{S}$, 
 $X$, $Y$,
$D\in[\omega_2]^{\aleph_0}$, $\mapping{\tau}{X\times Y}{2}$ and
$\mapping{\phi_\xi}{X}{X}$, $\mapping{\psi_\xi}{Y}{Y}$ for $\xi\in D$ \st\ 
\begin{xitemize}
	\xitem[po-3] $X^\zeta$, $\zeta\in \tilde{S}_0$ form a $\Delta$-system with 
	root $X$ and  
	$Y^\zeta$, $\zeta\in\tilde{S}_0$ form a $\Delta$-system with root $Y$\,;
	\xitem[po-3A] $\sup(Y)<\min(X^\zeta\setminus X)$ for all $\zeta\in\tilde{S}_0$\,;
	\xitem[po-4] $\tau^\zeta\restr X\times Y=\tau$ for all $\zeta\in\tilde{S}_0$\,;
	\xitem[po-5] $D^\zeta$, $\zeta\in\tilde{S}_0$ form a $\Delta$-system with root $D$\,;
	\xitem[po-6] $\phi^\zeta_\xi\restr X=\phi_\xi$ and
	$\psi^\zeta_\xi\restr Y=\psi_\xi$ for all $\zeta\in\tilde{S}_0$ and
	$\xi\in D$ \,;
	\xitem[po-7] $\tau^\zeta\restr(X^\zeta\setminus X)\times Y\equiv 1$ for all
	$\zeta\in\tilde{S}_0$  (this follows from \xitemof{po-3A} and \xitemof{(A)-7})\,;
	\xitem[po-8] $X^\zeta$, $\zeta\in\tilde{S}_0$ are order-isomorphic and 
	$Y^\zeta$, $\zeta\in\tilde{S}_0$ are order-isomorphic; 
\end{xitemize}
Note that the 
order-isomorphisms of $X^\zeta$'s and $Y^\zeta$'s do not 
move elements of $X$ and $Y$, respectively. 
\begin{xitemize}
	\xitem[po-9] the order-isomorphism sending $X^\zeta$ to $X^\eta$ sends 	
	$\tau^\zeta\restr((X^\zeta\setminus X)\times Y)$ to
	$\tau^\eta\restr((X^\eta\setminus X)\times Y)$ while the 
	order-isomorphism sending $Y^\zeta$ to $Y^\eta$ sends 
	$\tau^\zeta\restr(X\times(Y^\zeta\setminus Y))$ to
	$\tau^\eta\restr(X\times(Y^\eta\setminus Y))$. These order-isomorphisms 
	together send
	$\tau^\zeta\restr((X^\zeta\setminus X)\times(Y^\zeta\setminus Y))$ to
	$\tau^\eta\restr((X^\eta\setminus X)\times(Y^\eta\setminus Y))$;
	\xitem[po-10] the order-isomorphism sending $X^\zeta$ to $X^\eta$ sends
	$A_\zeta$ to $A_\eta$, 
   and the order-isomorphism sending $Y^\zeta$ to $Y^\eta $ sends $B_\zeta$ 
	 to $B_\eta$. 
\end{xitemize}
Note that $\bar{p}=\thecond{}$ is a condition in $\poP_0$ and 
$p_\zeta\leq_{\poP_0}\bar{p}$ for all $\zeta\in\tilde{S}_0$ (the condition 
\xxitemof{(A)-8}{e} for $\bar{p}$ and $p_\zeta$ holds by \xitemof{po-7}). 

Let $\dot{T}$ be a $\poP_0$-name \st\ 
\begin{xitemize}
	\xitem[] $\forces{\poP_0}{\dot{T}=\setof{\zeta\in\tilde{S}_0}{p_\zeta\in\dot{G}\,}}$
\end{xitemize}
where $\dot{G}$ is the standard $\poP_0$-name of the generic set. 
\begin{Claim}
	\label{homo-1}
	$\bar{p}\forces{\poP_0}{\dot{T}\mbox{ is a stationary subset of }\omega_2}$.
\end{Claim}
\prfofClaim
Since $\poP_0$ satisfies the $\aleph_2$-c.c.\ by \Lemmaof{aleph-2-c.c.}, 
for any $\poP_0$-name $\dot{C}$ of a club subset of $\omega_2$, there is a 
club subset $C$ of $\omega_2$ (in the ground model) \st\
$\forces{\poP_0}{C\subseteq\dot{C}}$. Hence it is enough to show the 
following:
\begin{xitemize}
	\xitem[po-11] For any $q\leq_{\poP_0}\bar{p}$ and any club subset $C$ of
	$\omega_2$, there are $p\leq_{\poP_0}q$ and $\zeta\in C\cap\tilde{S}_0$ 
	\st\ $p\leq_{\poP_0}p_\zeta$. 
\end{xitemize}
To show \xitemof{po-11}, let $q=\thecond{q}$ and let 
$\zeta\in C\cap\tilde{S}_0$ be \st\ 
\begin{xitemize}
	\xitem[po-12] $(X^\zeta\setminus X)\cap X^q=\emptyset$\,,
$(Y^\zeta\setminus Y)\cap Y^q=\emptyset$ and 
$(D^\zeta\setminus D)\cap D^q=\emptyset$. 
\end{xitemize}
This is possible by \xitemof{po-3} and since $C\cap\tilde{S}_0$ is 
stationary. 

Let $X^*=X^q\cup X^\zeta$, $Y^*=Y^q\cup Y^\zeta$, $D^*=D^q\cup D^\zeta$. 
For $\xi\in D^*$, let $\phi^*_\xi$ and $\psi^*_\xi$ be partial functions 
from $X^*$ to $X^*$ and from $Y^*$ to $Y^*$ respectively defined by
\begin{xitemize}
	\xitem[po-13] $\phi^*_\xi=\left\{ 
	\begin{array}[c]{@{}ll}
		\phi^q_\xi, &\mbox{if }\xi\in D^q\setminus D,\\[\jot]
		\phi^q_\xi\cup \phi^\zeta_\xi, &\mbox{if }\xi\in D,\\[\jot]
		\phi^\zeta_\xi, &\mbox{otherwise}
	\end{array}\right.$
	\quad\ and\quad\  
	$\psi^*_\xi=\left\{ 
	\begin{array}[c]{@{\ }ll}
		\psi^q, &\mbox{if }\xi\in D^q\setminus D,\\[\jot]
		\psi^q_\xi\cup \psi^\zeta_\xi, &\mbox{if }\xi\in D,\\[\jot]
		\psi^\zeta_\xi, &\mbox{otherwise.}
	\end{array}\right.$
\end{xitemize}
Finally, let $\mapping{\tau^*}{X^*\times Y^*}{2}$ be defined by 
\begin{xitemize}
	\xitem[po-14] $\tau^*(\alpha,\beta)=\left\{
	\begin{array}[c]{@{\ }ll}
		\tau^q(\alpha,\beta), &\mbox{if }\pairof{\alpha,\beta}\in X^q\times Y^q,
		\\[\jot]
		\tau^\zeta(\alpha,\beta), 
		&\mbox{else if }\pairof{\alpha,\beta}\in X^\zeta\times Y^\zeta,\\[\jot]
		1, &\mbox{otherwise}
	\end{array}\right.
	$
\end{xitemize}
for $\pairof{\alpha,\beta}\in X^*\times Y^*$. 

Then $p^*=\thecond{*}$ is an extendable condition and we have
$p^*\leq_1 q$, $p_\zeta$ : \xitemof{7-9b} for $p^*$ and $p_\zeta$ follows 
from 
\begin{xitemize}
\item[] $\tau^*\restr((X^*\setminus X^\zeta)\times Y^\zeta)=
	\tau^*\restr((X^q\setminus X)\times (Y^\zeta\setminus Y))\cup
	\tau^*\restr((X^q\setminus X)\times Y)\equiv 1$
\end{xitemize}
where we have
$\tau^*\restr((X^q\setminus X)\times (Y^\zeta\setminus Y))\equiv 1$ by the 
definition \xitemof{po-14} of $\tau^*$ and
$\tau^*\restr((X^q\setminus X)\times Y)\equiv 1$ by $q\leq_{\poP_0}\bar{p}$ 
(in particular, by the condition \xxitemof{(A)-8}{e} in the definition of 
$\leq_{\poP_0}$). 

By Extension Lemma (\Lemmaof{extension-lemma}) it follows that there is 
$p\in\poP_0$ with $p\leq_1 p^*$ and hence $p\leq_{\poP_0} q$, $p_\zeta$. 
\qedofClaim
%%%
\begin{Claim}
	$\bar{p}$ forces that $\dot{T}$ is as in \xitemof{homo-0} for
	$\seqof{\dot{A}_\zeta}{\zeta\in\dot{S}}$ and 
	$\seqof{\dot{B}_\zeta}{\zeta\in\dot{S}}$. 
\end{Claim}
\prfofClaim By \Claimof{homo-1} it is enough to prove the following:
\begin{xitemize}
	\xitem[homo-2] For any $q\leq_{\poP_0}\bar{p}$ and $n\in\omega$, if
	$\zeta_i$, $\eta_i\in\tilde{S}_0$ are pairwise distinct and
	$q\forces{\poP_0}{\zeta_i,\,\eta_i\in\dot{T}\mbox{ for }i\in n}$, then 
	there is $p\leq_{\poP_0}q$ with 
	\begin{xxitemize}
		\item[] $p=\thecond{p}$ 
	\end{xxitemize}
and $\xi_0\in D^p$ \st\ 
	$\mapping{\phi^p_{\xi_0}\restr X^{\zeta_i}}{X^{\zeta_i}}{X^{\eta_i}}$ and 
	$\mapping{\psi^p_{\xi_0}\restr Y^{\zeta_i}}{Y^{\zeta_i}}{Y^{\eta_i}}$ are 
	order-isomorphisms for all $i<n$. 
\end{xitemize}
\Wolog\ we may assume that 
\begin{xitemize}
	\xitem[homo-3] $q\leq_{\poP_0}p_{\zeta_i}$, $p_{\eta_i}$ for 
all $i<n$. 
\end{xitemize}
Let 
\begin{xitemize}
\item[] $q=\thecond{q}$. 
\end{xitemize}
Take $\xi_0\in\omega_2\setminus (D^q\cup\sup(X^q)\cup\sup(Y^q))$ and let
$D^*=D^q\cup\ssetof{\xi_0}$. 
Let $X^*=X^q$,  $Y^*=Y^q$  and
$\tau^*=\tau^q$. Let 
$\mapping{\phi^*_{\xi_0}}{\bigcup_{i<n}X^{\zeta_i}\cup\bigcup_{i<n}X^{\eta_i}}%
{\bigcup_{i<n}X^{\zeta_i}\cup\bigcup_{i<n}X^{\eta_i}}$ 
be the involution sending $X^{\zeta_i}$ order-isomorphically 
to $X^{\eta_i}$ and vice versa for all $i<n$ and 
$\mapping{\psi^*_{\xi_0}}{\bigcup_{i<n}Y^{\zeta_i}\cup\bigcup_{i<n}Y^{\eta_i}}%
{\bigcup_{i<n}Y^{\zeta_i}\cup\bigcup_{i<n}Y^{\eta_i}}$ 
be the involution sending $Y^{\zeta_i}$ order-isomorphically 
to $Y^{\eta_i}$ and vice versa for all $i<n$. Let 
$\phi^*_\xi=\phi^q_\xi$ and $\psi^*_\xi=\psi^q_\xi$ for $\xi\in D^q$.  

Then 
$p^*=\thecond{*}$ is an extendable condition with $p^*\leq_1 q$. 
To see this,  
we have to check $\pairof{X^*,Y^*,\tau^*,\phi^*_{\xi_0},\psi^*_{\xi_0}}$ 
satisfies \xitemof{7-7} and \xitemof{7-9}.
But this follows from \xitemof{po-7}, \xitemof{homo-3} 
and \xitemof{po-10}.

By Extension Lemma (\Lemmaof{extension-lemma}) there is $p\in\poP_0$ with
$p\leq_1 p^*$. Clearly $p$ forces that $\xi_0$ as above satisfies 
\xitemof{homo-0} together with $\dot{T}$, 
	$\seqof{\dot{A}_\zeta}{\zeta\in\dot{S}}$ and
	$\seqof{\dot{B}_\zeta}{\zeta\in\dot{S}}$.
\qedofClaim
\qedskip

Since the argument above can be repeated below arbitrary element 
of $\poP_0$, it follows that $\poP_0$ forces \xitemof{homo-0}. 
\qedofLemma
\section{A summary of consistency results and some open problems}
\label{zusammenfassung-und-offene-frage}
The following is a summary of consistency results 
in connection with the combinatorial principles in {fig.\,5} where  
\assertof{1} $\sim$ 
\assertof{7} below correspond to the separation lines \assertof{1} $\sim$ 
\assertof{7} drawn in {fig.\,9}. 
\medskip

(1): By adding random reals. More precisely, 
start from a model $V$ of $\CH$ and force with (the positive elements 
of) the measure algebra $\poB$ of, say, Maharam type $\aleph_2$. 
$\poB$ can be seen as a (measure theoretic) product of random 
forcing and inherits thus some of the homogeneity property of finite support 
product. This is used to prove $\donum=\aleph_1$ in the generic extension. It 
is also well-known that the ground model functions from $\omega$ 
to $\omega$ dominate the functions in a generic extension by a measure 
algebra. Hence we have $\dominatingno=\aleph_1$ in the model. 
K.\ Kunen proved that there is a $\kappa$-Lusin gap for an uncountable $\kappa$ 
in such a model. On the other hand, I.\ Juh\'asz, L.\ Soukup and Z.\ 
Szentmikl\'ossy proved in \cite{ju-so-sz} that there is no $\aleph_2$-Lusin 
gap under $\Cs(\aleph_2)$. This proves that $\Cs(\aleph_2)$ does not hold 
in the generic extension. 

This observation may be also interpreted as pinning down of the difference in 
the extent of 
homogeneity of product forcing and the forcing by measure theoretic 
products in terms of whether the principle $\Cs(\aleph_2)$ holds. 

(2): A model constructed by J.\ Brendle and T.\ LaBerge in 
\cite{brendle-laberge} realizes this separation. 

(3): By the model in Theorem 3.8 of I.\ Juh\'asz and K.\ Kunen 
\cite{juhasz-kunen} in which 
$\Cs(\aleph_2)$ and $\donum>\aleph_1$ hold. The model is obtained by a 
finite support product of $\aleph_2$ \pos\ of cardinality $\aleph_1$ 
starting from a model of \CH. From 
this, it  
follows easily that $\boundingno^*=\aleph_1$ and $\dominatingno=\aleph_2$. 

(4): By adding Cohen reals. More exactly, start from a model $V$ of $\CH$ 
and then add, say, $\aleph_2$ Cohen reals (by $\Fn(\aleph_2,2)$). Then by 
\Corof{fs-prod-IP},\,\assertof{c} we have $\IP(\aleph_2,\aleph_1)$ in the generic 
extension. Just as in (3), we have 
$\dominatingno=\aleph_2$ in such a generic extension and it is shown in 
S.\ Fuchino, S.\ Koppelberg and S.\ Shelah \cite{fu-ko-sh} that \WFN\ holds there. 

(5): By a model of Hechler. 

(6): By \Thmof{new-model}. 

(7): By \Corof{chang-conj-then-hechler}. 

(8): By \Thmof{prikry-silver}. See \cite{geschke} for the proof of
$\forces{\poS^\kappa}{\neg\WFN}$.

\mbox{}\\
\mbox{}\hspace{-1cm}
%%\input{coloring-x-fig3A}
%WinTpicVersion3.08
\unitlength 1cm
\begin{picture}( 15.0750, 10.0750)(  0.7000,-10.0125)
% VECTOR 0 0 3 0
% 2 2126 2873 1585 2873
% 
\special{pn 20}%
\special{pa 2093 2828}%
\special{pa 1561 2828}%
\special{fp}%
\special{sh 1}%
\special{pa 1561 2828}%
\special{pa 1626 2848}%
\special{pa 1613 2828}%
\special{pa 1626 2809}%
\special{pa 1561 2828}%
\special{fp}%
% VECTOR 0 0 3 0
% 2 3387 2873 2846 2873
% 
\special{pn 20}%
\special{pa 3334 2828}%
\special{pa 2802 2828}%
\special{fp}%
\special{sh 1}%
\special{pa 2802 2828}%
\special{pa 2868 2848}%
\special{pa 2854 2828}%
\special{pa 2868 2809}%
\special{pa 2802 2828}%
\special{fp}%
% VECTOR 0 0 3 0
% 2 4649 2873 4108 2873
% 
\special{pn 20}%
\special{pa 4576 2828}%
\special{pa 4044 2828}%
\special{fp}%
\special{sh 1}%
\special{pa 4044 2828}%
\special{pa 4110 2848}%
\special{pa 4096 2828}%
\special{pa 4110 2809}%
\special{pa 4044 2828}%
\special{fp}%
% VECTOR 0 0 3 0
% 2 3743 2130 3743 2670
% 
\special{pn 20}%
\special{pa 3685 2097}%
\special{pa 3685 2628}%
\special{fp}%
\special{sh 1}%
\special{pa 3685 2628}%
\special{pa 3704 2563}%
\special{pa 3685 2576}%
\special{pa 3665 2563}%
\special{pa 3685 2628}%
\special{fp}%
% VECTOR 0 0 3 0
% 2 3743 1318 3743 1859
% 
\special{pn 20}%
\special{pa 3685 1298}%
\special{pa 3685 1830}%
\special{fp}%
\special{sh 1}%
\special{pa 3685 1830}%
\special{pa 3704 1764}%
\special{pa 3685 1778}%
\special{pa 3665 1764}%
\special{pa 3685 1830}%
\special{fp}%
% STR 2 0 3 0
% 3 3880 1900 3880 1990 5 0
% $\donum=\aleph_1$
\put(9.7000,-4.9750){\makebox(0,0){$\donum=\aleph_1$}}%
% STR 2 0 3 0
% 3 3740 1090 3740 1180 5 0
% $\HP(\aleph_2)$
\put(9.3500,-2.9500){\makebox(0,0){$\HP(\aleph_2)$}}%
% STR 2 0 3 0
% 3 3747 2782 3747 2873 5 0
% $\boundingno^h=\aleph_1$
\put(9.3675,-7.1825){\makebox(0,0){$\boundingno^h=\aleph_1$}}%
% STR 2 0 3 0
% 3 2486 2782 2486 2873 5 0
% $\boundingno^\uparrow=\aleph_1$
\put(6.2150,-7.1825){\makebox(0,0){$\boundingno^\uparrow=\aleph_1$}}%
% STR 2 0 3 0
% 3 1270 2782 1270 2873 5 0
% $\boundingno=\aleph_1$
\put(3.1750,-7.1825){\makebox(0,0){$\boundingno=\aleph_1$}}%
% VECTOR 2 0 3 0
% 2 4986 956 4986 2676
% 
\special{pn 8}%
\special{pa 4908 941}%
\special{pa 4908 2634}%
\special{fp}%
\special{sh 1}%
\special{pa 4908 2634}%
\special{pa 4928 2568}%
\special{pa 4908 2582}%
\special{pa 4888 2568}%
\special{pa 4908 2634}%
\special{fp}%
% STR 2 0 3 0
% 3 4990 2792 4990 2882 5 0
% $\boundingno^*=\aleph_1$
\put(12.4750,-7.2050){\makebox(0,0){$\boundingno^*=\aleph_1$}}%
% STR 2 0 3 0
% 3 4986 671 4986 762 5 0
% $\WFN$
\put(12.4650,-1.9050){\makebox(0,0){$\WFN$}}%
% VECTOR 0 0 3 0
% 2 3549 1305 3107 1855
% 
\special{pn 20}%
\special{pa 3494 1285}%
\special{pa 3059 1826}%
\special{fp}%
\special{sh 1}%
\special{pa 3059 1826}%
\special{pa 3115 1788}%
\special{pa 3091 1785}%
\special{pa 3084 1763}%
\special{pa 3059 1826}%
\special{fp}%
% VECTOR 0 0 3 0
% 2 3194 2102 3583 2689
% 
\special{pn 20}%
\special{pa 3144 2069}%
\special{pa 3527 2647}%
\special{fp}%
\special{sh 1}%
\special{pa 3527 2647}%
\special{pa 3507 2581}%
\special{pa 3499 2604}%
\special{pa 3475 2603}%
\special{pa 3527 2647}%
\special{fp}%
% STR 2 0 3 0
% 3 3117 1882 3117 1972 5 0
% $\Cs(\aleph_2)$
\put(7.7925,-4.9300){\makebox(0,0){$\Cs(\aleph_2)$}}%
% VECTOR 2 0 3 0
% 2 4874 966 3266 1857
% 
\special{pn 8}%
\special{pa 4798 951}%
\special{pa 3215 1828}%
\special{fp}%
\special{sh 1}%
\special{pa 3215 1828}%
\special{pa 3282 1813}%
\special{pa 3261 1803}%
\special{pa 3263 1779}%
\special{pa 3215 1828}%
\special{fp}%
% VECTOR 0 0 3 0
% 2 5840 2870 5309 2870
% 
\special{pn 20}%
\special{pa 5749 2825}%
\special{pa 5226 2825}%
\special{fp}%
\special{sh 1}%
\special{pa 5226 2825}%
\special{pa 5292 2845}%
\special{pa 5278 2825}%
\special{pa 5292 2806}%
\special{pa 5226 2825}%
\special{fp}%
% STR 2 0 3 0
% 3 6190 2790 6190 2870 5 0
% $\dominatingno=\aleph_1$
\put(15.4750,-7.1750){\makebox(0,0){$\dominatingno=\aleph_1$}}%
% LINE 2 2 3 0
% 6 6310 1640 3430 1640 3430 1640 3430 2440 3430 2440 1030 2440
% 
\special{pn 8}%
\special{pa 6211 1615}%
\special{pa 3376 1615}%
\special{dt 0.045}%
\special{pa 3376 1615}%
\special{pa 3376 2402}%
\special{dt 0.045}%
\special{pa 3376 2402}%
\special{pa 1014 2402}%
\special{dt 0.045}%
% STR 2 0 3 0
% 3 840 2360 840 2440 5 0
% (1)
\put(2.1000,-6.1000){\makebox(0,0){(1)}}%
% STR 2 0 3 0
% 3 3140 3700 3140 3780 5 0
% (2)
\put(7.8500,-9.4500){\makebox(0,0){(2)}}%
% LINE 1 2 3 0
% 8 2340 1696 3380 1696 3380 1696 3380 2192 3380 2192 5536 2192 5536 2192 5536 3332
% 
\special{pn 13}%
\special{pa 2304 1670}%
\special{pa 3327 1670}%
\special{dt 0.045}%
\special{pa 3327 1670}%
\special{pa 3327 2158}%
\special{dt 0.045}%
\special{pa 3327 2158}%
\special{pa 5449 2158}%
\special{dt 0.045}%
\special{pa 5449 2158}%
\special{pa 5449 3280}%
\special{dt 0.045}%
% STR 2 0 3 0
% 3 2148 1616 2148 1696 5 0
% (3)
\put(5.3700,-4.2400){\makebox(0,0){(3)}}%
% LINE 1 1 3 0
% 8 4388 3296 4298 1526 4298 1526 3498 1526 3498 1526 3498 2496 3498 2496 1968 2496
% 
\special{pn 13}%
\special{pa 4319 3245}%
\special{pa 4231 1502}%
\special{da 0.050}%
\special{pa 4231 1502}%
\special{pa 3443 1502}%
\special{da 0.050}%
\special{pa 3443 1502}%
\special{pa 3443 2457}%
\special{da 0.050}%
\special{pa 3443 2457}%
\special{pa 1938 2457}%
\special{da 0.050}%
% LINE 2 1 3 0
% 4 5730 3300 5730 3300 5730 260 5730 260
% 
\special{pn 8}%
\special{pa 5640 3249}%
\special{pa 5640 3249}%
\special{da 0.050}%
\special{pa 5640 256}%
\special{pa 5640 256}%
\special{da 0.050}%
% LINE 2 1 3 0
% 2 5730 420 5730 3620
% 
\special{pn 8}%
\special{pa 5640 414}%
\special{pa 5640 3563}%
\special{da 0.050}%
% STR 2 0 3 0
% 3 4420 3384 4420 3464 5 0
% (6)
\put(11.0500,-8.6600){\makebox(0,0){(6)}}%
% STR 2 0 3 0
% 3 5740 3710 5740 3790 5 0
% (4)
\put(14.3500,-9.4750){\makebox(0,0){(4)}}%
% LINE 2 1 3 0
% 2 1810 3520 1810 1200
% 
\special{pn 8}%
\special{pa 1782 3465}%
\special{pa 1782 1182}%
\special{da 0.050}%
% STR 2 0 3 0
% 3 1810 3700 1810 3800 5 0
% (5)
\put(4.5250,-9.5000){\makebox(0,0){(5)}}%
% LINE 2 1 3 0
% 6 5620 3910 5620 1250 5620 1250 4720 1250 4720 1250 4720 30
% 
\special{pn 8}%
\special{pa 5532 3849}%
\special{pa 5532 1231}%
\special{da 0.050}%
\special{pa 5532 1231}%
\special{pa 4646 1231}%
\special{da 0.050}%
\special{pa 4646 1231}%
\special{pa 4646 30}%
\special{da 0.050}%
% STR 2 0 3 0
% 3 5620 3990 5620 4090 5 0
% (7)
\put(14.0500,-10.2250){\makebox(0,0){(7)}}%
% LINE 2 2 3 0
% 4 3140 2640 3140 3540 3140 2650 2140 1050
% 
\special{pn 8}%
\special{pa 3091 2599}%
\special{pa 3091 3485}%
\special{dt 0.045}%
\special{pa 3091 2609}%
\special{pa 2107 1034}%
\special{dt 0.045}%
% VECTOR 0 0 3 0
% 2 3740 680 3740 1051
% 
\special{pn 20}%
\special{pa 3682 670}%
\special{pa 3682 1035}%
\special{fp}%
\special{sh 1}%
\special{pa 3682 1035}%
\special{pa 3701 969}%
\special{pa 3682 983}%
\special{pa 3662 969}%
\special{pa 3682 1035}%
\special{fp}%
% STR 2 0 3 0
% 3 3740 480 3740 580 5 0
% $\IP(\aleph_2,\aleph_2)$
\put(9.3500,-1.4500){\makebox(0,0){$\IP(\aleph_2,\aleph_2)$}}%
% STR 2 0 3 0
% 3 3740 -40 3740 60 5 0
% $\IP(\aleph_2,\aleph_1)$
\put(9.3500,-0.1500){\makebox(0,0){$\IP(\aleph_2,\aleph_1)$}}%
% VECTOR 0 0 3 0
% 2 3740 180 3740 441
% 
\special{pn 20}%
\special{pa 3682 178}%
\special{pa 3682 435}%
\special{fp}%
\special{sh 1}%
\special{pa 3682 435}%
\special{pa 3701 369}%
\special{pa 3682 382}%
\special{pa 3662 369}%
\special{pa 3682 435}%
\special{fp}%
% POLYLINE 1 1 3 0
% 6 5450 3670 5450 1670 4500 1270 4500 310 1350 310 1350 310
% 
\special{pn 13}%
\special{pa 5365 3613}%
\special{pa 5365 1644}%
\special{pa 4430 1251}%
\special{pa 4430 306}%
\special{pa 1329 306}%
\special{pa 1329 306}%
\special{da 0.050}%
% STR 2 0 3 0
% 3 1150 210 1150 310 5 0
% (8)
\put(2.8750,-0.7750){\makebox(0,0){(8)}}%
\end{picture}%
\hfill\\
\mbox{}\hfill {\bf fig.\ 9}\hfill\mbox{}\bigskip\bigskip\\
Finally, we shall mention some open problems. 

In \cite{fu-ge-so} it is shown that $\gma=\aleph_1$ follows from \WFN\  
where $\gma$ is the 
almost disjoint number. 
In \cite{fuchino-geschke}, it is then shown that, under some 
additional assumptions, $\gma=\aleph_1$ already follows 
from \SEP\ which is a weakening of \WFN. Therefore, it seems natural to ask 
the following 
question:
\begin{Problem}
	Does $\gma=\aleph_1$ follow from $\HP(\aleph_2)$ or
	$\IP(\aleph_2,\lambda)$ for $\lambda=\aleph_1$,$\aleph_2$\,? 
\end{Problem}
\begin{Problem}
	Does $\WFN$ imply $\HP(\aleph_2)$ or $\donum=\aleph_1$ \,?
\end{Problem}
The model of $\boundingno^*=\aleph_2$ and $\donum=\aleph_1$ satisfies a 
strong form of negation of $\Cs(\aleph_2)$. This suggests the following problem:
\begin{Problem}
	Does $\HP(\aleph_2)$ (or even $\Cs(\aleph_2)$) 
imply $\boundingno^*=\aleph_1$\,?
\end{Problem}

At the moment, we do not have any model separating $\HP(\kappa)$ and
$\IP(\kappa,\kappa)$ for $\kappa>\kappa_1$. 
\begin{Problem}
	Is $\HP(\kappa)$ $+$ $\neg\IP(\kappa,\kappa)$ consistent for some (or 
	any) $\kappa>\aleph_1$\,?
\end{Problem}
%% \begin{Problem}
%% 	Does $\IP(\aleph_2,\aleph_1)$ hold in the side-by-side Sacks model? 
%% 	(cf.\  \Corof{fs-prod-IP}, \assertof{e}). 
%% \end{Problem}

In \Corof{chang-conj-then-hechler} which realizes the separation (7) in 
fig.9, a very strong large cardinal property is assumed.
\begin{Problem}
	Can we construct a model realizing (7) in fig.9 starting from \ZFC\ 
	without any large cardinal? 
\end{Problem}
The property (A) used in the proof of \Thmof{new-model} and proved to be 
consistent with \CH\ in \sectionof{forcing(A)} seems to be of its own interest.
\begin{Problem}
	Is $\neg${\rm(A)} consistent with \ZFC\ $+$ \CH\ (or even with \ZFC\ $+$ \GCH)\,?
\end{Problem}

\parbox[t]{10cm}{
\noindent
{\bf Authors' addresses}\bigskip\medskip\\
{\rm J\"org Brendle}\smallskip\\
Graduate School of Science and Technology\\ 
Kobe University\\ 
Rokko-dai 1-1, \\
Nada KOBE 657-8501 
Japan \\
{\tt  brendle@kurt.scitec.kobe-u.ac.jp}
\medskip\\
{\rm Saka\'e Fuchino}
\smallskip\\
Dept.\ of Natural Science and Mathematics\\
  College of Engineering,
  Chubu University,\\
  Kasugai AICHI 487-8501 Japan.\\
{\tt fuchino@isc.chubu.ac.jp}} 


\begin{thebibliography}{xxx}
\newcommand{\bysame}{\underline{\ \ \ \ \ \ }}
%% \bibitem{brendle-fuchino} J.\ Brendle and S.\ Fuchino, Coloring ordinals 
%% 	by reals, preprint.
%------------------------
\bibitem{brendle-laberge} 
%------------------------
	J.\ Brendle and T.\ LaBerge, 
	Forcing tightness in products of fans, 
	Fundamenta Mathematicae, 150,  211-226 (1996).
%--------
\bibitem{CP}
%--------
K.\ Ciesielski and J.\ Pawlikowski, 
The Covering Property Axiom, CPA, Cambridge University Press (2004).
%-----------------------
\bibitem{eda-kada-yuasa} 
%-----------------------
	K.\ Eda, M.\ Kada and Y.\ Yuasa, 
	The tightness about sequential fans and combinatorial properties, 
	Journal of Mathematical Society of Japan 49, 181--197 (1997).
%% %---------------------
%% \bibitem{fuchino-LIMS} 
%% %---------------------
%% 	S.\ Fuchino, 
%% 	On combinatorial principles
%% 	$\mbox{\sc Princ}(\kappa,\lambda)$, 
%%   $\mbox{C}^s(\kappa)$, $\mbox{HP}(\kappa)$ etc., 
%% 	RIMS Kokyuroku No.\ 1423, 13-27 (2005).
%------------------------
\bibitem{fuchino-geschke} 
%------------------------
	S.\ Fuchino and S.\ Geschke, 
	Some combinatorial principles defined in terms of elementary submodels, 
	Fundamenta Mathematicae 181, 233-255 (2004).
%-----------------
\bibitem{fu-ge-so} 
%-----------------
	S.\ Fuchino, S.\ Geschke and L.\ Soukup, 
	On the weak Freese-Nation property of ${\cal P}(\omega)$, 
	Archive for Mathematical Logic, Vol.40, 425-435 (2001). 
%------------------
\bibitem{fuchino-x}
%------------------
	S.\ Fuchino, S.\ Geschke and L.\ Soukup, 
	Principles 	capturing features of generic extensions by 
	almost side-by-side product, in preparation. 
%--------------------
\bibitem{fu-ge-sh-so} 
%--------------------
	S.\ Fuchino, S.\ Geschke, S.\ Shelah and L.\ Soukup, 
	On the weak Freese-Nation property of complete Boolean algebras,
	Annals of Pure and Applied Logic, 110 (1-3) (2001) 89-105.
%-----------------
\bibitem{fu-ko-sh} 
%-----------------
	S.\ Fuchino, S.\ Koppelberg and S.\ Shelah, Partial 
	orderings with the weak Freese-Nation property, Annals of Pure and 
	Applied Logic 80, 35--54 (1996).  
%% \bibitem{chubu} S.\ Fuchino and L.\ Soukup, Chubu-2001 notes.
%---------------
\bibitem{FuShSo}
%---------------
	S.\ Fuchino, S.\ Shelah and L.\ Soukup, Sticks and clubs, 
Annals of Pure and Applied Logic 90, no.1 (1997), 57-77. 
%--------------
\bibitem{fu-so} 
%--------------
	S.\ Fuchino and S.\ Soukup, More set-theory around the weak 
	Freese-Nation property, Fundamenta Matematicae 154, 159--176 (1997). 
%----------------
\bibitem{geschke} 
%----------------
	S.\ Geschke, On $\sigma$-Filtered Boolean Algebras, Dissertation, 
	Freie Universit\"at Berlin (1999). 
%---------------------
\bibitem{juhasz-kunen}
%---------------------
	I.~Juh\'asz and K.~Kunen, The Power Set of
	$\omega$, Elementary submodels and
	weakenings of CH, Fundamenta Mathematicae 170, 257--265 (2001).
%-----------------
\bibitem{ju-so-sz} 
%-----------------
	I.\ Juh\'asz, L.\ Soukup and Z.\ Szentmikl\'ossy, 
	Combinatorial principles from adding Cohen reals, 
	Proceedings of Logic Colloquium (1995). 
%-------------------
\bibitem{kada-yuasa}
%-------------------
	M. Kada and Y. Yuasa, 
	Cardinal invariants about shrinkability of unbounded sets, 
  Topology and its Applications 74 (1996), 215--223.
%-------------------
\bibitem{kunen-book}
%-------------------
	K.~Kunen: {Set Theory}, North-Holland (1980).
%---------------
\bibitem{miller}
%---------------
A.\ Miller, Infinite Combinatorics and Definability, 
Annals of Pure and Applied Mathematical Logic,
41 (1989), 179-203.
%---------------
\bibitem{shelah} 
%---------------
	S.\ Shelah,  
	a(n unpublished?) note for I.\ Juh\'asz  (2002).
%% %--------------
%% \bibitem{yuasa}
%% %--------------
%% Y.\ Yuasa, Shrinkability of unbounded sets in the Cohen extension, 
%% Mathematical Logic and Applications `94, Research Institute for 
%% Mathematical Sciences 930, Kyoto University (1995), 55--58.
%% %---------------------
%% \bibitem{yuasa-thesis}
%% %---------------------
%% Y.\ Yuasa, Infinitary combinatorics in the Cohen extension, 
%% Thesis, Waseda University (1996).
\end{thebibliography}
\end{document}